%% file: paper_tensors.tex
\documentclass[a4paper, notitlepage, 12pt]{article}
\usepackage{geometry}
\usepackage{ucs}                % Unicode support in LaTeX
\usepackage[utf8x]{inputenc}    % We use UTF-8 encoded files only

\usepackage[T1]{fontenc}     % for the right font encoding
\usepackage{lmodern}         % fonts for good pdf files

\usepackage[british]{babel}
\usepackage{lineno}
\usepackage{color}
\definecolor{mygreen}{rgb}{0.45, 0.65, 0.0}
\newcommand{\rvsi}[1]{\textcolor{black}{#1}}
\definecolor{myblue}{rgb}{0.1, 0.2, 0.9}
\newcommand{\rvsn}[1]{\textcolor{black}{#1}}

\usepackage{authblk}

\usepackage[table,x11names]{xcolor} %always add at the top

\usepackage{latexsym}
\usepackage{amsmath}
\usepackage{amssymb}

\usepackage{upgreek}
\newcounter{bla}

\newenvironment{keyword}{%
\noindent \emph{Keywords:}}{\newline}

\newcommand{\olsi}[1]{\,\overline{\!{#1}}} % overline short italic

\usepackage{graphicx}
\usepackage{calc}
\usepackage{multicol}
\usepackage{multirow}

\usepackage{float}
\usepackage[font=normalsize]{caption}
\usepackage{subcaption} % (avoid using outdated \subfig package which enables \subfloat)
\usepackage{mathtools}
\usepackage{siunitx}
\usepackage{booktabs}
\usepackage{rotating}
\usepackage{multirow}
\usepackage{wrapfig}
\usepackage[colorlinks]{hyperref}
\usepackage{etoolbox}
\usepackage{tikz}
\usepackage{url}
\usepackage{cleveref}
\usepackage[export]{adjustbox}
\usepackage{mathrsfs}
\usepackage[toc,page]{appendix}
\usepackage{pgfplots} 
\usepackage{pgfgantt}
\usepackage{pdflscape}
\pgfplotsset{compat=newest} 
\pgfplotsset{plot coordinates/math parser=false}

\input{macros}

\begin{document}
	
\sloppy % for better right margin text alignment
%\RaggedRight
%\raggedright

%\begin{frontmatter}

% =============================================================================
% Setup title, author, date and keywords.

% \thetitle, \theauthor, \thesubject and \thekeywords are used in pdfsetup.tex
%\newcommand{\authorcr}{\\}
\newcommand{\autheadcr}{\authorcr}

\newcommand{\authorsks}{Sharana Kumar Shivanand}
\newcommand{\authorbr}{Bojana Rosi\'{c}}
\newcommand{\authorhgm}{Hermann~G. Matthies}
\newcommand{\theauthor}{\authorsks \hspace{2.0em} \authorbr \hspace{2.0em} \authorhgm}

\newcommand{\affilwire}{{\small Institute of Scientific Computing, 
   Technische Universit\"at Braunschweig, Germany}}
\newcommand{\affiltwen}{{\small Applied Mechanics and Data Analysis, University of Twente, The Netherlands}}

% another title: Representation of random symmetric positive-definite material tensors
\newcommand{\thetitle}{Stochastic Modelling of\autheadcr Symmetric Positive Definite Material Tensors}

\newcommand{\textdate}{\today}

% ============================================================================
% connect to default LaTeX values
\title{\thetitle}

%\author{\theauthor}
\author[a]{\authorsks\thanks{Corresponding author. {E-mail:} \texttt{sshivanand@turing.ac.uk}}}
\author[b]{\authorbr}
\author[a]{\authorhgm}
%\affil{\affilwire}
\affil[a]{\affilwire}
\affil[b]{\affiltwen}

%\thanks{Corresponding author: D-38092 Braunschweig, 
%       Germany, e-mail: \texttt{wire at tu-bs.de}}

\date{\textdate}

\title{\thetitle}
%\author{\authorsks$^{\text{a},*}$}
%\author{\authorbr$^{\text{b}}$}
%\author{\authorhgm$^{\text{a}}$}
%\cortext[author] {Corresponding author. {E-mail:} \texttt{s.shivanand@tu-bs.de}}
%\address{$^{\text{a}}$ Institute of Scientific Computing, Technische Universit\"{a}t Braunschweig, Germany \\
%$^\text{b}$ Applied Mechanics and Data Analysis, University of Twente, The Netherlands}

\maketitle

%===========================ABSTRACT and KEYWORDS=====================================================

\input{abstract}

%=====================================================================================================

%\end{frontmatter}

%=======================================INTRODUCTION==================================================

\input{introduction}
\input{problem}

%=======================================MODELLING==================================================	

\input{modelling}

%========================================APPLICATION==================================================

% \input{application}
\input{application-y}

%==============================================RESULTS================================================ 

\input{results}

%==========================================CONCLUSION================================================

\input{conclusion}

%\newpage

\appendix

\input{abstract-metrics}

\vspace{2em}
\thanks{{\noindent \bf Acknowledgement:} Partly supported by the Deutsche Forschungsgemeinschaft 
(DFG) and a Gay-Lussac Humboldt prize.}

%==========================================BIBLIOGRAPHY===============================================

% \clearpage
%\bibliographystyle{abbrv}
\bibliographystyle{hgmplain-1}
%\bibliography{references}
\bibliography{paper_tensors}

\end{document}

%% file: macros.tex
\usepackage[sfdefault=cmbr,OMLmathsans]{isomath}  % isomath for slanted sans-serif

\newcommand{\ignore}[1]{}

\newcommand{\vek}[1]{\mathchoice{\displaystyle\boldsymbol#1}
{\textstyle\boldsymbol#1}{\scriptstyle\boldsymbol#1}
{\scriptscriptstyle\boldsymbol#1}}

\newcommand{\F}{\mathfrak} % math Fraktur
\newcommand{\C}{\mathcal}  % math Calligraphic
\newcommand{\mrm}{\mathrm}     % math roman upright

\newcommand{\di}{\mathop{}\!\mathrm{d}}

\newcommand{\ip}[2]{\langle #1 , #2 \rangle}

\newcommand{\nd}[1]{\| #1 \|}

\newcommand{\Ho}{\mathring{\mrm{H}}}
\newcommand{\Lp}{\mrm{L}}

\newcommand{\Hp}{\mrm{H}}

\DeclareMathOperator{\tr}{tr}
\DeclareMathOperator{\diag}{diag}

\newcommand{\meansymb}[1]{\widehat{#1}}
\newcommand{\wsymb}[1]{\mathring{#1}}	
\newcommand{\lsymb}[1]{\widetilde{#1}}	
\newcommand{\tsymb}[1]{\olsi{#1}}

\usepackage[nopostdot,nonumberlist,acronym,toc,section]{glossaries}
\glsdisablehyper %used to disable hyperef for acronyms
\newacronym{bmd}{BMD}{bone mineral density}
\newacronym{qct}{QCT}{quantitative computed tomography}
\newacronym{ct}{CT}{computed tomography}
\newacronym{hu}{HU}{Hounsfield units}
\newacronym{dxa}{DXA}{dual energy X-ray absorptiometry}
\newacronym{fea}{FEA}{finite element analysis}
\newacronym{fem}{FEM}{finite element method}
\newacronym{mc}{MC}{Monte Carlo method}
\newacronym{mpp}{MPP}{most probable point}
\newacronym{amv}{AMV}{advanced mean-value method}
\newacronym{form}{FORM}{first order reliability method}
\newacronym{sorm}{SORM}{second order reliability method}
\newacronym{rsms}{RSMs}{response surface methods}
\newacronym{pce}{PCE}{polynomial chaos expansion}
\newacronym{bc}{BC}{boundary condition}
\newacronym{pdf}{PDF}{probability density function}
\newacronym{mse}{MSE}{mean square error}
\newacronym{rmse}{RMSE}{root mean square error}
\newacronym{lhs}{LHS}{Latin hypercube sampling}
\newacronym{ism}{ISM}{importance sampling method}
\newacronym{kle}{KLE}{Karhunen-Lo\`{e}ve expansion}
\newacronym{spde}{SPDE}{stochastic partial differential equation}
\newacronym{pde}{PDE}{partial differential equation}
\newacronym{mlmc}{MLMC}{multilevel Monte Carlo method}
\newacronym{dof}{DOF}{degrees of freedom}
\newacronym{iso}{iso}{isotropy}
\newacronym{ortho}{ortho}{orthotropy}
\newacronym{scl}{scl}{scaling}
\newacronym{dir}{dir}{direction}
\newacronym{nt}{NT}{nodal temperature}
\newacronym{thfl}{THFL}{total heat flux}
\newacronym{nhfl}{NHFL}{normalized heat flux}

\newcommand{\spatialdomain}{\mathcal{G}}

\newcommand{\Euclideanspace}{\mathbb{R}^d}

\newcommand{\dispvector}{\vek{u}}

\newcommand{\conducttensor}{\vek{\kappa}}

\newcommand{\solspace}{\mathcal{U}}
\newcommand{\finitespace}{\mathcal{U}_h}
\newcommand{\solspaceheat}{\mathcal{V}}

\newcommand{\samplespace}{\varOmega}
\newcommand{\event}{\omega}

\newcommand{\sigmaalgebra}{\mathfrak{F}}
\newcommand{\probability}{\mathbb{P}}

\newcommand{\DFcond}{\conducttensor(x)}

\newcommand{\samplesize}{N}
\newcommand{\Expvalue}{\mathbb{E}}

\newcommand{\normaldist}{\mathcal{N}}

\newcommand{\SymMatSpace}{\mrm{Sym}}
\newcommand{\PosSymMatSpace}{\SymMatSpace^+}
\newcommand{\SymSecRankSpace}{\SymMatSpace(d)}
\newcommand{\PosSymSecRankSpace}{\PosSymMatSpace(d)}
\newcommand{\PosSymSecRankSpacetwoD}{\PosSymMatSpace(2)}
\newcommand{\PosSymSecRankSpacethreeD}{\PosSymMatSpace(3)}
\newcommand{\meanconducttensorscl}{\meansymb{\vek{\kappa}}^{iso}_{s}}
\newcommand{\meanconducttensorortho}{\meansymb{\vek{\kappa}}^{ortho}_{r}}
\newcommand{\conducttensorscl}{\vek{\kappa}_{s}}
\newcommand{\conducttensorrot}{\vek{\kappa}_r}
\newcommand{\genSecTensor}{\vek{C}}
\newcommand{\symtensor}{\vek{H}}
\newcommand{\symdiag}{\vek{Y}}
\newcommand{\SecTensorB}{\vek{B}}

\newcommand{\meanY}{\meansymb{\symdiag}}

\newcommand{\meanH}{\meansymb{\symtensor}}

\newcommand{\wmeanH}{\wsymb{\symtensor}}
\newcommand{\meanC}{\meansymb{\genSecTensor}}
\newcommand{\tmeanC}{\tsymb{\genSecTensor}}
\newcommand{\wmeanC}{\wsymb{\genSecTensor}}
\newcommand{\lmeanC}{\lsymb{\genSecTensor}}

\newcommand{\rotmat}{\vek{R}}
\newcommand{\rotangle}{\phi}

\newcommand{\diffoperator}{\mathcal{A}}
\newcommand{\Real}{\mathbb{R}}
\newcommand{\PosRealStrict}{\mathbb{R}^{+}}

\newcommand{\eigvect}{\vek{Q}}
\newcommand{\eigval}{\vek{\Lambda}}

\newcommand{\eigvalsclele}{{{\lambda}}_s^{(i)}(\event_s)}
\newcommand{\meaneigvect}{\meansymb{\eigvect}}
\newcommand{\meaneigval}{\meansymb{\eigval}}

\newcommand{\eigvalscleleone}{{{\lambda}}^{(1)}_s}
\newcommand{\eigvalscleletwo}{{{\lambda}}^{(2)}_s}
\newcommand{\eigvalsclelethree}{{{\lambda}}^{(3)}_s}

\newcommand{\meaneigvalscleleone}{{\meansymb{{\lambda}}}^{(1)}_s}
\newcommand{\meaneigvalscleletwo}{{\meansymb{{\lambda}}}^{(2)}_s}
\newcommand{\meaneigvalsclelethree}{{\meansymb{{\lambda}}}^{(3)}_s}
\newcommand{\meaneigvalroteleone}{{\meansymb{{\lambda}}}^{(1)}_r}
\newcommand{\meaneigvalroteletwo}{{\meansymb{{\lambda}}}^{(2)}_r}
\newcommand{\meaneigvalrotelethree}{{\olsi{{\lambda}}}^{(3)}_r}

\newcommand{\concparam}{\eta}
\newcommand{\unitcircle}{\mathbb{S}^1}
\newcommand{\circmeananly}{\olsi{\rotangle}}
\newcommand{\rotvect}{\vek{v}}

\newcommand{\temperature}{T}
\newcommand{\temperatureApprox}{{T}_h}
\newcommand{\stotemp}{{T}_{h,\samplesize}}
\newcommand{\meanMC}{\mu^{(\mrm{MC})}}
\newcommand{\meanMCbold}{\vek{\mu}^{(\mrm{MC})}}

\newcommand{\heatflux}{\vek{q}}

\newcommand{\stototalflux}{q^{(t)}_{h,\samplesize}}
\newcommand{\stofluxunit}{{\hat{\vek{q}}}_{h,\samplesize}}

\newcommand{\threeD}{\mathbb{R}^3}
\newcommand{\axis}{\vek{w}}
\newcommand{\skewsym}{\vek{W}}
\newcommand{\unitsphere}{\mathbb{S}^{d-1}}
\newcommand{\std}{\mathbb{S}\mrm{td}}
\newcommand{\stdMC}{\std^{(\mrm{MC})}}
\newcommand{\stdMCcirc}{\stdMC_{\mrm{circ}}}
\newcommand{\mucirc}{\vek{\mu}_{\mrm{circ}}}
\newcommand{\muMCcirc}{\mucirc^{(\mrm{MC})}}

\newcommand{\feq}[1]{Eq.(\ref{#1})} %[1]{#1} : [no of input parameters]{Parameter-1}
\newcommand{\feqs}[2]{Eqs.(\ref{#1}) and (\ref{#2})} %\ffigs{arg1}{arg2}
\newcommand{\feeqs}[2]{Eqs.(\ref{#1})--(\ref{#2})} %\ffigs{arg1}{arg2}
\newcommand{\fsec}[1]{Section~\ref{#1}}
\newcommand{\fsecs}[2]{Sections~\ref{#1} and \ref{#2}} %\ffigs{arg1}{arg2}

\newcommand{\ffig}[1]{Fig.\ref{#1}} %\ffig{arg1}
\newcommand{\ffigs}[2]{Figs.\ref{#1} and \ref{#2}} %\ffigs{arg1}{arg2}
\newcommand{\ffiigs}[2]{Figs.\ref{#1}--\ref{#2}} %\ffigs{arg1}{arg2}
\newcommand{\ftbl}[1]{Table~\ref{#1}}

%% file: abstract.tex
% !TEX root = ../paper_tensors.tex
% !TEX encoding = UTF-8 Unicode

\begin{abstract}
%\textit{
%In this paper, we model and generate random material tensors that are symmetric and positive definite (SPD). In particular, we limit ourselves to spatially constant second-order SPD tensor fields. The tensors are modelled such to have a fine independent control over its scaling and directional aspects, given that the material symmetry is fixed. Furthermore, another scenario of transitioning from the expected higher symmetry class in the mean to subsequent lower symmetry class in the realisations is also explored.
%As an example, we model random thermal conductivity tensor of a steady-state heat conduction of 2D and 3D human proximal femur.
%The numerical results showcase the distinct impact of incorporating different material uncertainties---scaling, orientation and material symmetry---independently into the constitutive model on the desired quantities of interest, such as nodal temperature and heat flux.
%}
Spatial symmetries and invariances play an important role in the behaviour of materials and
should be respected in the description and modelling of material properties.
The focus here is the class of physically symmetric and positive definite tensors,
as they appear often in the description of materials, and
%When modelling random material properties, 
%Here we discuss how to model and generate random ensembles of tensors where
one wants to be able to prescribe certain classes of spatial symmetries and invariances for
each member of the whole ensemble, while at the same time demanding that the mean or expected value
of the ensemble be subject to a possibly `higher' spatial invariance class.
We formulate a modelling framework which not only respects 
these two requirements---positive definiteness and invariance---but also allows
a fine control over orientation on one hand, and strength / size on the other.
As the set of positive definite tensors is not a linear space, but rather an open convex cone
in the linear space of physically symmetric tensors, we consider it advantageous
to widen the notion of mean to the so-called Fr\'{e}chet mean on a metric space,
which is based on distance measures or metrics between positive definite tensors other than
the usual Euclidean one. %, which is the cause of some undesirable 
%  Here the open convex cone is considered as a Riemannian manifold,
%which then yields a geodetic distance to be used in the definition of the Karcher or 
%Fr\'{e}chet mean; 
%We formulate certain desiderata for such means and metrics and discuss different possibilities.
%For the sake of simplicity, as well as to expose
%the main idea as clearly as possible, we limit ourselves here to
%second order tensors.  
It is shown how the random ensemble can be modelled and generated,
%with a fine control of the spatial symmetry or invariance of the whole ensemble, as
%well as its Fr\'{e}chet mean, 
independently in its scaling and orientational or directional aspects,
with a Lie algebra representation via a memoryless transformation.
The parameters which describe the elements in this Lie algebra are then to be considered
as random fields on the domain of interest.
As an example, a 2D and a 3D model of steady-state heat conduction in a human proximal femur, a
bone with high material anisotropy, is modelled with a random thermal conductivity tensor,
and the numerical results show the distinct impact of incorporating into the constitutive model 
different material uncertainties---scaling, orientation, and prescribed material symmetry---on 
the desired quantities of interest. %, such as temperature distribution and heat flux.
\end{abstract}

%======================================KEYWORDS=========================================================


%% file: introduction.tex
% !TEX root = ../paper_tensors.tex
% !TEX encoding = UTF-8 Unicode
\section{Introduction}\label{Introduction}
Heterogeneous materials and their statistical description are of great interest in many
fields.  The properties of such materials often vary considerably on the spatial scales
of interest, and additionally in many cases the material properties are only described
statistically, reflecting some underlying uncertainty as to the exact values.  
This leads to the idea of a stochastic/random representation of these properties.
The material properties are naturally collected in tensor quantities, 
%due to Curie's principle on spatial symmetries or invariances, and
%due to Onsager's reciprocity relations \cite{deGrootMazur1984}.  
and these tensors are often 
of even order and physically symmetric---in elasticity this is often termed the \emph{major}
symmetry.  Furthermore, when these tensorial properties appear in
the definitions of stored energy or entropy production of stable systems,
they are additionally positive definite.  It is this class of even-order physically
symmetric and positive definite (SPD) tensors which we are interested in describing.
Well known examples of such tensors are the second-order tensors describing
diffusion phenomena, mapping the gradient of some quantity to the flux of some related
quantity, e.g.\ the thermal conductivity tensor mapping temperature gradients to
fluxes of thermal energy.  An example of a fourth-order tensor
is the elasticity tensor of a linearly elastic material, which maps
the strain tensor to the stress tensor.  Any even-order
tensor may be viewed as a linear mapping on the space of tensors of half the order
by contracting over half the indices.
By choosing a basis in the space of tensors of half the order, each such even-order tensor can
hence be represented by a matrix; so the class we are interested in can
thus be represented by SPD matrices.  We shall switch to
this matrix representation whenever it is convenient to do so.

Even though frequently the exact values of the tensorial entries may not be known,
some properties resulting from general principles often are known; namely,
often it is known what kind of symmetries or invariances are to be expected
\cite{nye_physical_1984, Cowin2013, malgrange_symmetry_2014}.  These invariances---we
prefer the term invariance over symmetry here to avoid confusion
with the afore mentioned physical symmetry, resp.\ the symmetry of the 
tensor as a linear mapping---under the operation
of some symmetry group are well known, and here we concentrate on spatial point operations
such as those defining isotropy.  These symmetry groups are represented as linear transformation
groups on the relevant space of even order tensors, or, equivalently, on the space of representing
symmetric matrices mentioned above.  Such transformation groups then define invariant subspaces
which represent the tensors with the given invariance property.
As the number of such point groups is larger in the case of fourth-order spatial tensors 
(e.g.\ \cite{cowin_identification_1987, bona_coordinate-free_2007, Cowin2013})
than in the case of second-order spatial tensors \cite{malyarenko_tensor-valued_2019},
in our explicit treatment here we shall confine ourselves to second-order spatial tensors
for the sake of simplicity and clarity.  
But the general idea applies to any even-order SPD tensor field.

So, after a choice of bases, what will be considered here are 
random SPD matrices, where certain invariances are known for the whole population,
as well as for the mean.
There is already quite a bit of previous work in this field, drawing on several sources,
which will be reviewed only very cursory and briefly here.
Random SPD matrices (e.g.\ \cite{gupta_matrix_1999}) arise in a number of fields,
e.g.\ \cite{schwartzman_random_nodate, soize_nonparametric_2000, GuilleminotSoizeGhanemR.2011,
Grigoriu2016}, to name just \rvsi{some}.

\rvsi{Each SPD matrix can be spectrally decomposed with positive eigenvalues
and orthogonal eigenvectors.  This separates size / strength aspects as encoded 
in the eigenvalues, and orientation as encoded in the eigenvectors.}  
We see these two items as vital information about
how the tensor / matrix acts, and we shall base our treatment on this 
fundamental decomposition.
On the other hand, for the purpose of modelling or representing a tensor, the most
comfortable and desirable is to find a representation on a linear space.  It
is one of the aims to show here how this may be achieved for the spectral decomposition.
\rvsi{Geometrically speaking, the}
positive definite symmetric even-order tensors resp.\ SPD matrices are an open convex
cone in the vector space of symmetric even-order tensors resp.\ symmetric matrices,
and are thus a differentiable manifold, but not a vector space. 
\rvsi{The spectral decomposition with positive eigenvalues represents this cone,
and it will be shown that this decomposition is connected to another differentiable
manifold, which is even a Lie group---a manifold with a differentiable group structure.  
This is one of the uses of group theory here, and in the actual 
numerical modelling it just boils down to the use of the matrix exponential.}

Differentiable manifolds can be represented---mostly only locally---on their 
tangent spaces.  The analysis of manifolds is 
typically facilitated if they can be adorned with a Riemannian structure
\cite{PennecFillardAyache2004, ThanwerdasPennec2019, ThanwerdasPennec2021},
and even more so if the manifold carries \rvsi{the structure of a Lie group}, 
as then everything can be played down to
the tangent space at the neutral element, i.e.\ the associated Lie algebra.
\rvsi{The Lie algebra is mapped to the Lie group via the exponential map.}  
Unfortunately, the manifold of SPD matrices is not a Lie group under
normal matrix multiplication, although it is possible \cite{ArsignyFillardPennecEtAl2007}
to give the SPD manifold a new commutative or Abelian multiplicative group structure, 
which coincides with matrix multiplication on commuting matrices.
%One may use \rvsi{the map from the Lie algebra to the group} to also
%carry any Euclidean structure on the Lie algebra to a Riemannian one on the Lie group.
%\rvsi{Our approach here is similar in the use of the exponential map, but has a 
%notable difference, partly following \cite{schwartzman_random_nodate, Grigoriu2016},
%where SPD matrices are represented through their spectral decomposition
%on the product of the group of orthogonal and
%the Abelian group of diagonal SPD matrices.}

\rvsi{The second use of group theory is in the material symmetry induced group reduction
\cite{faessStiefel1992} ---introducing as many zeros in the tensor as possible 
through spatial rotations, or, equivalently, through reflections \cite{Cowin2013}. 
This well-known subject defines the different symmetry classes
\cite{Cowin2013, malyarenko_tensor-valued_2019}.}
For second-order spatial tensors this spatial group reduction %eigenvector decomposition 
completely determines the eigenvector transformation of the associated SPD matrix
\rvsi{into classes defined by the number of degenerate or coinciding eigenvalues}, whereas for
e.g.\ for fourth-order tensors there is an iso-spectral sub-group in the 
relevant Lie group of orthogonal transformations, cf.\ e.g.\
\cite{silvestrov_spectral_2016, malyarenko_random_2017, altenbach_tensor_2018, 
malyarenko_tensor-valued_2019}, which determines the eigen-strain distribution 
in the relevant elasticity class
\cite{nye_physical_1984, cowin_identification_1987, cowMehr1995, 
malgrange_symmetry_2014, bona_coordinate-free_2007}.
This is one more reason to confine ourselves here initially to second-order
spatial tensors for the sake of brevity and simplicity of exposition,
the corresponding more involved results for fourth-order spatial tensors 
(i.e.\ elasticity tensors) will be published elsewhere.

The spectral decomposition of a random SPD matrix has thus two components,
the SPD diagonal matrix of eigenvalues and the orthogonal orientation matrix.
In particular---as will be shown---this allows the independent control of the
strength (or size or norm) of the tensor and of its directional properties.
Both components may be random, and thus random orthogonal matrices have some
relevance here as well 
\cite{mardia_directional_2000, jammalamadaka_topics_2001, Mezzadri2007, Grigoriu2016}.
%It is well known that the manifold of orthogonal matrices is a non-commutative Lie group
%with the associated Lie algebra the space of skew-symmetric matrices.
%But the SPD diagonal matrices---i.e. diagonal matrices with positive entries
%on the diagonal---are also a commutative or Abelian Lie group under matrix
%multiplication with the associated Lie algebra the diagonal matrices.  
\rvsi{The manifold of SPD matrices may thus be modelled on the product of the two Lie groups 
\cite{JungSchwartzmanGroisser2015, GroisserJungSchwartzman2017, groissJungSchwman2017}
of orthogonal and of positive diagonal matrices, and ultimately on their
Lie algebras---the vector space of skew-symmetric and the space of diagonal matrices.
A somewhat similar approach is also explored in \cite{Grigoriu2016}, where rotation
angles are used for the representation of rotations, which is one way to parametrise
the Lie algebra of skew matrices.  Here the well-known Rodrigues formula will be used
for that purpose.}
%Thus the spectral decomposition leads to a representation on a product
%of Lie algebras---the orthogonal matrices and the positive definite diagonal ones.

\rvsi{The Euclidean structure on the direct sum of the Lie algebras just mentioned
may then} be mapped \rvsi{via the exponential map} onto a product Riemannian
structure on the product manifold of the two Lie groups which represents the
SPD matrices.  \rvsi{Such} a Riemannian structure on a 
manifold as just alluded to allows one to define
the length of paths, and thus shortest paths resp.\ geodesics, and hence to 
measure distances on the manifold.  As the variational definition of mean, 
resp.\ average, or expectation---we will use these terms interchangeably---depends 
on the distance measure (e.g.\ \cite{NielsenBhatia2013, GuiguiEtal2023}, cf.\ also 
\fsec{Sec:Problem}), this leads to the question as to what is an adequate definition 
of the mean for SPD tensors.  
%This is connected on how the set of SPD tensors---which may be
%seen as an open convex cone in the linear space of all symmetric tensors---is to be regarded.
%If viewed as a general metric space or more specifically as a Riemannian manifold, other possibilities
%arise which allow to consider further \emph{intrinsic} properties of SPD tensors.
%As regards the mean or expected value, questions arise
%about what is the \emph{proper} or \emph{adequate} mean to use.  
%This is due to the fact that the 
The usual arithmetic mean is wedded to a flat (vector-)space and a Euclidean distance, as
is well known and will be explained in more detail later in \fsec{Ssec:means-spd}.
On a metric space the mean of two points is the point half-way along the shortest path between 
them.  This generalisation is known as the Fréchet mean.  
This connection between the mean
(and variance to be defined later) and a distance metric may be used to switch the
discussion of what is desired from the mean to the underlying distance metric,
as on a non-flat manifold the Euclidean metric of the ambient space 
is not advantageous.
For example on the Lie group of orthogonal matrices it is completely common 
to use such a Riemannian metric
\cite{mardia_directional_2000, jammalamadaka_topics_2001, Mezzadri2007}.

%As already mentioned, the SPD matrices, when viewed as a subset
%of the ambient linear vector space of symmetric matrices, are geometrically an open cone.
%And from this embedding obviously 
On the manifold of SPD matrices the usual Euclidean metric of the ambient space of symmetric
matrices is derived from the Frobenius resp.\ Hilbert-Schmidt inner product and norm.
It has long been recognised that this metric has some undesirable properties
\cite{PennecFillardAyache2004, AndoLiMathias2004, Moakher2005, ArsignyFillardPennecEtAl2006, 
ArsignyFillardPennecEtAl2007, drydenEtal2009, DrydenEtal2010, FujiiSeo2015, FeragenFuster2017}, in
particular the so-called \emph{swelling, fattening, and shrinking effects} in interpolation
resp.\ averaging \cite{schwartzman_random_nodate, JungSchwartzmanGroisser2015, 
schwartzman_lognormal_2016, GroisserJungSchwartzman2017, groissJungSchwman2017}.
This refers to the ellipsoid representing the SPD matrix, and is connected to the 
non-monotonic interpolation of the eigenvalues or some of their 
functions (e.g.\ the determinant) along Euclidean geodesics, i.e.\ straight lines.
In the context of stochastic material modelling this means a loss of anisotropy
for the whole ensemble, which may be undesirable in some situations.
In the literature just cited, there is an intensive discussion---particularly 
coming from the medical field of diffusion tensor imaging, where these undesirable
effects are not acceptable---on how to make the 
SPD tensors into a Riemannian manifold with 
different metrics \cite{PennecSommerFletcher2020}.  From this one obtains variatonally
the Fréchet or Karcher mean---the generalisation of the arithmetic resp.\ Euclidean 
mean to general metric spaces \cite{NielsenBhatia2013, GuiguiEtal2023, groissJungSchwman2023}.
This is an area of active research \cite{ThanwerdasPennec2019, ThanwerdasPennec2019-2, 
ThanwerdasPennec2021}, the `best' metric and mean seem to be application dependent,
and here we shall formulate some desiderata for a metric---and thus
mean---for SPD tensors of material properties, and propose one which in our view is most suited
for this kind of field.  \rvsi{Subsequently, we use this proposed Fréchet mean,
and it agrees nicely with the modelling in terms of a product Lie group, see also
\cite{JungRooksEtAl2023}.}

As to the modelling of random tensor fields with given invariance properties,
the culmination of a long series of works \cite{malyarenko_random_2017, 
ostoja-starzewski_microstructural_2007, silvestrov_spectral_2016, malyarenko_statistically_2014, altenbach_tensor_2018} is reported in \cite{malyarenko_tensor-valued_2019},
which looks at stochastically homogeneous and isotropic random tensor fields of any order 
in a three dimensional domain.  The results reported there allow a fine
control about what invariance the mean of such a tensor field has, and what kind
of invariance may be required of each realisation.  Using the spectral theorem for
the covariance of such homogeneous and isotropic fields, as well as the representation
theory of the appropriate groups, here one may find the spectral resolution of the
covariance function and from it a Karhunen-Loève like representation 
(e.g.\ \cite{hgm07}) of the random tensor field with the desired invariance properties.  
In these deliberations, no special attention is payed to 
the topic of positive definiteness, as the methods used
are purely vector space like linear combinations or generalisations such as series and integrals.
And although the limiting expression may well represent a positive definite tensor,
in practical computational situations these series and integrals have to be
truncated resp.\ numerically approximated, and this step may lead to a numerical loss
of positive definiteness.

A reduced non-parametric approach to generate random SPD matrices, which takes explicit account 
of the topic of positive definiteness, was presented in \cite{soize_nonparametric_2000, soize_2006}.  
Here an algebraic property, namely that positive
elements in the algebra are squares of other elements in the algebra, is used to ascertain
that each generated tensor field is indeed SPD.
%, and the theory of C$^*$-algebras 
%\cite{segalKunze78, yosida-fa-1980} ensures that each SPD tensor(-field) has a square root.
\rvsi{This approach to ensure positivity is also examined in \cite{Grigoriu2016}.} 
Based on this, the generation of
random elasticity tensors with specified invariance in the mean and fully anisotropic invariance
in each realisation, resp.\ controlled elasticity class invariance with constant spatial 
orientation of the symmetry axes in each realisation, is shown in a series of papers
\cite{GuilleminotSoize2011, GuilleminotSoize2012,  GuilleminotSoize2012a, GuilleminotSoize2013, 
GuilleminotSoize2013a, GuilleminotSoize2014, NouySoize2014, StaberGuilleminot2017, Soize2021a}.  

Another example for a frequent approach to ensure the SPD property
%, based on spectral and C$^*$-algebra properties, 
is to use the exponential function, resp.\ the logarithm
in the opposite direction.  
This is employed in e.g.\ \cite{GuilleminotSoize2013a, NouySoize2014, 
JungSchwartzmanGroisser2015, schwartzman_lognormal_2016,
GroisserJungSchwartzman2017, groissJungSchwman2017}, 
and it can even be combined with the previous squaring
approach \cite{GuilleminotSoize2013a, NouySoize2014, Grigoriu2016} to ensure compliance with
a particular invariance class.  This technique is used quite
widely  when describing measured diffusion and conductivity tensors,
namely instead of addressing the SPD tensor itself, one focuses with the description
on its logarithm.
%---and C$^*$-algebra theory assures us that a SPD tensor has a unique main-branch logarithm.  
To obtain the tensor, one finally takes the exponential, thus
ensuring positive definiteness.  This approach can take advantage of the fact that the
invariance classes are the same for a tensor and its exponential. 
We shall follow a similar path here and essentially use the exponential function
to ensure positive definiteness numerically in each realisation.

%random tensor field with same symmetry 
%useful in higher order tensors with lower material symmetries like elasticity where the identification of stochastic parameters becomes simpler \\
%cons of reduced approach: lack of control over different kinds of uncertainties \\
%classical Full parametric approach provides full control\\

%To recap, what is desired is a method to model and generate random SPD matrices, and that have
%a defined invariance under specified spatial symmetry groups, both for the mean or
%expectation and for each single realisation.  
The plan of the paper is as follows.
In \fsec{Sec:Problem} we formulate the overall problem more precisely and spell
out the desiderata for the modelling and generation of symmetric positive definite (SPD) 
matrices, as well as for the mean, based on requirements for the metric.
The desiderata for the mean or expectation will be set forth in \fsec{Ssec:means-spd}.
As we shall use the spectral decomposition in \fsec{Sec:StoModel}
to separate orientation and strength / size, 
and an appropriate definition of mean and corresponding modelling approach,
one has also to be able to generate random orthogonal transformations.   
%Additionally, we formulate some desiderata for a average or mean for SPD matrices,
%and base these on requirements for a metric for SPD tensors of material properties.
In the end the SPD matrix at each point in the domain for each realisation
will be a function of a few new numerical parameters, here the logarithms of
the eigenvalues as well as a description for the orthogonal eigenvector matrix,
which we shall base on the Rodrigues formula via the so called rotation vector.

It is these eigenvalue parameters as well as the rotation vector which vary
randomly and in the domain, \rvsi{a similar approach as in \cite{Grigoriu2016},
where the angles describing the eigenvectors are random fields.}
On this level, the step to spatially random fields is relatively simple, and one may borrow
the idea formulated in \cite{GuilleminotSoize2011, GuilleminotSoize2012,  GuilleminotSoize2012a, 
GuilleminotSoize2013, GuilleminotSoize2013a, GuilleminotSoize2014, Grigoriu2016},
and regard the representation
parameters as random fields without any positivity or invariance constraints.  This then 
automatically produces random SPD tensor fields when played back via the spatially point-wise
exponential map on the product Lie group, and from these to the SPD tensors.

As the spatial variation of the random tensor fields is well known from the
%spectral decomposition of the 
corresponding correlation functions---no matter what
kind of mean was used for the computation of the correlation function, this can be
converted to a joint correlation function for the log eigenvalues and the
rotation vector.  Spatially homogeneous and invariant-wise isotropic tensor fields,
which is the subject of \cite{malyarenko_tensor-valued_2019}, are a special case
of this more general approach.  As may be gleaned from this short overview, 
ensuring positive definiteness and the proper invariance at each point and realisation
is a main concern, so we focus mainly on this part.
Random fields of the parameters can then be generated through well known synthesis
techniques connected to the Karhunen-Loève representation (e.g.\ \cite{hgm07}),
or even more general representations with tight frames \cite{LuschgyPages2009},
where no nonlinear constraints like positivity have to be observed,  
%Although, here we will not operate directly with SPD tensors, but instead introduce 
%a convenient representation based on Lie algebras, as this is a linear space,
%and the random fields may then be introduced on the level of the representation,
yielding the SPD matrix through a memoryless transformation.

We also have an application in mind, which will be used in numerical examples to demonstrate
the workings of the different modelling assumptions.  This is the thermal conductivity
of bones  --- a highly anisotropic material ---
which is important in some surgical procedures 
\cite{mediouni_optimal_2017, augustin_cortical_2012}.  
It is described in \fsec{Sec:Application}.
The results of the computations are shown in \fsecs{Sec:2DResults}{Sec:3DResults}, which then
displays the effect on actual computational results of different modelling assumptions
about the anisotropy.  While \fsec{Conclusion} concludes, some additional material on
distance metrics for symmetric positive definite matrices is collected in 
Appendix~\ref{App:metrics-SPD}.

%% file: problem.tex
\section{Problem description}\label{Sec:Problem}
Let $\spatialdomain\subset\Euclideanspace$ be a bounded domain in a $d$-dimensional Euclidean
space $\Euclideanspace$ (here $d=2$ or $d=3$), such that the behaviour of a physical 
system in the domain $\spatialdomain$ is governed by an abstract equation
\begin{linenomath}
\begin{equation}\label{Eq:deterministic}
    \diffoperator(\genSecTensor,\dispvector) = \vek{f}.
\end{equation}
\end{linenomath}
Here $\dispvector \in \solspace$ describes the state of the system lying in a Hilbert 
space $\solspace$ (for the sake of simplicity), 
$\diffoperator$ is a---possibly non-linear---operator modelling the physics of the system, 
and $\vek{f}\in {\solspace}^*$ is some external influence (action/excitation/loading).

Furthermore, the coefficient $\genSecTensor$ is regarded as a material tensor which 
represents  intrinsic physics-based properties of materials, such as
thermal conductivity, magnetic permeability, chemical diffusivity, and similar
\cite{nye_physical_1984, malgrange_symmetry_2014}.  It is modelled as a tensor-valued 
field, at each point in the domain $\spatialdomain$ 
taking values in the open convex cone of real-valued second order positive-definite tensors: 
%\footnote{Note that, in 
%this article, $\genSecTensor$ has only 4 components and hence, does not account for 
%reduced elasticity tensor, which has 9 components for the considered geometry $\spatialdomain$.}
% 
%(\twoD\otimes\twoD) C_{ij} e_i\otimes e_j (\twoD)^{\otimes2}
\begin{linenomath}
\begin{equation}
  \PosSymSecRankSpace:=\{ \genSecTensor\in (\Euclideanspace\otimes\Euclideanspace) \mid 
    \genSecTensor=\genSecTensor^T, \,\vek{z}^T \genSecTensor  \vek{z}>\vek{0}, \forall 
    \vek{z}\in\Euclideanspace\setminus\{\vek{0}\} \}.
\end{equation}
\end{linenomath}

As a simple example of such an equation, the reader may think of the well known
thermal conduction equation --- the example to be used in the numerical experiments:
\begin{linenomath}
\begin{equation}\label{Eq:deterministic-thm}
 -\nabla\cdot(\DFcond \nabla T(\vek{x})) = f(\vek{x})\quad \forall \vek{x}\in\spatialdomain,
\end{equation}
\end{linenomath}
plus some appropriate boundary conditions on $\partial \spatialdomain$.  Here the system
state is described by
the temperature $T(\vek{x})$ at each point $\vek{x}\in\spatialdomain$, the `loading' 
$f(\vek{x})$ are thermal sinks and sources, and $\genSecTensor(\vek{x})=
\DFcond\in\PosSymSecRankSpace$ is the second order tensor field of thermal conductivities.
%As already mentioned in \fsec{Introduction}, for the sake of simplicity we shall only
%deal with spatially constant random tensors, so that we would take here 
%$\DFcond \equiv \conducttensor\in\PosSymSecRankSpace$.

\subsection{Stochastic model} \label{Ssec:stoch-model}
Coming back to the general model \feq{Eq:deterministic},
due to uncertainties in the exact value of $\genSecTensor$, which may be due to
a highly heterogeneous material, or simply a lack of knowledge about the
exact value, we assume that a probabilistic model is adopted, so that it 
is modelled as a  second-order tensor valued random field, i.e.\ a mapping  
\begin{linenomath}
\begin{equation}\label{Eq:Tmapping}
\genSecTensor(\vek{x},\event):\spatialdomain\times\samplespace\rightarrow\PosSymSecRankSpace
\end{equation}
\end{linenomath}
on a probability space $(\samplespace,\sigmaalgebra,\probability)$, where $\samplespace$ 
represents the sample space containing all outcomes $\event\in\samplespace$, $\sigmaalgebra$
is the $\sigma$-algebra of measurable subsets of $\samplespace$, and $\probability$ is 
a probability measure. 
Similarly, one may assume the boundary or initial conditions, as well as external actions 
to be uncertain. However, this is secondary to the purpose of this paper and will 
not be covered further.  To avoid unnecessary clutter in the notation the
argument $\vek{x}\in\spatialdomain$ will be omitted, unless it is necessary for an
understanding of the expressions.  The same will often be true of the argument
$\event\in\samplespace$.

Evidently, this change transforms \feq{Eq:deterministic} to its 
stochastic counterpart, taking the form
\begin{linenomath}
\begin{equation}  \label{Eq:deterministic-stoch}
\diffoperator(\genSecTensor(\event),\dispvector(\event)) = \vek{f},
\end{equation} 
\end{linenomath}
in which $\dispvector(\event)$ %:\samplespace\rightarrow\Euclideanspace$ 
is a stochastic response, due to the randomness in $\genSecTensor(\event)$.
Here we shall not be concerned on how this stochastic response may be computed
(see e.g.\ \cite{hgm07}),
but primarily on how to represent the symmetric positive definite (SPD) tensor
$\genSecTensor(\event)\in\PosSymSecRankSpace$.

\subsection{Representation requirements} \label{Ssec:group-inv}
One important point in the description of material properties of 
$\genSecTensor(\event) \in\PosSymSecRankSpace$ are symmetries in the
sense of invariance against certain transformations.  It is well known that the collection
of all such transformations form a group $G$ \cite{nye_physical_1984,
malyarenko_tensor-valued_2019}, and here we shall only be concerned with so-called point
groups.  The representation of any such subgroup $G \subseteq \mrm{O}(d)$ of the spatial 
orthogonal group $\mrm{O}(d)$ of $\Euclideanspace$ is fairly direct, and invariance
here means that for all $\vek{R}\in G$ one has 
$\vek{R}^T\genSecTensor(\event) \vek{R} = \genSecTensor(\event)$, a homogeneous 
linear constraint on $\genSecTensor(\event)$.  

The set of tensors $\genSecTensor \in\SymSecRankSpace$ which satisfy such an invariance
\rvsi{is} easily seen to form a subspace in $\SymSecRankSpace$ 
(e.g.\ \cite{Cowin2013, malyarenko_tensor-valued_2019}), the invariant subspace 
$\SymSecRankSpace_G \subseteq \SymSecRankSpace$ of the action of $G$.  Hence the 
distribution of the  $\genSecTensor(\event)$ lives only on this subspace, or more precisely on
$\PosSymSecRankSpace_G := \PosSymSecRankSpace \cap \SymSecRankSpace_G$, 
as there will be no realisations outside of it.  

For the \emph{mean} 
$\tmeanC_{\vartheta} := \Expvalue_{\vartheta}(\genSecTensor(\cdot))$, which is an
averaging operation, it is possible (e.g.\ \cite{Cowin2013, malyarenko_tensor-valued_2019}) that
it is invariant w.r.t.\ a larger group $G_m\supseteq G$.  
\rvsi{This property of the mean possibly having a larger symmetry group is commonly
envoked e.g.\ when considering random assemblages of anisotropic crystals \cite{Cowin2013}.}
As was noted already in
\fsec{Introduction}, it is not all obvious which kind of mean to take here,
as it depends on how $\PosSymSecRankSpace$ is seen as a metric manifold
with metric $\vartheta$; more about that will follow further below in
\fsec{Ssec:means-spd} and in Appendix~\ref{App:metrics-SPD}.  This 
implies that 
\begin{linenomath}
\[
\tmeanC_{\vartheta}\in\SymSecRankSpace_{G_m} \subseteq \SymSecRankSpace_G, \; \text{ whereas each } \;
\genSecTensor(\event) \in \SymSecRankSpace_G \subseteq \SymSecRankSpace, 
\]
\end{linenomath}
as due to $G_m\supseteq G$ the $G_m$-invariant subspace $\SymSecRankSpace_{G_m}$
must be a subspace of the $G$-invariant subspace $\SymSecRankSpace_G$.  All this
leads to the following requirements for the parametrisation of $\genSecTensor(\event)$:
\begin{enumerate}
\item Each  $\genSecTensor(\event)$ has to be symmetric, i.e.\ 
      $\genSecTensor(\event)=\genSecTensor^T(\event)\in\SymSecRankSpace$.
\item Even under numerical approximations like truncation of series resp.\ 
      numerical quadrature, each realisation $\genSecTensor(\event)$ has to be SPD, i.e.\ 
      for all $\vek{z}\in\Euclideanspace\setminus\{\vek{0}\}:\, 
      \vek{z}^T\genSecTensor(\event)\vek{z} > 0$.
\item $\forall \event\in\samplespace: \genSecTensor(\event)\in\PosSymSecRankSpace_G$ 
      has to be invariant under some group of transformations
      $G \subseteq \mrm{O}(\Euclideanspace)$, the invariance or symmetry class of 
      each realisation.
\item The appropriate mean $\tmeanC_{\vartheta} := \Expvalue_{\vartheta}(\genSecTensor(\cdot))
      \in\PosSymSecRankSpace$ has to be invariant under a possibly larger group of transformations
      $G_m$, with $\mrm{O}(d) \supseteq G_m\supseteq G$, the invariance or symmetry class of 
      the mean $\tmeanC_{\vartheta}$.
%      , where $G\subseteq G_m \subseteq \mrm{N}(G) \subseteq \mrm{O}(d)$;
%      here $G_m$ has to be `between' the symmetry group $G$ for each realisation and the
%      normaliser subgroup $\mrm{N}(G)$ ($ =  \{ \vek{Q}\in\mrm{O}(d) \mid
%      \forall \vek{R} \in G: \vek{Q}\vek{R}\vek{Q}^T \in G \}$, see also
%      \cite{malyarenko_tensor-valued_2019}).
\end{enumerate}
The first, third, and fourth item are not difficult to satisfy, as they are linear constraints
in the vector space of all $d \times d$ matrices, all one has to do is to pick an appropriate
basis.  The first item is almost trivial,  and the second item is elaborated
on in \fsec{Ssec:ensuring-spd}.  The third and fourth point are discussed in 
\fsec{Ssec:group-inv-detl}, whereas the question of the mean in the fourth item
will be addressed in the following \fsec{Ssec:means-spd}.
%$\F{gl}(d):=\Euclideanspace \otimes \Euclideanspace$, and $\SymSecRankSpace\subset\F{gl}(d)$ 
%is a linear subspace.  
%One may note that while  $\dim (\F{gl}(d))=d^2$, for the subspace $\SymSecRankSpace$
%one only needs $\dim (\SymSecRankSpace) = d(d+1)/2$ parameters.

\subsection{Means for symmetric positive definite matrices} \label{Ssec:means-spd}
The second and third items on the above requirements list will be discussed in more detail
later, but  to discuss the fourth item, one has to first define the kind of mean
or average resp.\ expected value one wants to use.  As was already alluded to in \fsec{Introduction},
this is very much connected to the way in which distances on $\PosSymSecRankSpace$ are measured. 
And computing a weighted mean between two items is actually a way of interpolating. 
The next few paragraphs are meant to illustrate and remind the reader of this connection.
 
It has been well known since antiquity that there are different ways to compute a 
mean or an average, one just has to look at the classical Pythagorean means.   
The arithmetic or Euclidean mean $\bar{\vek{x}}$ of a set of vectors 
$\vek{x}_1,\dots, \vek{x}_n \in \Real^k$ ---best known by almost 
anyone is the case $k=1$ ---with standard Euclidean squared norm
$\nd{\vek{x}}^2 = \vek{x}^T\vek{x}$ and corresponding Euclidean distance
$\vartheta_E(\vek{y},\vek{x}) = \nd{\vek{y}-\vek{x}}$, is given by
\begin{linenomath}
\begin{equation}  \label{Eq:vec-var-w-min}    
    \Expvalue_{\vartheta_E}[\{\vek{x}_j\}] :=
   \bar{\vek{x}} := \frac{1}{n} \sum_{j=1}^n \vek{x}_j = \arg \min_{\vek{x}\in\Real^k} 
   \Psi(\vek{x}) = \arg \min_{\vek{x}\in\Real^k} \frac{1}{n} \sum_{j=1}^n \nd{\vek{x}-\vek{x}_j}^2 ,
\end{equation}
\end{linenomath}
which displays right away the variational characterisation of the mean $\bar{\vek{x}}$
%\feq{Eq:arith-mean}, one defines the `$x$-based variance' as
%\begin{equation}  \label{Eq:var-def}
%    \Psi(x) =  \frac{1}{n} \sum_{j=1}^n (x_j- x)^2 ,
%\end{equation}
by elementary calculus. % shows that the arithmetic or Euclidean mean \feq{Eq:arith-mean} has the 
%following variational characterisation:
%\begin{equation}  \label{Eq:arith-mean-var}
%   \bar{x} = \arg \min_{x\in\Real}  \Psi(x) = \frac{1}{n} \sum_{j=1}^n x_j,
%\end{equation}
The minimal value $\Psi(\bar{\vek{x}})$ of the function $\Psi$ is called the 
variance of the set $\vek{x}_1,\dots, \vek{x}_n \in \Real^k$.
Observe that the terms $\nd{\vek{x}_j- \vek{x}}^2 =\vartheta_E(\vek{x}_j,\vek{x})^2$ 
in \feq{Eq:vec-var-w-min} are the squares of the Euclidean distances 
$\vartheta_E(\vek{x}_j,\vek{x})$ between $\vek{x}_j$ and $\vek{x}$.
One may also note that for $n=2$ one gets from \feq{Eq:vec-var-w-min} that $\bar{\vek{x}} = 
\vek{x}_1 + (1/2)(\vek{x}_2-\vek{x}_1)$
is the half-way interpolation on the `shortest path'---the straight line---between 
$\vek{x}_1$ and $\vek{x}_2$, which is an elementary
example of the connection between averaging and interpolation.

The \feq{Eq:vec-var-w-min} is easily generalised to not necessarily equal (probability) weights
$w_j > 0$, $\sum_{j=1}^n w_j = 1$, instead of the constant $w_j \equiv 1/n$,
and vectors $\vek{x}_j \in \C{V}$ in any Euclidean vector space $\C{V}$ with
inner product $\ip{\vek{x}}{\vek{y}}_{\C{V}}$ and corresponding norm, giving 
$\bar{\vek{x}}^w  = \sum_{j=1}^n w_j \vek{x}_j$.
Hereby it was assumed that---as is the case for an isotropic Euclidean distance---the gradient  
$\nabla_{\vek{x}}\ip{\vek{x}}{\vek{x}}_{\C{V}}$ of the inner product is $2 \vek{x}$.

The generalisation of \feq{Eq:vec-var-w-min} %as in \feq{Eq:var-P-min} 
corresponding to a general 
probability distribution $\probability$ of a ${\C{V}}$-valued random variable $\vek{Y}$ would be
\begin{linenomath}
\begin{multline}  \label{Eq:vec-var-P-min}
    \Psi_{\probability}(\vek{x}) = \int_{\C{V}}  \vartheta_E(\vek{y},\vek{x})^2 \,\probability(\di\vek{y})
    = \int_{\C{V}}  \nd{\vek{y}-\vek{x}}_{\C{V}}^2\,\probability(\di\vek{y})
     =  \int_{\C{V}}  \ip{\vek{y}-\vek{x}}{\vek{y}-\vek{x}}_{\C{V}} \,\probability(\di\vek{y}) 
      \\ \Longrightarrow\quad
     \olsi{\vek{Y}} = \Expvalue_{\vartheta_E}(\vek{Y}) =  \arg \min_{x\in{\C{V}}}  
     \Psi_{\probability}(\vek{x}) =  \int_{\C{V}} \vek{y} \, \probability(\di \vek{y}),
\end{multline}
\end{linenomath}
and again the minimum $\Psi_{\probability}(\olsi{\vek{Y}})$ 
is the variance of the distribution $\probability$.
%Obviously the expressions \feeqs{Eq:vec-var-w-min}{Eq:vec-var-P-min} can easily be
%generalised for an abstract inner product space.

In case the underlying set does not have a vector space structure %like e.g.\ $\Real^N$, 
but only that of a metric space, i.e.\ when only a distance metric $\vartheta$ is defined,
or, like in our case, another distance than the Euclidean is to be used,
one may define the so-called Fréchet mean or average attached to that metric 
\cite{NielsenBhatia2013, PennecSommerFletcher2020, GuiguiEtal2023, groissJungSchwman2023} through a 
variational characterisation similar to \feqs{Eq:vec-var-w-min}{Eq:vec-var-P-min},
by replacing the Euclidean metric $\vartheta_E$ in those expressions by $\vartheta$.
Thus one may \emph{define} the \emph{Fréchet-Karcher mean} on 
such a metric space $\C{S}$ with a distance function $\vartheta$ 
%is again one where there is a collection of $n$ items 
%$s_1,\dots,s_n\in\mathcal{S} $ with `weights' or probabilities $w_1,\dots,w_n$,
%the mean or expected value is defined as the minimiser of the `variance functional'
%\begin{equation}  \label{Eq:FrechetMean-gen-w}
%    \Psi_w(s) =  \sum_{j=1}^n w_j \vartheta(s_j,s)^2 
%    \quad \Longrightarrow \quad    
%    \bar{s}^w = \Expvalue_{\vartheta}[\{w_j,s_j\}] :=  \arg \min_{s\in\mathcal{S}}  \Psi_w(s) ,
%\end{equation}
%with the variance defined by $\Psi_w(\bar{s}^w)$.
%In case $\mathcal{S}$ is a metric space and 
which is at the same time a measure space with a general probability measure  
$\probability$ for an $\C{S}$-valued random variable $R$ as
\begin{linenomath}
\begin{equation}  \label{Eq:FrechetMean-gen-P}
    \Psi_{\probability}(s) =  \int_{\mathcal{S}} \vartheta(r,s)^2 \, \probability(\di r)
    \quad \Longrightarrow \quad    
    \olsi{R} = \Expvalue_{\vartheta}[R] :=  \arg \min_{s\in\mathcal{S}}  \Psi_{\probability}(s) ;
\end{equation}
\end{linenomath}
and again the variance is defined  by the minimum value $\Psi_{\probability}(\olsi{R})$.
If a unique minimiser exists, for general metric spaces this is called the Fréchet mean and
variance, and for Riemannian manifolds a local minimiser of the variance functionals
in \feq{Eq:FrechetMean-gen-P} is also sometimes called the Karcher mean 
\cite{NielsenBhatia2013, PennecSommerFletcher2020, GuiguiEtal2023, groissJungSchwman2023}.

%As an aside, one may observe that on the positive numbers $\PosRealStrict$ with
%distance function $\vartheta_L(x,y) = |\log x - \log y|$, one obtains the geometric
%mean, whereas with the distance function $\vartheta_H(x,y) = |x^{-1} - y^{-1}|$ on
%$\PosRealStrict$ one obtains the harmonic mean.
%If $\mathcal{S}$ is also a Riemannian manifold, so that $\vartheta(r,s)$ can be defined
%as the length of a shortest geodesic joining $r$ and $s$, then the case $n=2$ in
%\feq{Eq:FrechetMean-gen-w}, setting $w_2=\alpha$ and $w_1=1-w_2=1-\alpha$, leads as
%before for \feq{Eq:vec-var-w-min} to the conclusion that the Fréchet-($\alpha$-)mean 
%$\bar{s}^\alpha \in \mathcal{S}$ of $s_1, s_2 \in \mathcal{S}$ is a point on that
%shortest geodesic satisfying
%\begin{equation}  \label{Eq:Fechet-w-interp}
%    \vartheta(s_1,\bar{s}^\alpha) = \alpha \, \vartheta(s_1,s_2); \quad
%    \vartheta(\bar{s}^\alpha,s_2) =  (1-\alpha) \, \vartheta(s_1,s_2),
%\end{equation}
%and is the interpolation of the `fraction' $\alpha$ of the path length on the `shortest path', 
%again showing the connection between weighted averaging and weighted interpolation.
%

These foregoing observations show the intimate connection between distance metric, mean or average
and expectation, and interpolation.  The properties of the Fréchet mean or average
depend completely on the underlying metric, see \feqs{Eq:vec-var-w-min}{Eq:FrechetMean-gen-P}.
Turning now to positive definite tensors $\genSecTensor\in\PosSymSecRankSpace$,
as mentioned already in \fsec{Introduction}, the Euclidean or Frobenius resp.\ Hilbert-Schmidt norm
$\nd{\genSecTensor}^2_F := \ip{\genSecTensor}{\genSecTensor}_F := \tr(\genSecTensor^T\genSecTensor)$
(cf. \feq{Eq:Frob-norm} in Appendix~\ref{App:metrics-SPD}) has long been recognised to
result in some undesirable properties
\cite{PennecFillardAyache2004, AndoLiMathias2004, Moakher2005, ArsignyFillardPennecEtAl2006, 
ArsignyFillardPennecEtAl2007, drydenEtal2009, DrydenEtal2010, FujiiSeo2015}, 
especially the  \emph{swelling, fattening, and shrinking effects} in interpolation
resp.\ averaging \cite{schwartzman_random_nodate, JungSchwartzmanGroisser2015, 
schwartzman_lognormal_2016, GroisserJungSchwartzman2017, groissJungSchwman2017, FeragenFuster2017}.
It is therefore of considerable interest to consider different metrics on $\PosSymSecRankSpace$
\cite{PennecSommerFletcher2020, NielsenBhatia2013, FeragenFuster2017}.    
In order not to distract the reader's
attention too much by reviewing different metrics used on $\PosSymSecRankSpace$, we
relegate this material to Appendix~\ref{App:metrics-SPD}.  Having chosen some distance metric
$\vartheta_D$ with desirable properties on $\PosSymSecRankSpace$, 
one defines the $\SecTensorB$-based Fréchet variance analogous to 
\feqs{Eq:vec-var-w-min}{Eq:FrechetMean-gen-P}:
\begin{linenomath}
\begin{equation}\label{Eq:FrechetVar}
  \Psi(\SecTensorB) = \sum_{k=1}^n w_k \, \vartheta_D(\SecTensorB,\genSecTensor_k)^2, \;
  \text{ or } \;
   \Psi(\SecTensorB) = \int_\Omega  \vartheta_D(\SecTensorB,\genSecTensor(\event))^2 \, 
   \probability(\di \event) .
\end{equation}
\end{linenomath}
From this one obtains the definition of the corresponding Fréchet or Karcher mean resp.\ expectation
\begin{linenomath}
\begin{equation}\label{Eq:FrechetMean}
  \Expvalue_{\vartheta_D}[\genSecTensor] :=
  \tmeanC_{\vartheta_D} := \underset{\SecTensorB \in \PosSymSecRankSpace}{\arg\min} \Psi(\SecTensorB) .
\end{equation}
\end{linenomath}
Inserting this minimiser---if it exists---into the functional \feq{Eq:FrechetVar}, one obtains
the variance  of the $\PosSymSecRankSpace$-valued random variable
$\genSecTensor(\event)$ as the minimum value $\Psi(\tmeanC_{\vartheta_D})$.

% One may recall here that interpolation
%and the mean or average are closely related on metric spaces, i.e.\ sets with a distance measure.
%After all, the mean or average of two items is midway between them on a shortest path---a
%geodesic, and along the geodesic one may interpolate between two items.
%Obviously a vector space like $\SymSecRankSpace$ with the usual Frobenius norm
%$\nd{\genSecTensor}^2_F := \ip{\genSecTensor}{\genSecTensor}_F := \tr(\genSecTensor^T\genSecTensor)$
%is a metric space with the usual Euclidean metric, and straight lines are geodesics,
%where one takes for the Euclidean metric $\vartheta_F(\genSecTensor_1,\genSecTensor_2) 
%=\nd{\genSecTensor_1-\genSecTensor_2}_F$.  

Let us collect here some, in our view desirable, properties of a possible metric for material
tensors $\genSecTensor(\event)\in\PosSymSecRankSpace$, which then may point to
a certain kind of mean to use.  As an extreme case one could argue that the mean
$\tmeanC_{\vartheta_D}$ should be such, that if used in \feq{Eq:deterministic}, one obtains
the mean $\bar{\dispvector}$ of the stochastic solution of \feq{Eq:deterministic-stoch}:
\begin{linenomath}
\begin{equation}  \label{Eq:deterministic-stoch-homogen}
\diffoperator(\tmeanC_{\vartheta_D},\bar{\dispvector}) = \vek{f},
\end{equation} 
\end{linenomath}
This is in fact the goal of homogenisation techniques, 
see e.g.\ \cite{KozlovOleinikZhikov1994, Hornung1997},
but it is probably too much to require from a general kind of average, as it would 
depend on the model \feqs{Eq:deterministic}{Eq:deterministic-stoch} and thus be
highly problem dependent.  Hence we try and settle with some more general requirements
based on invariance considerations.  The invariance requirements for the metric
translate immediately into invariance requirements for the corresponding Fréchet mean.

From the point of physical modelling, if $\genSecTensor(\event)$ appears in a constitutive
model, there is typically an inverse formulation involving $\genSecTensor^{-1}(\event)$.
This means that one normally knows as much about the distribution of $\genSecTensor(\event)$
as about  the distribution of $\genSecTensor^{-1}(\event)$, as the
set $\PosSymSecRankSpace$ is stable under inversion.  So one property of a 
desired metric $\vartheta_D$ on $\PosSymSecRankSpace$ would be
%$\vartheta_D(\genSecTensor_1,\genSecTensor_2)=\vartheta_D(\genSecTensor_1^{-1},\genSecTensor_2^{-1})$,
%i.e.\ 
invariance under inversion.
Changing physical units would correspond to a scaling by a positive factor, so the desired metric
should be invariant under scaling.
%: $\vartheta_D(\alpha\genSecTensor_1,\alpha\genSecTensor_2)= 
%\vartheta_D(\genSecTensor_1,\genSecTensor_2)$ for all $\alpha>0$.
And as a change of the coordinate frame by an orthogonal transformation $\rotmat\in\mrm{O}(d)$
changes the constitutive tensor $\genSecTensor$ to $\rotmat^T\genSecTensor \rotmat$,
one would like the desired metric $\vartheta_D$ to be invariant under such transformations
in order to respect the principle of frame indifference.
%, i.e.\
%for all $\rotmat\in\mrm{O}(d)$: $\vartheta_D(\genSecTensor_1,\genSecTensor_2)= 
%\vartheta_D(\rotmat^T\genSecTensor_1 \rotmat,\rotmat^T\genSecTensor_2 \rotmat)$.
Collecting everything, the desiderata for a metric $\vartheta_D$, and hence for a 
Fréchet mean based on it via \feq{Eq:FrechetMean}, then are to satisfy for all 
 $\genSecTensor_1,\genSecTensor_2 \in \PosSymSecRankSpace$, $\alpha>0$, and all
 $\rotmat\in\mrm{O}(d)$:
\begin{itemize}
\item invariance under inversion: $\vartheta_D(\genSecTensor_1^{-1},\genSecTensor_2^{-1})=
                  \vartheta_D(\genSecTensor_1,\genSecTensor_2)$;
\item invariance under scaling: $\vartheta_D(\alpha\genSecTensor_1,\alpha\genSecTensor_2)=                       
                           \vartheta_D(\genSecTensor_1,\genSecTensor_2)$;
\item invariance under frame transformations: $\vartheta_D(\rotmat^T\genSecTensor_1 \rotmat,
        \rotmat^T\genSecTensor_2 \rotmat)=\vartheta_D(\genSecTensor_1,\genSecTensor_2)$.
\end{itemize}

By choosing a metric $\vartheta_D$ on $\PosSymSecRankSpace$ which satisfies those requirements, the
Fréchet-Karcher mean $\Expvalue_{\vartheta_D}[\genSecTensor]$ based on that metric will
automatically inherit all those properties, as can be easily gleaned from the very definition
above in \feqs{Eq:FrechetVar}{Eq:FrechetMean}.

The usual Frobenius resp.\ Euclidean distance which leads to the normal arithmetic
mean \feq{Eq:vec-var-P-min} satisfies the last requirement, but not the first two.
Additionally, the swelling, fattening, and shrinking effects that result from
the straight line interpolation or averaging in the vector space $\SymSecRankSpace$ of two 
elements of $\PosSymSecRankSpace$ are highly undesirable  
\cite{ArsignyFillardPennecEtAl2006, schwartzman_random_nodate, JungSchwartzmanGroisser2015, 
schwartzman_lognormal_2016, GroisserJungSchwartzman2017, groissJungSchwman2017, FeragenFuster2017}.  
Any Fréchet-Karcher mean based on a distance function which exhibits these
undesirable effects in interpolation would also exhibit these effects in the averaging
procedure, i.e.\ the thus defined mean would be more swollen, fatter, or shrunk than the
population it is derived from.  This short consideration
shows that the normal arithmetic or Euclidean mean based on the Frobenius distance
$\vartheta_F$ (\feq{Eq:Euclid-dist} in Appendix~\ref{App:metrics-SPD}) is not ideal 
for SPD material tensors.

What comes out of from the more detailed considerations in Appendix~\ref{App:metrics-SPD} 
regarding possible metrics on $\PosSymSecRankSpace$, is that all the other
distance metrics considered there, namely the geometric or affine-invariant distance
$\vartheta_G$ (\feq{Eq:aff-inv-dist} in Appendix~\ref{App:metrics-SPD}), 
the logarithmic or log-Euclidean distance
$\vartheta_L$ (\feq{Eq:log-Euclid-dist} in Appendix~\ref{App:metrics-SPD}), 
and the scaling-rotation distance---which
we will often refer to as just the scaling distance---
$\vartheta_S$ (\feq{Eq:scal-dist} in Appendix~\ref{App:metrics-SPD}), 
and which in turn is based on the product
Lie group distance $\vartheta_P$ (\feq{Eq:prod-Lie-dist} in Appendix~\ref{App:metrics-SPD}), 
all satisfy all of the above three requirements.  
This means that for a choice between them additional considerations are necessary.
These are how to ensure the requirement listed above in \fsec{Ssec:group-inv}, namely how to
ensure positive definiteness through an appropriate representation, and at the same time
ensure the symmetry class of the corresponding mean and each single realisation.
\rvsi{One may note that we have defined $\PosSymSecRankSpace$ as an \emph{open convex cone}
in $\SymSecRankSpace$, which means that all elements $\genSecTensor\in\PosSymSecRankSpace$ 
are positive definite.  Thus any average or mean, which in a vector space is located
in the convex hull of the support of the probability distribution will also be positive definite.
In our approach in \fsec{Sec:StoModel} this is ensured by modelling the logarithms of the
eigenvalues.
A slightly different situation may arise when \cite{drydenEtal2009, DrydenEtal2010}
experimental or measured covariance matrices are averaged.  Such matrices may be rank 
deficient, i.e.\ singular, or symmetric positive but \emph{not} definite, and thus geometrically 
on the boundary $\partial\,\PosSymSecRankSpace$ of the convex cone of positive definite matrices.
In this case the arithmetic mean may also be singular and thus not positive definite, but the 
other Fréchet means may more readily give an average in the interior of the cone, and
hence a positive definite mean.}

\rvsi{In the following \fsec{Sec:StoModel}, where the modelling of SPD tensors
based on the chosen representation is described in more detail, such a Fréchet mean
will be used.  To avoid the afore mentioned swelling, fattening, and shrinking effects,
the distance function we shall employ for this mean will be the scaling distance---
$\vartheta_S$ (\feq{Eq:scal-dist} in Appendix~\ref{App:metrics-SPD}), and often
for simplicity just the product Lie group distance $\vartheta_P$ (\feq{Eq:prod-Lie-dist} in 
Appendix~\ref{App:metrics-SPD}) on which it is based, and which agrees with it for
not too large variations---as it allows the desired fine control of
size / strength and orientation, as well as making it easy to ensure positive
definiteness.}

%% file: modelling.tex
% !TEX root = ../paper_tensors.tex
% !TEX encoding = UTF-8 Unicode

\section{Stochastic modelling}\label{Sec:StoModel}
Methods to ensure positive definiteness will be briefly considered in \fsec{Ssec:ensuring-spd}.
We shall settle for an approach based on the exponential map, which is most easily handled in
the spectral decomposition representation,  %This kind of representation 
\rvsi{and which} is also well suited
to address the problem of symmetry invariance. %, which will be done in \fsec{Ssec:group-inv-detl}.
\rvsi{The tensor $\genSecTensor(\event)$ ---at each location is modelled in \fsec{Ssec:group-inv-detl}
}
%The outcome of this approach is that the tensor $\genSecTensor(\event)$ at each spatial location
%may be described 
by a number---depending on the symmetry class---of parameters, which can be
clearly separated into scaling or `size' parameters, and \rvsi{those} which determine the orientation
of the symmetry axes.  \rvsi{This description is restricted to  second-order spatial tensors 
like the heat conduction tensor, higher order tensors like the fourth-order elasticity
tensor require additional parameters for the eigen-strain distribution.}
Although this is not the main topic of this paper, in \fsec{Ssec:spatial-model}
it is sketched how, borrowing an idea from \cite{GuilleminotSoize2011, GuilleminotSoize2012,
GuilleminotSoize2012a, GuilleminotSoize2013, GuilleminotSoize2013a, GuilleminotSoize2014, 
Grigoriu2016}, these parameters can be modelled as random fields without any positivity 
or invariance constraints,
thus giving a random tensor field $\genSecTensor(\vek{x},\event)$ through a memoryless 
point-wise transformation.

\subsection{Ensuring positive definiteness and group invariance} \label{Ssec:ensuring-spd}
The second item on the list of requirements for representation above in \fsec{Ssec:group-inv}, 
that each $\genSecTensor(\event)$ is SPD,  can pose some problems, as 
$\PosSymSecRankSpace\subset \SymSecRankSpace$ is not a subspace, but geometrically an open
convex cone.  
%This means that still $\dim (\PosSymSecRankSpace) = d(d+1)/2$ parameters are 
%needed, but they are not necessarily free resp.\ independent.  
Therefore we shall look only at
representations which can ensure that the resulting $\genSecTensor(\event)$ is SPD,
even under numerical approximations.  Some of these representations are actually closely
connected to some of the metrics described in Appendix~\ref{App:metrics-SPD}; and although the kind of
metric and associated mean and the representation are formally independent, some representations
make it easier to compute certain means.  
%All of these representations
%take advantage of the fact that $\SymSecRankSpace$ is an algebra as already mentioned
%in \fsec{Introduction}, which means that one may compute polynomials, power series and 
%in fact any real function defined on $\Real$ ---as the algebra is finite dimensional---applied 
%to elements of the algebra \cite{segalKunze78, yosida-fa-1980}.

The first representation we shall look at is the exponential map,  e.g.\ see 
\cite{ArsignyFillardPennecEtAl2006, ArsignyFillardPennecEtAl2007, GuilleminotSoize2013a, 
NouySoize2014, JungSchwartzmanGroisser2015, schwartzman_lognormal_2016, Grigoriu2016,
GroisserJungSchwartzman2017, groissJungSchwman2017} for some recent examples, which is 
analytic and one-to-one and onto, i.e.\ it can be taken as a global chart.  For a 
$\symtensor \in \SymSecRankSpace$ we set $\genSecTensor = \exp(\symtensor)\in\PosSymSecRankSpace$.
%\begin{equation}    \label{Eq:exp-gen-def}
%   \exp: \SymSecRankSpace\ni\symtensor  \to \genSecTensor = \exp(\symtensor)\in\PosSymSecRankSpace ,
%\end{equation}
This can be computed
%which can be defined 
either via spectral calculus, or simply through the everywhere convergent
power series.  One may also note that thanks to the Cayley-Hamilton 
theorem the series can actually be written as a finite sum.  
This representation also makes it very easy to compute the logarithmic or log-Euclidean distance
$\vartheta_L$, \feq{Eq:log-Euclid-dist} in Appendix~\ref{App:metrics-SPD}.

The main result is that whatever errors may be committed by approximating
$\symtensor(\event)$, by computing the exponential $\exp(\symtensor(\event))$ 
in the end the result is guaranteed to be SPD.  One  advantage of
the representation is that $\SymSecRankSpace$ is a linear vector space, 
and $\symtensor\in\SymSecRankSpace$ is completely free in that vector space, which 
makes stochastic modelling a lot easier.
%, and it has obviously the right number of parameters $\dim(\SymSecRankSpace) = d(d+1)/2=
%\dim(\PosSymSecRankSpace)$.
One should also note that ensuring any kind of group invariance
is very easy, as $\genSecTensor$ and its logarithm
$\symtensor$ have exactly the same group invariance.  For a fixed orientation, this
then manifests itself as requiring $\symtensor$ to be in some linear subspace of
$\SymSecRankSpace$. 

One may note that by wanting to model $\symtensor=\log(\genSecTensor)$, 
the tensor $\genSecTensor$ has physical units, so that it does
not directly make sense to take its logarithm.  This is usually dealt with by some kind
of scaling with a representative inverse tensor $\meanC^{-1}$  which has the right invariance 
properties---this could be the mean $\tmeanC_{\vartheta}$ one wants to use---as a scaling, 
and using its spectral decomposition (cf.\ also \feq{Eq:SpecDecSPD}) 
$\meanC = \meaneigvect\meaneigval\meaneigvect^T$
one may produce a  dimensionless tensor
\begin{linenomath}
\begin{equation} \label{Eq:Scaling}
  \wmeanC(\event)   := \meaneigvect\meaneigval^{-1/2}\genSecTensor(\event) \meaneigval^{-1/2}
  \meaneigvect^T , 
%     \quad \text{or even }\quad 
%  \genSecTensor(\event)  := \rotmat(\event)\widehat{\genSecTensor}^{1/2}
%      \mathring{\genSecTensor}(\event) \widehat{\genSecTensor}^{1/2} \rotmat(\event),
\end{equation}
\end{linenomath}
so one may take its logarithm $\wmeanH(\event)=\log(\wmeanC(\event))$.  
Usually it is tacitly assumed
that such a transformation or scaling to a dimensionless form was performed before
computing functions of $\genSecTensor$.
%As this problem may appear with other representations, it will be discussed further later
%in \fsec{Ssec:group-inv-detl}.

The second kind of representation we shall consider, is the approach (see e.g.\ 
\cite{soize_nonparametric_2000, soize_2006, GuilleminotSoize2011, GuilleminotSoize2012,
Grigoriu2016}) 
where the problem to ensure that
each $\genSecTensor(\event)$ is SPD is solved by representing the tensor as a generalised square.  
The very definition of being positive in an algebra  
is that the element is a generalised square, i.e.\ 
$\genSecTensor(\event) = \vek{G}^T(\event)\vek{G}(\event) $ with $\vek{G}(\event)\in \mrm{GL}(d)$ 
non-singular.
In this approach as described in the cited literature, the matrix $\vek{G}(\event)$ is chosen as
upper triangular, and thus one is certain to use only the right number of parameters.
In \cite{GuilleminotSoize2013a, GuilleminotSoize2014}   this representation as a generalised
square is combined with the exponential representation described above.
\rvsi{Group invariance due to material symmetries in this approach has only been considered
for the mean \cite{soize_nonparametric_2000, soize_2006, GuilleminotSoize2011, GuilleminotSoize2012,
Grigoriu2016}.  A construction similar to \feq{Eq:Scaling} has been used for this purpose,
in the form of factoring the mean $\meanC = \widehat{\vek{G}}^T\widehat{\vek{G}}$ 
which has the required group invariance (e.g.\ Choleski factorisation), 
and then setting $\genSecTensor(\event) = \widehat{\vek{G}}^T 
\lmeanC(\event) \widehat{\vek{G}}$. %$\tmeanC, \wmeanC, \meanC, \lmeanC$
Here $\lmeanC(\event) \in \PosSymSecRankSpace$ is completely anisotropic, 
but SPD and with mean equal to the identity $\vek{I}$; 
cf.\ also \feqs{Eq:SpecDec-rep}{Eq:ExpSpecDec-rep}.
Now $\PosSymSecRankSpace$ is not a vector space, and any \rvsi{nonsingular} $\lmeanC \in 
\PosSymSecRankSpace$ has to satisfy the non-linear constraint $\det \lmeanC \neq 0$. 
By controlling the diagonal elements of an upper 
triangular matrix $\vek{L}(\event) \in \mrm{GL}(d)$ such that $\det \vek{L}(\event) > 0$
---the matrix $\vek{L}(\event)$ may thus be seen as a Cholesky-like factor of 
$\lmeanC (\event) = \vek{L}(\event)^T\vek{L}(\event) \in \PosSymSecRankSpace$
}
---one can indeed make sure that it is invertible, and this otherwise troublesome 
nonlinear constraint of non-vanishing determinant becomes simpler.

The third kind of representation, and this will be the one we will eventually use
in our computations, allows a fine control of orientation and `size'.
This is achieved by looking at the spectral decomposition of a 
$\genSecTensor\in\PosSymSecRankSpace$ and its logarithm $\symtensor=\log(\genSecTensor)$:
\begin{linenomath}
\begin{equation} \label{Eq:SpecDecSPD}
  \genSecTensor = \eigvect\eigval\eigvect^T,\quad
  \log(\genSecTensor) = \eigvect\log(\eigval)\eigvect^T = \symtensor
  =: \eigvect\symdiag\eigvect^T ,
\end{equation} 
\end{linenomath}
with $\eigvect\in\mrm{SO}(d)$ and the diagonal positive definite matrix $\eigval=\diag(\lambda^{(i)})\in
\mrm{Diag}^+(d)\subset\PosSymSecRankSpace$ of positive 
eigenvalues $\{\lambda^{(i)}> 0 \mid i=1,\dots,d \}$.
Here the eigenvector matrix  $\vek{Q}$ controls the orientation, whereas the diagonal 
matrix $\eigval$ (resp.\ $\symdiag = \diag(y^{(i)}) = \diag(\log(\lambda^{(i)})) =
 \log(\eigval)$) controls the scale, size, or strength of the tensor $\genSecTensor$. 
%resp.\ the diagonal matrix $\symdiag = \diag(y^{(i)}) = \diag(\log(\lambda^{(i)})) = 
%\log(\eigval) \in \F{diag}(d)\subset\SymSecRankSpace$ controls the `size' of 
%$\symtensor=\log\genSecTensor$.

In Appendix~\ref{Asec:scale-rot} the spectral decomposition \feq{Eq:SpecDecSPD} is also
used to investigate the scaling-rotation distance $\vartheta_S$ on $\PosSymSecRankSpace$.
Informally, \feq{Eq:SpecDecSPD} says that it induces the \rvsi{formal} representation
\cite{schwartzman_random_nodate, JungSchwartzmanGroisser2015, schwartzman_lognormal_2016, 
Grigoriu2016, GroisserJungSchwartzman2017, groissJungSchwman2017}
\begin{linenomath}
\begin{equation}  \label{Eq:repr-prod-Lie-2}
  \mrm{Diag}^+(d) \times \mrm{SO}(d) \ni (\eigval,\eigvect)
    \mapsto \genSecTensor=\eigvect\eigval\eigvect^T \in \PosSymSecRankSpace .
\end{equation}
\end{linenomath}
One should point out though
that the representation map \feq{Eq:repr-prod-Lie-2} is onto, but not necessarily one-to-one;
an example where this can be seen is the occurrence of multiple eigenvalues.

%While it was pointed out already that $\PosSymSecRankSpace$
%is not a Lie group under matrix multiplication, the subset $\mrm{Diag}^+(d)\subset\PosSymSecRankSpace$
%of positive definite diagonal matrices obviously \emph{is} a commutative Lie group with matrix 
%multiplication.  And as $\mrm{SO}(d)$ is well known to be a Lie group, the product on the left 
%hand side of \feq{Eq:repr-prod-Lie-2} is a product of Lie groups, and thus is again a Lie group.  
%And the advantage of having a representation on a Lie group is that the Lie group can 
%itself be easily represented on its Lie algebra---the tangent space at the identity and 
%thus importantly an unconstrained linear vector space---via the exponential map.
%In the vicinity of the origin the exponential map is an analytic diffeomorphism, and for
%the Lie group $\mrm{Diag}^+(d)$ it is actually a globally valid chart, similarly
%as for $\PosSymSecRankSpace$ in $\genSecTensor = \exp(\symtensor)$ before.

\rvsi{Unlike $\PosSymSecRankSpace$, the positive definite diagonal matrices
$\mrm{Diag}^+(d)\subset\PosSymSecRankSpace$ are a commutative Lie group under normal 
matrix multiplication.  Together with the Lie group of proper rotations $\mrm{SO}(d)$,
on the left hand side of \feq{Eq:repr-prod-Lie-2} is a product of Lie groups, and thus 
again a Lie group.  As already indicated, a Lie group can easily be represented via
the exponential map on its Lie algebra---the tangent space at the identity and 
an unconstrained linear vector space.  For $\mrm{Diag}^+(d)$ the Lie algebra are
the diagonal matrices $\F{diag}(d)\subset\SymSecRankSpace$, and for $\mrm{SO}(d)$
it is well known to be the skew symmetric matrices $\F{so}(d)$ 
($\skewsym = -\skewsym^T$), see the following \fsec{Ssec:group-inv-detl}.
Here the exponential map reduces to the matrix exponential,
trivial for diagonal matrices, and for proper orthogonal matrices it may be realised via 
the well known Rodrigues formula; see \fsec{Ssec:group-inv-detl}, 
\feeqs{Eq:Exp-rand-ax-rotation}{Eq:Rodrigues-inv}.  But the connection with Lie group
theory not only gives a theoretical underpinning for the representation,
at the same times it allows the definition of Riemannian structures and thus the 
definition of distances and Fréchet means based on it, see \fsec{Ssec:means-spd}
and Appendix~\ref{App:metrics-SPD}.}

%The Lie algebra of the commutative Lie group $\mrm{Diag}^+(d)\subset\PosSymSecRankSpace$ 
%is easily seen to be the subspace $\F{diag}(d)\subset\SymSecRankSpace$
%of diagonal (and thus symmetric) matrices.
%The Lie algebra of the special orthogonal 
%group $\mrm{SO}(d)$ is known to be the Lie algebra of skew symmetric matrices $\F{so}(d)$ 
%($\skewsym = -\skewsym^T$, see the following \fsec{Ssec:group-inv-detl}), and 
%thus the 
The spectral synthesis representation \feq{Eq:repr-prod-Lie-2} can be extended to the 
left with the exponential map, giving a two step modelling process:
\begin{linenomath}
\begin{align}   \label{Eq:repr-prod-Lie-exp} 
    \F{diag}(d) \times \F{so}(d) \; & \overset{\exp \times \exp}{\longrightarrow} \quad 
    \mrm{Diag}^+(d) \times \mrm{SO}(d) \hspace{1.7em}  \longrightarrow \; \PosSymSecRankSpace , \\  \nonumber
    (\symdiag,\skewsym) \;  & \;\; \longmapsto  \; (\eigval,\eigvect) = (\exp \symdiag,\exp  \skewsym) 
   \longmapsto  \genSecTensor=\eigvect\eigval\eigvect^T . %\hspace{3em}
\end{align} 
\end{linenomath}
Obviously $\eigval = \diag(\lambda^{(i)}) = \diag(\exp y^{(i)}) = \symdiag$ ensures that
$\genSecTensor$ is SPD. %, and the exponential of a diagonal matrix is easy to compute.
%The two steps of exponential map and spectral synthesis can actually be interwoven, 
%explicitly showing that the exponential map is used in the last step:
\rvsi{Interweaving the exponential and the spectral synthesis in \feq{Eq:repr-prod-Lie-exp},
the last step then shows also formally that the logarithm
$\symtensor = \log \genSecTensor$ is being modelled:}
\begin{linenomath}
\begin{align}   \label{Eq:repr-prod-Lie-exp-var}
    \F{diag}(d) \times \F{so}(d) \; & \overset{(\cdot,\exp)}{\longrightarrow}
   \;  \F{diag}(d) \times \mrm{SO}(d) \; \longrightarrow 
   \qquad  \SymSecRankSpace \quad \overset{\exp}{\longrightarrow} \quad \PosSymSecRankSpace, \\
   (\symdiag,\skewsym) \; & \, \longmapsto \;  (\symdiag,\eigvect=\exp \skewsym) \,
     \longmapsto \; \symtensor=\eigvect\symdiag\eigvect^T \,
     \longmapsto \, \genSecTensor=\exp \symtensor . \nonumber
\end{align} 
\end{linenomath}
%The important advantage we see in
%this last representation, vis-a-vis the first two considered before, 
\rvsi{The advantage we see in this last representation}
is that it separates
out the eigenvalues $\lambda^{(i)}$ resp.\ their logarithms $y^{(i)} = \log \lambda^{(i)}$
---which physically may be conductivities or resistances---from the angle-like parameters
defining $\skewsym\in\F{so}(d)$, which one may use for fine-tuning in the stochastic modelling.
\rvsi{A similar approach as regards the modelling of rotations through their rotation angles
has been proposed in \cite{Grigoriu2016}.}

As was already indicated earlier, and as easily seen
from \feq{Eq:SpecDecSPD}, the group invariance properties of of a symmetric tensor 
$\symtensor$ and its exponential $\genSecTensor = \exp \symtensor$ are the same, 
a fact which was also used in \cite{GuilleminotSoize2013a, GuilleminotSoize2014}.  
This means that all the invariance constraints in the Lie representation 
\feqs{Eq:repr-prod-Lie-exp}{Eq:repr-prod-Lie-exp-var}
to be proposed here may be applied on the Lie algebra level.  Thus, building 
on the spectral decomposition \feq{Eq:SpecDecSPD}, one  uses the
inverse of the scaling equation \feq{Eq:Scaling} above together with the
spectral decomposition \feq{Eq:SpecDecSPD}.  

\rvsi{The Lie representation also effortlessly combines with the task of ensuring group
invariance resulting from material symmetry.}  It was mentioned already that for 
second order symmetric tensors in 2D or 3D the situation is relatively simple
\cite{altenbach_tensor_2018, malyarenko_tensor-valued_2019},
\rvsi{as through group reduction} \cite{faessStiefel1992} the matrix
representation can be fully diagonalised.  Thus in 3D there are only three symmetry classes
\cite{malyarenko_tensor-valued_2019}, and
depending on the number of distinct eigenvalues, the symmetry subgroup relations are
\begin{linenomath}
\begin{equation}   \label{Eq:3D-symm-Hasse}
   [\text{orthotropic: 3 dEVs}] \hookrightarrow  [\text{plan isotropic: 2 dEVs}]
  \hookrightarrow [\text{isotropic: 1 EV}] .
\end{equation}
\end{linenomath}
Here ``orthotropic: 3 dEVs'' means that there are 3 distinct eigenvalues, and the
symmetry group for the orthotropic class---\rvsi{the symmetries of a general right rectangular 
prism in 3D}---is a subgroup of the symmetry group for the plan isotropic 
class---\rvsi{the symmetries of a circular cylinder in 3D}---which has only 2 distinct 
eigenvalues, and similarly for the other cases.  This means that it should be possible 
e.g.\ that each realisation is orthotropic, but the mean is isotropic---\rvsi{and thus
has the symmetries of a sphere in 3D}---or any other
such combination up the inclusion chain in \feq{Eq:3D-symm-Hasse}, such that the
mean has the same or a larger symmetry group.
The case in 2D is even simpler, the symmetry subgroup relations are
\begin{linenomath}
\begin{equation}   \label{Eq:2D-symm-Hasse}
  [\text{orthotropic: 2 dEVs}]  \hookrightarrow [\text{isotropic: 1 EV}] ,
\end{equation}
\end{linenomath}
\rvsi{where the orthotropic class has the symmetries of a general rectangle in 2D,
and the isotropic class has the symmetries of a circle in 2D.}

\subsection{Point-wise modelling} \label{Ssec:group-inv-detl}
To show the different possibilities, % \rvsi{for second order spatial tensors},
this will be done in three steps: first  \rvsi{with fixed} 
spatial orientation only the eigenvalues (i.e.\ the scaling) will be random.
In the case of second order tensors, \rvsi{like the heat conduction tensor},
which is considered here, this determines the group invariance class by the number
of coinciding eigenvalues.  In a second step, the eigenvalues \rvsi{are} fixed but
the spatial orientation of the symmetry axes \rvsi{varies} randomly; and in the last step both
of sources of randomness will be modelled together.

\paragraph{Fixed spatial orientation:}
The modelling approach is easiest to see by modelling only the anisotropic scaling as
random, and keep a fixed spatial orientation of the anisotropy.
\begin{linenomath}
\begin{multline} \label{Eq:SpecDec-rep}
  \genSecTensor(\vek{x},\event) = \meaneigvect(\vek{x})   \meaneigval^{1/2}(\vek{x}) 
  \wsymb{\eigval}(\vek{x},\event) \meaneigval^{1/2}(\vek{x}) \meaneigvect^T(\vek{x}) \\
  = \meaneigvect(\vek{x})   \meaneigval^{1/2}(\vek{x}) \exp(\wsymb{\symdiag}(\vek{x},\event)) 
  \meaneigval^{1/2}(\vek{x}) \meaneigvect^T(\vek{x}) .
\end{multline} 
\end{linenomath}
With this scaling $\meaneigval$, which together with the fixed orientation $\meaneigvect$ and
the random element $\wsymb{\eigval}(\event)= \exp(\wsymb{\symdiag}(\event))$
determine the symmetry class for each realisation corresponding to the group $G$, it can be
arranged that %$\Expvalue_{\vartheta_F}[\wsymb{\symdiag}(\vek{x},\cdot)] = \vek{0}$, and thus
\rvsi{
$
  \Expvalue_{\vartheta_S}[\wsymb{\eigval}] =  \Expvalue_{\vartheta_P}[\wsymb{\eigval}] = 
  \Expvalue_{\vartheta_L}[\wsymb{\eigval}] =  \vek{I}
$
for the scaling mean based on the scaling-rotation resp.\ scaling
distance $\vartheta_S$ (\feq{Eq:scal-dist} in Appendix~\ref{App:metrics-SPD}) preferred here, 
as well as the Lie product metric $\vartheta_P$ (\feq{Eq:prod-Lie-dist} on which it is based, 
and in this special case even for the log-Euclidean metric $\vartheta_L$ (\feq{Eq:log-Euclid-dist}.
}
%for both for the log-Euclidean resp.\ logarithmic mean based on the log-Euclidean distance metric
%$\vartheta_L$ (\feq{Eq:log-Euclid-dist} in Appendix~\ref{App:metrics-SPD}),
%as well as for the scaling mean based on the scaling-rotation resp.\ scaling
%distance $\vartheta_S$ (\feq{Eq:scal-dist} in Appendix~\ref{App:metrics-SPD}), 
%which in turn is based on the product
%Lie group distance $\vartheta_P$ (\feq{Eq:prod-Lie-dist} in Appendix~\ref{App:metrics-SPD}).
Hence one has from \feq{Eq:SpecDec-rep}
\begin{linenomath}
\begin{multline} \label{Eq:ExpSpecDec-rep}
  \Expvalue_{\vartheta_S}[\genSecTensor(\vek{x},\cdot)] = 
  \Expvalue_{\vartheta_P}[\genSecTensor(\vek{x},\cdot)] = 
  \Expvalue_{\vartheta_L}[\genSecTensor(\vek{x},\cdot)] \\
   = \meaneigvect(\vek{x})   \meaneigval^{1/2}(\vek{x}) \;
  \vek{I} \; \meaneigval^{1/2}(\vek{x}) \meaneigvect^T(\vek{x})  
  = \meaneigvect(\vek{x})   \meaneigval(\vek{x}) \meaneigvect^T(\vek{x}) = \meanC(\vek{x}) ,
\end{multline} 
\end{linenomath}
which one can arrange to have the invariance properties of the group $G_m\supseteq G$,
\rvsi{i.e.\ the desired number of coinciding eigenvalues in $\meaneigval(\vek{x})$.}.

The representation equation \feq{Eq:SpecDec-rep} may just as well be written for
$\symtensor(\vek{x},\event)$:
\begin{linenomath}
\begin{multline} \label{Eq:SpecDec-rep-log-gen}
  \symtensor(\vek{x},\event)   = \log \genSecTensor(\vek{x},\event)
  = \meaneigvect(\vek{x})  \log\left( \meaneigval^{1/2}(\vek{x}) \wsymb{\eigval}(\vek{x},\event)  
  \meaneigval^{1/2}(\vek{x}) \right) \meaneigvect^T(\vek{x}) \\
  = \meaneigvect(\vek{x})  \left(\log \wsymb{\eigval}(\vek{x},\event) + 
      \log \meaneigval(\vek{x})\right) \meaneigvect^T(\vek{x})
  = \meaneigvect(\vek{x})  (\wsymb{\symdiag}(\vek{x},\event) + \meanY(\vek{x})) \meaneigvect^T(\vek{x}) 
\end{multline} 
\end{linenomath}
as diagonal matrices commute, with $\meanY(\vek{x}) = \log \meaneigval(\vek{x})$. 
Thus one establishes a link to the exponential representation.  
From the last expression in \feq{Eq:SpecDec-rep-log-gen} one may observe 
that in the unconstrained linear vector space $\SymSecRankSpace$ of the log
the scaling is merely a shift of origin, 
and as the non-dimensional form in \feq{Eq:Scaling} typically just makes the formulas
look more complicated, we shall drop it in the notation 
and simply assume that this has been carried out in the background.  
The Lie algebra representation \feqs{Eq:repr-prod-Lie-exp}{Eq:repr-prod-Lie-exp-var} 
thus can be simplified from \feq{Eq:SpecDec-rep} ---\rvsi{using $\symdiag_s(\vek{x},\event) =
\diag(y^{(j)}_s(\vek{x},\event)) = \wsymb{\symdiag}(\vek{x},\event) + \meanY(\vek{x})$} ---to
\begin{linenomath}
\begin{equation} \label{Eq:SpecDec-rep-spl}
  \genSecTensor_s(\vek{x},\event) = \meaneigvect(\vek{x})  \eigval_s(\vek{x},\event)  
  \meaneigvect^T(\vek{x}) 
  = \meaneigvect(\vek{x})   \exp(\symdiag_s(\vek{x},\event)) \meaneigvect^T(\vek{x}) ,
\end{equation} 
\end{linenomath}
where an index ``s'' has been attached to signal that only random anisotropic scaling,
i.e.\ a random $\symdiag_s$ was used.  Additionally, one has to require that 
$\Expvalue_{\vartheta_F}[\symdiag_s(\vek{x},\cdot)] = \meanY(\vek{x}) = 
\log \meaneigval(\vek{x})$, \rvsi{i.e.\ for all $j$: $y^{(j)}_s(\vek{x},\event) = 
\wsymb{y}^{(j)}(\vek{x},\event) + \bar{y}^{(j)}(\vek{x})$ and $\Expvalue[\wsymb{y}^{(j)}]=0$.}
Here we are using the Euclidean mean in the unconstrained vector space  $\SymSecRankSpace$, which
corresponds to the log-Euclidean, \rvsi{Lie product,} and scaling means 
of $\eigval_s(\vek{x},\event)$ ---and thus
$ \Expvalue_{\vartheta_L}[\eigval_s] = \rvsi{\Expvalue_{\vartheta_P}[\eigval_s] =} 
    \Expvalue_{\vartheta_S}[\eigval_s] = \meaneigval .
$ 
Obviously there is an analogous equation to \feq{Eq:SpecDec-rep-spl} 
involving the reference tensor $\meanC(\vek{x})$, 
i.e.\ from \feq{Eq:ExpSpecDec-rep} one has:
\begin{linenomath}
\begin{equation} \label{Eq:SpecDec-ref}
  \meanC  = \exp(\meanH) = \meaneigvect  \meaneigval \meaneigvect^T  
  \Longleftrightarrow \quad  \meanH = \log(\meanC) = 
\meaneigvect \meanY  \meaneigvect^T.
\end{equation} 
\end{linenomath}

\rvsi{In case it is only known that the $y^{(j)}(\vek{x},\event)$ have finite
variance for each $\vek{x}$, the maximum entropy distribution \cite{Jaynes-2003} for the scalar
random fields $\wsymb{y}^{(j)}(\vek{x},\event)$ are Gaussian scalar random fields 
$\wsymb{y}^{(j)}(\vek{x},\event)=\gamma^{(j)}(\vek{x},\event)$, possibly correlated.
This would result in a log-normal distribution of $\eigval(\vek{x},\event)$ and thus for
$\genSecTensor_s(\vek{x},\event)$.  Such a distribution is unbounded, both for
$\wsymb{y}^{(j)}$ and $1/\wsymb{y}^{(j)}$, i.e.\ both for $\eigval(\vek{x},\event)$
and its inverse $\eigval^{-1}(\vek{x},\event)$,  and hence is unbounded for
$\genSecTensor_s(\vek{x},\event)$ as well as for its inverse $\genSecTensor^{-1}_s(\vek{x},\event)$.
This poses some difficulties in the formulation of the stochastic boundary value problem (SBVP)
\feq{Eq:deterministic}, or its concrete example of heat conduction \feq{Eq:deterministic-thm}.
In \fsec{Ssec:spatial-model} this will be discussed further.}

%Further, to model a specific invariance or material symmetry class in the mean and each realisation,
%%where the anisotropy axes may be random as well,
%the correlation or dependence of the scaling random variables in %\feq{Eq:rndscl} 
%has to be mathematically described for each type of symmetry class.
%  This is described in the next \fsec{Ssec:spatial-model}.
%For second order tensors, 
%there are three types of invariance classes in 3D ($d=3$), and two in 2D ($d=2$)
%\cite{malyarenko_tensor-valued_2019}, depending on how many eigenvalues coincide.  

Let us first consider the simpler case of 2D in \feq{Eq:2D-symm-Hasse}, 
where the tensor can be isotropic with two 
equal eigenvalues $y_s^{(1)} = y_s^{(2)}$ ---a multiple of the identity---or orthotropic 
with distinct eigenvalues 
$y_s^{(1)} \neq y_s^{(2)}$, in which case the orthogonal matrix $\meaneigvect_s$ determines
the directions of the major axes of the tensor.  If in the first case one demands that each
realisation is isotropic---and obviously also the reference $\meanC$ used in \feq{Eq:Scaling}
---then only one random variable is needed to model both eigenvalues.  On the other hand, if 
only the mean has to be isotropic but individual realisations may be orthotropic, then 
two identically distributed  random variables could be used for the eigenvalues.  
Turning now to the case where the reference $\meanC$ is orthotropic, obviously one
eigenvalue logarithm has to be larger than the other, say $\bar{y}_s^{(1)} > \bar{y}_s^{(2)}$.  
For the individual realisations, one may now choose to require that this relation holds
also for each realisation, or only in the mean, together with some possible correlation of
the eigenvalues.  In case one wants to require for each realisation
$y_s^{(1)} > y_s^{(2)}$, one trick to achieve this
is to take two unconstrained random variables, say $\xi_1$ and $\xi_2$, and set $y_s^{(2)} = \xi_1$ 
and $y_s^{(1)} = \xi_1 + \exp(\xi_2)$.  Hence, in case individual realisations may be orthotropic,
two log-stiffness like parameters  $y_s^{(1)}, y_s^{(2)}$ are needed, and in case the
individual realisations are always isotropic, just one log-stiffness like parameters  $y_s^{(1)}$
is sufficient.

After this somewhat broad discussion of the 2D case, the 3D case in \feq{Eq:3D-symm-Hasse}
may be treated in a similar way.  The three invariance
classes are determined again by the multiplicity of the eigenvalues, and the three cases
are isotropy---all eigenvalues equal $y_s^{(1)} = y_s^{(2)} = y_s^{(3)}$ ---or plan isotropy, 
where two eigenvalues are equal (defining the plane of isotropy) and the third eigenvalue
differs from these, or finally orthotropy, where all three eigenvalues are different.
So it is possible to require that each realisation is plan-isotropic (or orthotropic), but that the
mean is isotropic, or that each realisation may be orthotropic but the mean plan isotropic
(or isotropic).  \rvsi{This is often invoked when modelling random mixtures of anisotropic
crystals \cite{Cowin2013}.}  Thus, depending on which case one wants to model, from one up to three
log-stiffness like parameters  $y_s^{(j)}$ are needed.

%The Lie algebra representation \feq{Eq:SpecDec-rep-spl} is thus equivalent in expressive power
%to the exponential representation based on \feq{Eq:GenKLE} resp.\ the representation
%\feq{Eq:Soize-sq-1}. 
%One sees clearly that in this case of fixed anisotropy axes $\meaneigvect$, in the case
%of second order tensors $\genSecTensor(\vek{x},\event)$ only the eigenvalues $\eigval(\vek{x},\event)$
%are random.

\paragraph{Random spatial orientation with fixed scaling:}
Nevertheless, it is a well known fact that anisotropic materials often also exhibit 
uncertainty in orientations due to the randomly oriented fibres/grains/crystals. 
The general procedure for this was pointed out before in \feq{Eq:rand-ax-rep}.
The directional uncertainty in $\symtensor(\event)$---by keeping the material 
symmetry constant---is incorporated by subjecting the eigenvectors in $\meaneigvect$ in 
as in \feq{Eq:rand-ax-rep} to random rotations.

%To address the topic of a possible random orientation of anisotropy axes,
% we present here a procedure
%which works not only with the Lie group representation \feqs{Eq:repr-prod-Lie-exp}{Eq:SpecDec-rep-spl},
%but also with the exponential representation \feqs{Eq:exp-gen-def}{Eq:GenKLE}, as well as the
%square-exponential representation \feeqs{Eq:Soize-sq-1}{Eq:Soize-exp}.

So for the start we keep the eigenvalues constant, and only subject the eigenvectors
to random fluctuations; we use a subscript ``r'' to signify that only the eigenvectors are random. 
The stochastic tensor $\symtensor_r(\event_r)\in \SymSecRankSpace$---with random 
orientations only---is then analogously to \feq{Eq:rand-ax-rep} defined in the form:
\begin{linenomath}
\begin{align}\label{Eq:stomodrot}
\symtensor_r(\event_r) &= \rotmat(\event_r) \meanH_r \rotmat(\event_r)^T =
   \rotmat(\event_r) \meaneigvect_r  \meanY_r  \meaneigvect_r^T \rotmat(\event_r)^T; \\
     \label{Eq:stomodrot-exp}
   \genSecTensor_r(\event_r) &= \exp(\symtensor_r(\event_r)) = 
   \rotmat(\event_r)\meanC\rotmat(\event_r)^T,
\end{align}
\end{linenomath}
%---or equivalently 
%$\genSecTensor_r(\event_r) = \rotmat(\event_r) \exp(\meanH_r) \rotmat(\event_r)^T
%= \rotmat(\event_r) \meanC_r \rotmat(\event_r)^T$, 
where  
$\rotmat(\event_r) \in \mrm{SO}(d)$ is a random rotation matrix %as in \feq{Eq:Exp-rand-ax-rotation}, 
chosen such that each realisation
of $\symtensor_r(\event_r)$ resp.\ $\genSecTensor_r(\event_r)$ satisfies the appropriate
invariance requirements.  We also want
\begin{linenomath}
\begin{equation}  \label{Eq:Exp-rand-rot}
   \Expvalue_{\vartheta_S}[\genSecTensor] =
   \Expvalue_{\vartheta_p}[\genSecTensor] =  \meanC = \Expvalue_{\vartheta_R}[\rotmat]\;
    \meanC \; \Expvalue_{\vartheta_R}[\rotmat^T] 
   = \vek{I} \;  \meanC \vek{I}  ,
\end{equation}
\end{linenomath}
which requires $\Expvalue_{\vartheta_R}[\rotmat] = \vek{I}$.

It now becomes necessary to represent random rotations like $\rotmat(\event_r)\in \mrm{SO}(d)$
in \feqs{Eq:stomodrot}{Eq:stomodrot-exp}.
There are very many approaches to represent rotations, as well as a large amount
of easily accessible literature; this will not be reviewed presently as this would lead the 
discussion astray.  In order to arrive at a vector space setting preferred by us,
we resort to the well known correspondence  \feq{Eq:Exp-rand-ax-rotation} between
the Lie group $\mrm{SO}(d)$ of orthogonal matrices with unit determinant and its Lie algebra 
$\F{so}(d)$ of skew-symmetric matrices \cite{Mezzadri2007, Moakher2002, cardoso_2010},
also alluded to in Appendix~\ref{App:metrics-SPD}.
\begin{linenomath}
\begin{equation}  \label{Eq:Exp-rand-ax-rotation}
  \rotmat = \exp(\skewsym)\in \mrm{SO}(d), \quad\text{ and } 
     \quad \skewsym = \log(\rotmat) \in \F{so}(d).
\end{equation} 
\end{linenomath}
As one wants that $\Expvalue_{\vartheta_R}[\rotmat(\vek{x},\cdot)] = \vek{I}$, this means
that one needs $\Expvalue_{\vartheta_F}[\skewsym(\vek{x},\cdot)] = \vek{0}$ in $\F{so}(d)$,
where the Lie algebra $\F{so}(d)$ is an unconstrained %$d(d-1)/2$-dimensional 
linear vector space.
Any distribution of $\skewsym(\vek{x},\event)$ symmetric to the origin $\vek{0}\in \F{so}(d)$ 
will produce the desired result.

%Of course, it will not be possible to establish a one-to-one relation between the 
%Lie group $\mrm{SO}(d)$ and some real vector space, as $\mrm{SO}(d)$ is compact---this
%is easily seen in 2D where $\mrm{SO}(2)$ is homeomorphic to the one one-sphere, i.e.\ the unit
%circle in $\mathbb{R}^2$, or in 3D where the representation through unit quaternions
%shows that $\mrm{SO}(3)$ is homeomorphic to the three-sphere, i.e.\ the unit
%sphere in $\mathbb{R}^4$ ---and a real vector space is not compact.  But the relations in
%\feq{Eq:Exp-rand-ax-rotation} come pretty close to this desideratum, which makes the probabilistic
%modelling quite a bit easier.
\rvsi{Of course, a continuous one-to-one map between the compact
Lie group $\mrm{SO}(d)$ and some real vector space it is not possible,
as a real vector space is not compact.  The compactness of the special rotation
groups is obvious in 2D, where $\mrm{SO}(2)$ is homeomorphic to the unit
circle in $\mathbb{R}^2$, or in 3D, where the unit quaternion representation
shows that $\mrm{SO}(3)$ is homeomorphic to the sphere in $\mathbb{R}^4$.  
The elements in the Lie algebra $\F{so}(d)$ may be seen
as generalised angles, they appear in the arguments of trigonometric function
when mapped onto $\mrm{SO}(d)$ by the exponential map, cf.\ \feq{Eq:Rodrigues}.
In \cite{Grigoriu2016} a similar approach is taken, where the rotations are represented
by angular parameters.}

\rvsi{Coming back to the probability distributions on the Lie algebra $\F{so}(d)$ which
then are mapped onto the compact Lie group $\mrm{SO}(d)$, one has to remark that even
if they are unbounded, they are wrapped on the compact Lie group to produce so-called
circular statistics, e.g.\ see \cite{jammalamadaka_topics_2001} ---just think of the 2D
case, where the Lie algebra $\F{so}(2)$ is the real line, which is mapped via $\alpha
\mapsto (\cos(\alpha), \sin(\alpha))$ onto the circle.  
Thus probability distribution of the rotation matrix in $\mrm{SO}(d)$ is always bounded.}

Especially in 3D one may use an explicit version of \feq{Eq:Exp-rand-ax-rotation},
the well-known Rodrigues rotation formula, e.g.\ \cite{cardoso_2010}. 
There is a familiar correspondence between a skew-symmetric matrix $\skewsym\in\F{so}(3)$ and
the defining Euler vector $\axis = [w_1,w_2,w_3]^T \in \threeD$ along its unit 
rotation axis vector $\axis / \|\axis\|$
\begin{linenomath}
\begin{equation}  \label{Eq:def-skew-mat}
	\skewsym = \begin{bmatrix}
	0 & -w_3 & \phantom{-}w_2 \\ \phantom{-}w_3 & 0 & -w_1 \\ -w_2 & \phantom{-}w_1 & 0
	\end{bmatrix},
\end{equation}
\end{linenomath}
which is the unique rotation axis of a non-trivial rotation $\rotmat$ with
rotation angle $\rotangle=  \|\axis\|$.  Rodrigues's rotation formula
may then, thanks to the Cayley-Hamilton theorem, be written as
\begin{linenomath}
\begin{equation}\label{Eq:Rodrigues}
	\rotmat = \exp(\skewsym) = \vek{I} + \frac{\sin\rotangle}{\rotangle}\skewsym + 
	 \frac{1-\cos\rotangle}{\rotangle^2} \skewsym^2 =
	 \vek{I} + \frac{\sin\rotangle}{\rotangle}\skewsym + 
	 2\, \frac{\sin^2 (\rotangle/2)}{\rotangle^2} \skewsym^2,
\end{equation}
\end{linenomath}
and inversely \cite{Engoe2001}
\begin{linenomath}
\begin{equation}\label{Eq:Rodrigues-inv}
\skewsym = \frac{\arcsin \alpha}{\alpha}\vek{S}\quad\text{ with } 
   \vek{S} = \frac{1}{2} (\rotmat - \rotmat^T)\quad\text{ and }
   \alpha = \sqrt{\frac{1}{2}\, \tr(\vek{S}\vek{S}^T)}.
\end{equation}
\end{linenomath}
The defining vector $\axis(\event_r)$ is a $\threeD$-valued
random vector with values in the ball with radius $\uppi$, and the random angle of rotation is
$\rotangle(\event_r) = \|\axis(\event_r)\|$, defined on a probability space 
$(\samplespace_r, \sigmaalgebra_r,\probability_r)$.  
But in case each realisation is plan isotropic (this assumed to be the 1-2 plane), there is no point
in rotations with a $w_3$ component --- a rotation in the 1-2 plane.  In this instance one only
needs two angle-like parameters $w_1, w_2$.

Note that in the 2D case the axis of rotation $\axis = (0,0,\rotangle)$ is fixed
perpendicular to the plane; only the
rotational angle $\rotangle(\event_r)$ resp.\ the parameter $w_3$ is modelled as a circular/angular 
random variable \cite{mardia_directional_2000, jammalamadaka_topics_2001}.
And in case each realisation is isotropic, no rotation at all is required, as it will have no effect.

In any case one obtains a random $\rotmat(\event_r)\in\mrm{SO}(d)$, determined by 0--3
angle-like parameters $w_j$, defining the axis of rotation $\axis$.  Again, as long as the
distribution of $\axis(\event_r)$ resp.\ $\rotangle(\event_r)$ is symmetric around the
origin, % or rather $\Expvalue(\axis) = \vek{0}$ resp.\ $\Expvalue(\rotangle) = 0$, 
the expected value of $\rotmat(\event_r)$ of the here adequate Fréchet or Karcher
mean with the metric $\vartheta_R$ (defined in Appendix~\ref{App:metrics-SPD}) is the identity, 
$\vek{I}=\Expvalue_{\vartheta_R}(\rotmat)$.  Hence \feq{Eq:Exp-rand-rot} holds in this case.

\paragraph{Random spatial orientation and scaling:}
To account for both random orientation and scaling attributes in the tensor $\symtensor$,
one may combine the two approaches.  If these two aspects are considered independent,
one may define the probability space as a product space $(\samplespace,\sigmaalgebra,\probability)$ 
with  $\samplespace:=\samplespace_s\times \samplespace_r$, 
$\sigmaalgebra:=\sigmaalgebra_s\otimes \sigmaalgebra_r$, and a product measure
$\probability:=\probability_s\probability_r$, \rvsi{but as pointed out in \cite{Grigoriu2016},
often they will be not independent, as they both result from the probability distribution of 
$\genSecTensor$.}    In any case, denoting a combined
rotational-scaling uncertainty with a subscript ``rs'', and reusing the already defined subscripts
``s'' and ``r'',
let $\genSecTensor_s(\vek{x},\event_s)$ denote
a representation of an SPD tensor field   with given group invariance
and a fixed  orientation of anisotropy axes, and $\Expvalue_{\vartheta_L}[\genSecTensor_s(\vek{x},\cdot)]
= \meanC(\vek{x})$, like the one just mentioned in \feq{Eq:SpecDec-rep-spl},
and let $\rotmat(\vek{x},\event_r)
\in \mrm{SO}(d)$ be a random orientation %with $\Expvalue_{\vartheta_R}[\rotmat(\vek{x},\cdot)] = \vek{I}$
as in \feeqs{Eq:stomodrot}{Eq:Exp-rand-ax-rotation}.
%---for the definition of the `rotation distance' $\vartheta_R$ on the Lie group $\mrm{SO}(d)$ see the
%Appendix~\ref{Asec:Lie-grp-dist}.   
\begin{linenomath}
\begin{equation}  \label{Eq:rand-ax-rep}
   \genSecTensor_{rs}(\event_{rs}) = \rotmat(\event_r) \genSecTensor_s(\event_s) 
                 \rotmat(\event_r)^T; \text{ resp.} \quad
   \symtensor_{rs}(\event_{rs}) = \rotmat(\event_r) \symtensor_{s}(\event_{s})  \rotmat(\event_r)^T,
\end{equation}
\end{linenomath}
is then a representation of an SPD  tensor field $\genSecTensor(\vek{x},\event_{rs})$  with given 
group invariance of $\genSecTensor_s(\vek{x},\event_s)$ and random anistropic scaling, 
but in addition
with a random orientation of anisotropy axes---with the random orientation described by the 
$\mrm{SO}(d)$-valued random field $\rotmat(\vek{x},\event_r)$.  
One then has for the scaling mean, based on the scaling
distance $\vartheta_S$ (\feq{Eq:scal-dist} in Appendix~\ref{App:metrics-SPD}):
\begin{linenomath}
\begin{equation}  \label{Eq:Exp-rand-ax-rep}
   \Expvalue_{\vartheta_S}[\genSecTensor_{rs}] = 
   \Expvalue_{\vartheta_P}[\genSecTensor_{rs}] = \Expvalue_{\vartheta_R}[\rotmat]\;
   \Expvalue_{\vartheta_L}[\genSecTensor_s] \; \Expvalue_{\vartheta_R}[\rotmat^T] 
   = \vek{I} \; \Expvalue_{\vartheta_L}[\genSecTensor_s]\; \vek{I} =  \meanC .
\end{equation}
\end{linenomath}
\rvsi{This again shows how well the Lie representation fits with the Fréchet mean based on the
metric $\vartheta_S$, resp.\ $\vartheta_P$, and which satisfies all
our desiderata in \fsec{Ssec:means-spd}.}

\rvsi{This is also the appropriate place to summarise the possible properties of probability
measures on the parameters spaces, and their effects on the probability distribution
of  the tensor $\genSecTensor$, as well as on the computation of the Fréchet means
$\Expvalue_{\vartheta_S}$ and $\Expvalue_{\vartheta_P}$.}

\rvsi{Looking at \feqs{Eq:SpecDec-rep}{Eq:SpecDec-rep-log-gen}, \feqs{Eq:stomodrot}{Eq:stomodrot-exp},
 and \feq{Eq:rand-ax-rep}, one may express this in total as 
\begin{linenomath}
\begin{align}  \label{Eq:SpecDec-randC}
  \genSecTensor_{rs} &= \rotmat \meaneigvect  \meaneigval \wsymb{\eigval} 
   \meaneigvect^T \rotmat^T, \; \text{ with } \; \wsymb{\eigval} = \exp(\wsymb{\symdiag}),
   \; \rotmat = \exp(\skewsym)
  \\  \label{Eq:SpecDec-randH}
  \symtensor_{rs} &= \rotmat \meaneigvect  \meanY \wsymb{\symdiag} 
   \meaneigvect^T \rotmat^T,   \; \text{ and } \; \meaneigval = \exp(\meanY).
\end{align}
\end{linenomath}
The random parameters are $(\wsymb{\symdiag},\skewsym) \in \F{diag}(d) \times \F{so}(d)$,
and from the previous discussion it is clear that in order to have \feq{Eq:Exp-rand-ax-rep},
one requires 
\begin{linenomath}
\[
\Expvalue_{\vartheta_S}[\rotmat] = \Expvalue_{\vartheta_R}[\rotmat]=\vek{I}\; \text{ and } \;
\Expvalue_{\vartheta_S}[\wsymb{\eigval}] = \Expvalue_{\vartheta_L}[\wsymb{\eigval}] = \vek{I}.
\]
\end{linenomath}
As the actual computation of the general Fréchet means is more involved than the computation
of the normal arithmetic resp.\ Euclidean mean, it is worthwhile to point out once more
that this will be achieved if both $\Expvalue_{\vartheta_F}[\skewsym] = \vek{0}$ and
$\Expvalue_{\vartheta_F}[\wsymb{\symdiag}] = \vek{0}$, and additionally
the distributions are symmetric w.r.t.\ the origin.  This symmetry w.r.t.\ the origin 
of the logarithms of the perturbation of the mean, $\wsymb{\eigval}$, additionally
makes it quite easy to reflect the fact that one typically knows as much about $\symtensor_{rs}$
as about its inverse $\symtensor_{rs}^{-1}$, cf.\ the desiderata in \fsec{Ssec:means-spd}.}

\rvsi{It also may be gleaned from
the descriptions in Appendix~\ref{App:metrics-SPD} that if the distributions of
$\rotmat$ and $\wsymb{\eigval}$ do not deviate too much from the identity $\vek{I}$,
or, equivalently, the the distributions of
$\skewsym$ and $\wsymb{\symdiag}$ do not deviate too much from the origin $\vek{0}$,
the scaling metric $\vartheta_S$ and the Lie product metric $\vartheta_P$ actually
coincide, and thus also their respective means.  This makes the computation of
$\Expvalue_{\vartheta_S}$ much easier.
}

\rvsi{It was also already remarked, that whatever the distribution of the angle-like parameters
$\skewsym) \in  \F{so}(d)$, the rotation matrix $\rotmat  = \exp(\skewsym)$ is always
bounded---it lives on the compact set $\mrm{SO}(d)$---and hence has moments of any order.
This means that the question of whether $\symtensor_{rs}$ is bounded or has moments of
some order is completely determined by the scaling $\wsymb{\eigval} = \exp(\wsymb{\symdiag})$.
In particular, if the distribution of $\wsymb{\symdiag}$ is bounded, so are the distributions
of $\wsymb{\eigval}$ and $\wsymb{\eigval}^{-1} = \exp(-\wsymb{\symdiag})$, and thus also
those of $\symtensor_{rs}$ and $\symtensor_{rs}^{-1}$.  Unfortunately, the often used
centred normal distribution for $\wsymb{\symdiag}$ is not bounded---but has moments of any order---
and thus the distributions of both $\symtensor_{rs}$ and $\symtensor_{rs}^{-1}$ are unbounded,
although of second order.  We will comment on this specific case in the next 
\fsec{Ssec:spatial-model}.}

\subsection{Spatial modelling} \label{Ssec:spatial-model}
According to \feqs{Eq:repr-prod-Lie-exp}{Eq:repr-prod-Lie-exp-var} the random SPD matrix at
a point $\vek{x} \in \spatialdomain$ depends on the log-eigenvalues in 
$\symdiag(\vek{x},\event) = \log \eigval(\vek{x},\event)$ which are 1--3 free 
log-scale like parameters $\tilde{\vek{y}} = [\dots,y_j,\dots]^T$, and on the skew 
matrix $\skewsym(\vek{x},\event) = \log \eigvect(\vek{x},\event)$ determined 
by 0--3 angle-like variables $\tilde{\vek{w}} =[\dots,w_j,\dots]^T$.  We
collect those parameters in the column vector 
$\tilde{\vek{z}} = [\tilde{\vek{y}}^T,\tilde{\vek{w}}^T]^T \in \Real^m$.
Similarly, one may collect the parameters for the mean matrix 
\feq{Eq:SpecDec-ref} in the column vector $\hat{\vek{z}}(\vek{x})$.

\paragraph{Correlation structure:}
The joint correlation structure of the random parameters $\vek{z}(\vek{x},\event)$
is described through a joint covariance matrix $\vek{J} \in \Real^{m\times m}$:
\begin{linenomath}
\begin{multline}  \label{Eq:JointCorr}
  \vek{J}(\vek{x}_1,\vek{x}_2)  =  \Expvalue_{\vartheta_E}\left( 
  [\tilde{\vek{y}}(\vek{x}_1,\cdot), \tilde{\vek{w}}(\vek{x}_1,\cdot)]\otimes  
  [\tilde{\vek{y}}(\vek{x}_2,\cdot), \tilde{\vek{w}}(\vek{x}_2,\cdot)] \right)  \\ 
  =  \Expvalue_{\vartheta_E}  \begin{pmatrix} 
     \tilde{\vek{y}}(\vek{x}_1,\cdot) \otimes \tilde{\vek{y}}(\vek{x}_2,\cdot) &
     \tilde{\vek{y}}(\vek{x}_1,\cdot) \otimes \tilde{\vek{w}}(\vek{x}_2,\cdot) \\
     \tilde{\vek{w}}(\vek{x}_1,\cdot) \otimes \tilde{\vek{y}}(\vek{x}_2,\cdot) &
     \tilde{\vek{w}}(\vek{x}_1,\cdot) \otimes \tilde{\vek{w}}(\vek{x}_2,\cdot) 
     \end{pmatrix}  .
\end{multline}
\end{linenomath}
%but due to the constraint $\vek{K} \vek{z} = 0$ we require  
%$\vek{K}\widetilde{\vek{J}}(\vek{x}_1,\vek{x}_2)\vek{K}^T = 0$.  If one had random fields
%$\vek{\eta}(\vek{x},\event) \in \Real^{m-r}$, setting 
%$\vek{z}(\vek{x},\event) = \vek{B}\vek{\eta}(\vek{x},\event)$ and then computing
%$\widetilde{\vek{J}}(\vek{x}_1,\vek{x}_2)$ as in \feq{Eq:JointCorr}, it is seen that the requirement
%$\vek{K}\widetilde{\vek{J}}(\vek{x}_1,\vek{x}_2)\vek{K}^T = 0$ is automatically satisfied.
%
%A short calculation shows that the joint correlation structure 
%of the fields $\vek{\eta}(\vek{x},\event)$
%has to be related to $\widetilde{\vek{J}}(\vek{x}_1,\vek{x}_2)$ from \feq{Eq:JointCorr} through
%\begin{equation}  \label{Eq:JointCorr1}
%  \vek{J}(\vek{x}_1,\vek{x}_2)  =  \Expvalue_{\vartheta_E}\left( 
%  \vek{\eta}(\vek{x}_1,\cdot)\otimes  \vek{\eta}(\vek{x}_2,\cdot) \right) 
%  = \vek{B}^T \; \tilde{\vek{J}}(\vek{x}_1,\vek{x}_2) \vek{B} .
%\end{equation}
\rvsi{A good description of heterogeneous materials and their statistical description
may be found in \cite{torquato:2001} and references therein.  
An often employed family of covariance functions, which allow one to separately control
the smoothness and correlation length of the random field, are the family of Matérn functions
(e.g.\ \cite{Matern_1986, Cressie1993, GneitingEtAl2010}).  These kind of functions also
allow efficient approximation  via low-rank tensor methods \cite{Litvinenko_Keyes2019}, 
important in the solution of the Fredholm integral equation \feq{Eq:KLE-eigprob}
and field expansion \feq{Eq:GenKLE}.
On the plus side of this homogeneous covariance is its versatility and that its spectrum is 
analytically known, unfortunately, neither this covariance nor its continuous spctrum
are compactly supported.  One example of a covariance function with compact support is
given in \cite{Gneiting2002}.  This renders the integral operator in \feq{Eq:KLE-eigprob}
compact, hence it has a discrete spectrum and a proper Karhunen-Loève series expansion
\feq{Eq:GenKLE}.}

\paragraph{Random field expansion:}
Aiming at a Karhunen-Loève expansion of a corresponding spatial field (see e.g.\ \cite{hgm07}), 
one has to solve the eigenproblem
\begin{linenomath}
\begin{equation}  \label{Eq:KLE-eigprob}   J(\vek{\eta}_\ell)(\vek{x}_1) :=
  \int_\spatialdomain \vek{J}(\vek{x}_1,\vek{x}_2) \, \vek{\eta}_\ell(\vek{x}_2) \;\di \vek{x}_2 = 
    \mu_\ell \, \vek{\eta}_\ell(\vek{x}_1) ,
\end{equation}
\end{linenomath}
for eigenvalues $\mu_\ell$ and eigenvector fields $\vek{\eta}_\ell(\vek{x}) \in \Real^m$.
 
It should be noted that spatial homogeneity would be evident as usual from the fact
that $\vek{J}(\vek{x}_1,\vek{x}_2)$ would be a function of only the difference of
the arguments, $\vek{x}_1 - \vek{x}_2$.  In that case the eigenproblem \feq{Eq:KLE-eigprob} 
is a convolution equation, and it is well known how to solve those by Fourier
methods, as indeed the Fourier functions are (possibly generalised) eigenfunction
of the operator $J$ represented by \feq{Eq:KLE-eigprob}.
In any case, the orthonormal eigenfunctions $\vek{\eta}_\ell(\vek{x})$ determine uncorrelated
mean zero unit variance
random variables $\xi_\ell(\event)$ through the projection of $\vek{z}(\vek{x},\event)$:
\begin{linenomath}
\begin{equation}  \label{Eq:KLE-RV-proj}
      \xi_\ell(\event) = \frac{1}{\sqrt{\mu_\ell}}\,
  \int_\spatialdomain \vek{z}(\vek{x},\event)^T  \, \vek{\eta}_\ell(\vek{x}) \;\di \vek{x} .
\end{equation}
\end{linenomath}
But normally we want to perform the reverse task, a spectral synthesis, and thus the random variables
$\xi_\ell(\event)$ have to be chosen such that $\vek{z}(\vek{x},\event)$ has the proper
statistics at each $\vek{x} \in \spatialdomain$.  This problem is outside the scope of this
work and has been addressed elsewhere,
see e.g.\ \cite{hgm07} and the literature therein.

Assuming the proper mean zero unit variance
random variables $\xi_\ell(\event)$ have been chosen, 
for the spectral synthesis resp.\ Karhunen-Loève expansion one typically has an
infinite series---or an integral in case of a continuous spectrum---over 
the spectrum $\mu_\ell \in \sigma(J)$ of the integral  covariance operator in \feq{Eq:KLE-eigprob} 
(see e.g.\ \cite{hgm07, DLiuEtAl2020}):
\begin{linenomath}
\begin{equation}  \label{Eq:GenKLE}
  \vek{z}(\vek{x},\event)  = \hat{\vek{z}}(\vek{x}) + \tilde{\vek{z}}(\vek{x},\event) = 
  \hat{\vek{z}}(\vek{x}) + 
    \sum_\ell  \xi_\ell(\event)  \vek{\eta}_\ell(\vek{x}).
\end{equation}
\end{linenomath}
%Here $\hat{\vek{\eta}}_m$ is a basis of $\Real^{m - s}$ so that 
%$\widehat{\vek{B}}\, \hat{\vek{\eta}}_m$ is a basis for the nullspace of 
%$\widehat{\vek{K}}$, and thus the first term $\bar{\vek{z}}(\vek{x})$ in \feq{Eq:GenKLE}
%represents the parameters for the mean matrix $\widehat{\vek{C}}$ with the
%proper invariance, ascertained through $\widehat{\vek{K}} \bar{\vek{z}}_m(\vek{x})=0$.  
%A similar argument shows that the second term  $\tilde{\vek{z}}(\vek{x})$  in 
%\feq{Eq:GenKLE} satisfies ${\vek{K}} \tilde{\vek{z}}_m(\vek{x})=0$ and thus represents
%parameters for the random part of the matrix $\genSecTensor(\vek{x},\event)$ with
%the required invariance.
It is of course also possible to use more general expansion functions $\vek{\eta}_\ell$,
e.g.\ tight frames \cite{LuschgyPages2009}, an approach which may have advantages especially
in the lognormal case \cite{BachmayrCohenDeVoreEtAl2017, BachmayrCohenMigliorati2018}.

Note that through the use of a Lie group mapping onto $\PosSymSecRankSpace$, we are
able to essentially do all the modelling in the Lie algebra, which is a free vector
space, and the expression \feq{Eq:GenKLE} is an example of this.  This brings the
modelling of random spatial fields or temporal stochastic processes of $\PosSymSecRankSpace$ tensors
to the same level of difficulty as the usual modelling of random fields resp.\ stochastic processes.

\paragraph{Elliptic stochastic boundary value problems (SPDEs) and their solution:}
\rvsi{As already alluded to in the motivation for this work in \fsec{Sec:Problem}, the final
use of the modelling is to be able to investigate and solve equations like \feq{Eq:deterministic},
here we will focus specifically on the linear elliptic partial differential 
equation (PDE) \feq{Eq:deterministic-thm}, 
e.g.\ a boundary value problem (BVP) for the case of stationary heat conduction, or more generally
a stationary diffusion problem, where the random tensor $\genSecTensor(\vek{x}) = \DFcond$ is the
heat conductivity resp.\ diffusion tensor.  In case this tensor is uncertain and modelled
probabilistically, this becomes a stochastic BVP (SBVP), or a stochastic PDE (SPDE).}

\rvsi{The solution space for the temperature $T(\vek{x})$ in \feq{Eq:deterministic-thm} 
in the simplest instance in the deterministic
case is a closed subspace of the Hilbert-Sobolev space $\Hp^1(\spatialdomain)$ incorporating
the boundary conditions, say for simplicity $T(x) \in \Ho^1(\spatialdomain)$, i.e.\ $T(\vek{x})$
vanishes on the boundary $\partial \spatialdomain$.  The right-hand side $f(\vek{x})$
then typically has to be in the dual space $f(\vek{x}) \in \Hp^{-1}(\spatialdomain)$.
Thus the heat conductivity resp.\ diffusion tensor $\DFcond$ 
has to be in $\Lp_\infty(\spatialdomain)$, i.e.\ bounded,
in order for the map $T \mapsto f$ to be continuous.  For a well-posed problem one would also
want the inverse of that map to exits and be continuous, and this requires also the inverse
diffusion tensor $\vek{\kappa}^{-1}(\vek{x}) \in \Lp_\infty(\spatialdomain)$, i.e.\ to be also
bounded.  This is well known standard elliptic theory thanks to the
open mapping theorem \cite{segalKunze78, yosida-fa-1980}.}

\rvsi{When considering its stochastic extension with a random $\vek{\kappa}(\vek{x},\event)$,
the simplest is to require the random solution temperature field $T(\vek{x},\event)$ to
have finite second moments, i.e.\ to
live in the tensor product Hilbert space $\Ho^1(\spatialdomain) \otimes \Lp_2(\samplespace)$,
and the right-hand side random sources or sinks have to be in its dual space
$f(\vek{x},\event) \in \Hp^{-1}(\spatialdomain) \otimes \Lp_2(\samplespace)$.
The formulation of a well-posed problem in this case requires
the random diffusion tensor field $\vek{\kappa}(\vek{x},\event) \in 
\Lp_\infty(\spatialdomain \times \samplespace)$ to be bounded, as well as its inverse
$\vek{\kappa}^{-1}(\vek{x},\event) \in \Lp_\infty(\spatialdomain \times \samplespace)$,
e.g.\ see \cite{Matthies2005, hgm07, Matthies2008, LordPowell2014, DungSchwabX2022}.  
As already mentioned in
\fsec{Ssec:group-inv-detl}, this requires the probability densities of
the eigenvalues $\eigval(\vek{x},\event)$ to be bounded, and bounded away from singularity,
i.e.\ $\eigval^{-1}(\vek{x},\event)$ also to be bounded, at each point $\vek{x}\in \spatialdomain$.
Both of these requirements can be satisfied simultaneously by choosing a bounded distribution for
the log-eigenvalues $\symdiag(\vek{x},\event)$.
One possibility is the quite versatile beta distribution, see e.g.\ \cite{Grigoriu2016}.}

\rvsi{So far, this is the standard case.  But, as already mentioned above, 
the maximum entropy distribution for the
log-eigenvalues $\symdiag(\vek{x},\event)$ with known finite second moments is the normal
or Gaussian distribution \cite{Jaynes-2003}, so this is a highly desired choice
\cite{schwartzman_lognormal_2016}.
This makes the eigenvalues $\eigval$ and hence $\genSecTensor = \vek{\kappa}$ 
lognormal, and the density is unbounded for both $\eigval$ and $\eigval^{-1}$.
This means that neither the map $T \mapsto f$ nor its inverse (the solution map) are
continuous on the solution space for the bounded case above.  But it is possible
also to proceed in this case, with different techniques, which either make the soultion
space smaller, or change the solution concept slightly; see e.g.\
\cite{EspigHackbuschEtAl2014, LordPowell2014, HoangSchwab2014,
BachmayrCohenDeVoreEtAl2017, BachmayrCohenMigliorati2018, HerrmannSchwab2019, DungSchwabX2022}
for more details.}

\rvsi{Such SBVPs resp.\ SPDEs will be considered and solved as examples in \fsec{Sec:Application},
and we will not further dwell on the numerical solution techniques, the interested reader
may find more material in e.g.\ \cite{Matthies2005, hgm07, Matthies2008, 
EspigHackbuschEtAl2014, HoangSchwab2014, grigoriu2012, LordPowell2014, HerrmannSchwab2019, 
DungSchwabX2022} and the references therein.  The numerical solution of such SPDEs,
and more concretely such example situation as depicted in \fsec{Sec:Application}, is
also necessary when one wants to identify the properties of the tensor field from measurements
resp.\ observations; a topic we turn to next.}

\paragraph{Identification of tensorial random fields:}
\rvsi{The identification of the properties, and especially the identification of the
symmetry classes of materials is a vast topic beyond the scope of this work, 
so here we will only give some remarks.
The identification of symmetry classes may be based e.g.\ on the knowledge of the
crystal structure of the material \cite{nye_physical_1984}, or on its fabrication
\cite{Cowin2013}.  Otherwise one must rely on measurements or observations.  Much has been
doen here in the case of elasticity classes \cite{CowinMehr1989, FrancoisEtAl1998, 
DanekEtAl2015, DanekSlawinski2015, GierlachDanek2018}, which are more involved than
for the case of second order tensors considered here.
Specifically for the kind of stochastic modelling proposed here some pointers are given in 
section 4.2 of \cite{Grigoriu2016}.}

\rvsi{One very general technique for such identification resp.\ inverse problems
are the Bayesian methods, see e.g.\
\cite{mosegaard:2002:IP, tarantola:2005, marzouk:2007:JCP, marzouk:2009:CompPhys, AMStuart2010, 
RosicLitvinenkoEtAl2012, PajonkRosicEtAl2012, RosicKucerovaSykoraEtAl2013, 
MatthiesZanderEtAl2016, MatthiesZanderRosicEtAl2016, RosicSykoraPajonkEtAl2016, hgm18}
and the references therein.
These are probabilistic methods, which have already been used for anisotropic materials, e.g.\
\cite{FittEtAl2019, SavvasPapai2020}.}

\rvsi{If the measurements or tests, like those computed in the examples in \fsec{Sec:Application}, 
are designed properly, so that they are informative about the quantities to be identified,
quite complicated material behaviour can be identified and calibrated, see e.g.\
\cite{RosicKucerovaSykoraEtAl2013, MatthiesZanderEtAl2016, MatthiesZanderRosicEtAl2016, RosicSykoraPajonkEtAl2016}, and also e.g.\
\cite{DobrillaEtAl2023, AIbrahimbEtAl_SN-2022, AI-HGM-EK-2020}
for highly nonlinear and irreversible behaviour.
Whether the tests or experiments are informative may be ted out exactly with such computations
as in  \fsec{Sec:Application}.  This means the whole test or measurement and updating may be
performed first virtually in the computer.  From such computations one usually gets a good
idea whether the measurement is informative as regards what has to be identified, and, in
case it is not very informative, how to change the test or measurement.}

%% file: application-y.tex
% !TEX root = ../paper_tensors.tex
% !TEX encoding = UTF-8 Unicode

\section{Application: heat conduction in a proximal femur}\label{Sec:Application}
%
%Let us consider a two dimensional geometry $\spatialdomain \subset \Euclideanspace$ with smooth Lipschitz boundary $\Gamma$. The following governing equations
%\begin{align}\label{Eq:DetlForm}
%-\textrm{div}(\conducttensor \cdot \nabla\temperature(x))   &= \heatsource(x),  \quad \forall x \in \spatialdomain, \nonumber \\
%\temperature(x) &= {\temperature}_0,  \quad \forall x \in \Gamma_D, \\
%(\conducttensor \cdot \nabla\temperature(x))\cdot\vek{n}(x) &= {g}(x), \quad \forall x \in \Gamma_N, \nonumber
%\end{align}
%describe the physical phenomenon of steady state heat transfer. The aim is to determine the temperature $\temperature(x) \in \mathbb{R}$ at a spatial point $x\in \spatialdomain$.
%In \feq{Eq:DetlForm}, $\Delta$ defines the Laplacian operator; $\conducttensor\in\PosSymSecRankSpacetwoD$ represents the thermal conductivity tensor; $\heatsource(x) \in \mathbb{R}$ is the heat source;  ${g}(x) \in \mathbb{R}$ represents the surface heat flux on Neumann boundary $\Gamma_N \subset \Gamma$, and $\vek{n}(x) \in \mathbb{R}^2$ is the outward unit normal to $\Gamma_N$. For simplification, a homogeneous temperature ${\temperature}_0 \in \mathbb{R}$ is considered on the Dirichlet boundary $\Gamma_D \subseteq \Gamma $. We may also further assume that $\Gamma_D \cap \Gamma_N = 0$. 

The application we have in mind is the steady-state heat conduction in bones, a highly \rvsn{heterogeneous and} anisotropic material. More specifically, we shall look at both a 2D and 3D model of the \rvsn{human} proximal femur, this is the upper part of the thigh bone close to the hip.
\rvsn{This example aims to enhance our understanding of bone tissue response to heat exposure, crucial for guiding surgical procedures such as bone drilling or cutting \cite{mediouni_optimal_2017,augustin_cortical_2012}, or to analyse the frictional heat dissipation that occurs within the total hip replacement e.g. at the head-cup interface \cite{uddin_frictional_2013, davidson_heat_1988}.} 

\rvsn{It is crucial to emphasize that the bone heat conduction example in this research should be regarded merely as a simplified model, not representing a physiologically accurate scenario. Therefore, the chosen boundary conditions and material parameters in \fsecs{Sec:2DResults}{Sec:3DResults} remain largely experimental. The primary objective of the authors is to} investigate the effects of different models of the anisotropy of the heat \rvsn{conductivity} tensor
on the heat conduction and temperature distribution in the femur.

Referring to the deterministic model of heat conduction defined in \feq{Eq:deterministic-thm}, 
\rvsn{for the sake of simplicity, only spatially constant conductivity tensors are considered in this study, so that $\conducttensor(\vek{x})\equiv\conducttensor\in\PosSymSecRankSpace$.}
\rvsn{Following this, } let the conductivity tensor $\conducttensor\in\PosSymSecRankSpace$ be modelled as
a random tensor as described
%on probability space $(\samplespace,\sigmaalgebra,\probability)$
%as an exponential mapping of the model explained 
in \feq{Eq:rand-ax-rep}. 
The scaling terms $\vek{\Lambda_s}=\diag(\lambda_s^{(i)})\in\PosSymSecRankSpace$ of the 
tensor $\conducttensor$ are modelled as positive log-normal random variables:
\begin{linenomath}
\begin{equation}\label{Eq:lognormal}
 \vek{\Lambda_s}(\event_s) = \diag(\lambda_s^{(i)}(\event_s)) = \exp(\symdiag_s(\event_s))=
 \diag(\exp(y_s^{(i)}(\event_s))), 
\end{equation}
\end{linenomath}
i.e.\ a non-linear transformation of the Gaussian random variables
\begin{linenomath}
\begin{equation}
    y_s^{(i)}(\event_s) \sim \normaldist(\mu_i,\sigma_i),  
      \quad i= 1, \dots, d,
\end{equation}
\end{linenomath}
% \begin{equation}\label{Eq:lognormal}
% \textrm{diag}(\lambda_s(\event_s)^{(i)}) = \text{exp}({\mu_i} + \sigma_i\theta_i(\event_s)),  \quad i=\left\lbrace 1,..,d\right\rbrace,
% \end{equation}
% i.e., a non-linear transformation of standard Gaussian random variables $[\theta_i(\event_s)]^d_{i=1}$. 
where the parameters $\{\mu_i\}^d_{i=1}\in\mathbb{R}$ and $\{\sigma_i\}^d_{i=1}\in\PosRealStrict$ 
denote the mean and standard deviation corresponding to the Gaussian distribution.
% \begin{equation}
% 	 \normaldist(\mu_i,\sigma_i) \sim \log\eigvalsclele(\event_s), \quad i=\left\lbrace 1,..,d\right\rbrace.
% \end{equation}

Further, the directional uncertainty of $\conducttensor$ is accounted by modelling the 
eigenvectors in $\meaneigvect_r$ as von Mises-Fisher distributions defined on the unit sphere 
$\unitsphere$ \cite{mardia_directional_2000}. If $\rotvect$ is one of the eigenvectors of the matrix $\meaneigvect_r$, then a von Mises-Fisher unit vector $\rotvect(\event_r)$ with mean direction $\mucirc\in\unitsphere$ 
has a \acrfull{pdf} of the form:
%A von-Mises random variable $\rotangle(\event_r)$ with mean direction $\circmeananly=0$ has the \acrfull{pdf} of the form:
\begin{linenomath}
\begin{equation}\label{Eq:VMFPDF}
	f(\rotvect(\event_r)|\mucirc,\concparam) = 
	\frac{\concparam^{d/2-1}}{(2\pi)^{d/2}\mathcal{I}_{d/2-1}(\concparam)} 
	    \exp(\concparam \mucirc^T\rotvect(\event_r)), 
\end{equation}
\end{linenomath}
where 
%$\circmeananly\in[0,2\pi)$ (note that, in this study, $\circmeananly=0$) and 
$\concparam\in\PosRealStrict$ is the measure of concentration, and $\mathcal{I}_p$ is the 
modified Bessel function of order $p$.
When $d=2$, as only the
rotational angle $\rotangle(\event_r)$ is modelled as a random variable on the unit circle $\unitcircle$, the above PDF reduces to the 
von Mises distribution \cite{best_efficient_1979}:
\begin{linenomath}
\begin{equation}\label{Eq:VMPDF}
f(\rotangle(\event_r)|\circmeananly,\concparam) = 
 \frac{\exp(\concparam\cos(\rotangle(\event_r)-\circmeananly))}{2\pi \mathcal{I}_0(\concparam)},
\end{equation}
\end{linenomath}
in which $\circmeananly$ is the mean direction of the von Mises circular random variable 
$\rotangle(\event_r)$.
    
By rewriting the deterministic model in \feq{Eq:deterministic-thm} to a stochastic one 
\begin{linenomath}
\begin{equation}\label{Eq:stochastic-thm}
-\nabla\cdot(\conducttensor(\event) \nabla T(\vek{x},\event)) = 
      f(\vek{x})\quad \forall \vek{x}\in\spatialdomain, \event\in\samplespace,
\end{equation}
\end{linenomath}
which has to be understood in a weak sense, both as regards the spatial as well as the
stochastic dependence.
%\begin{align}\label{Eq:StoForm}
%-\textrm{div}(\conducttensor(\event) \cdot \nabla\temperature(x,\event))   &= \heatsource(x),  \quad \forall x \in \spatialdomain, \event\in\samplespace, \nonumber \\
%{\temperature}(x,\event) &= {\temperature}_0,  \quad \forall x \in \Gamma_D, \event\in\samplespace, \\
%(\conducttensor(\event) \cdot \nabla\temperature(x,\event))\cdot\vek{n}(x) &= {g}(x), \quad \forall x \in \Gamma_N, \event\in\samplespace, \nonumber
%\end{align}
The goal is to determine the random temperature field 
$\temperature(\vek{x},\event):\spatialdomain\times\samplespace\rightarrow\mathbb{R}$, assuming deterministic boundary conditions and heat source. 

% Note that here the boundary conditions and heat source stay deterministic. 
Given the weak form of \feq{Eq:stochastic-thm}, we use the finite element method to seek 
for an approximate solution $\temperatureApprox(\vek{x},\event)$ 
on a finite-dimensional subspace $\finitespace\subset\solspace$ of the solution space $\solspace$ of \feq{Eq:deterministic-thm}, where as usual $h$ is the discretisation parameter---
an indicator of element size.
Subsequently, the objective of this study is to determine the statistics like the mean 
and standard deviation of the solution $\temperatureApprox(\vek{x},\event)$ using the 
classical Monte Carlo method (\acrshort{mc}) \cite{metropolis_monte_1949, fishman_monte_1996}.
Given, i.i.d.\ $\{\conducttensor(\event_i)\}_{i=1}^{N}$ samples, where $\samplesize$ is the 
sample size, one may compute the statistics of discretised solution 
$\stotemp:=\{\temperatureApprox(\vek{x},\event_i)\}_{i=1}^{N}$.  The unbiased sample mean
 \cite{dekking_modern_2005} is then calculated as
\begin{linenomath}
\begin{equation}
	\meanMC(\stotemp) := \frac{1}{\samplesize}\sum_{i=1}^{\samplesize}
	     \temperatureApprox(\vek{x},\event_i), 
\end{equation}
\end{linenomath}
and the corrected sample standard deviation as
\begin{linenomath}
\begin{equation}
\stdMC(\stotemp) := \sqrt{\frac{1}{\samplesize-1}\sum_{i=1}^{\samplesize}
    (\temperatureApprox(\vek{x},\event_i)-\meanMC(\stotemp))^2}.
\end{equation}
\end{linenomath}

Along with the temperature field, we also evaluate the heat flux vector field 
$\heatflux(\vek{x},\event):\spatialdomain\times\samplespace \rightarrow \Euclideanspace$, whose approximate solution
% By determining an approximation to the random heat flux vector field 
$\heatflux_h(\vek{x},\event)$ is determined on a finite-dimensional subspace $\solspaceheat_h\subset\solspaceheat$, where $\solspaceheat$ is the solution space of heat flux, $\heatflux(\vek{x})=\conducttensor \nabla T(\vek{x})$.
% $\heatflux_h(\vek{x},\event):\spatialdomain\times\samplespace \rightarrow \Euclideanspace$, 
Further, we determine 
the second order statistics of the magnitude of approximate heat flux vector field in the Euclidean norm i.e.\ $q_h^{(t)}(\vek{x},\event) = 
\|\heatflux_h(\vek{x},\event)\|$.
More specifically, the statistics such as the mean and standard deviation are obtained 
by sampling based MC estimators $\meanMC(\stototalflux)$ and 
$\stdMC(\stototalflux)$, respectively.
% Therefore, 
% \acrshort{mc} is implemented to evaluate the statistics of heat flux $\stoflux:=[\heatflux_h(x,\event_i)]_{i=1}^{N}$. More specifically, the mean and variance of total heat flux $\stototalflux = \|\stoflux\|$, are determined by estimators $\meanMC(\stototalflux)$ and $\varMC(\stototalflux)$ respectively.

Apart from evaluating the statistics of magnitude of heat flux field $\heatflux_h(\vek{x},\event)$, 
one may also be interested in 
understanding the statistics of direction of heat flow.
% determining the  statistics of directional aspect of the vector field. 
To this, one may compute the second order circular statistics of the normalized heat flux, 
$\hat{\heatflux}_h(\vek{x},\event) = {\heatflux_h(\vek{x},\event)}/{q_h^{(t)}(\vek{x},\event)}$.
% In addition, to estimate the circular statistics of the sampled response 
% $\stoflux:=\{\heatflux_h(\vek{x},\event_i)\}_{i=1}^{N}$ 
% \cite{mardia_directional_2000,jammalamadaka_topics_2001}, one is interested in 
% determining the orientation and not only the magnitude of $\stoflux$. 
% Thus, for convenience, the quantity $\heatflux_h$ is normalised as 
% \begin{equation}\label{Eq:normheatflux}
% 	\hat{\heatflux}_h = \frac{\heatflux_h}{q_h^{(t)}}.
% \end{equation}
%Thereby, the unit heat flux vector field $\stofluxunit$ belongs to the unit sphere 
%$ \unitsphere$.
Subsequently, the \acrshort{mc} mean estimate of sampled unit vector field 
$\stofluxunit:=\{\hat{\heatflux}_h(\vek{x},\event_i)\}_{i=1}^{N}$ 
% $\hat{\heatflux}_h(\vek{x},\event)$ 
can be estimated accordingly \cite{mardia_directional_2000,jammalamadaka_topics_2001}:
\begin{linenomath}
\begin{equation}\label{Eq:meannormheatflux}
	\meanMCbold(\stofluxunit) = \frac{1}{\samplesize}\sum_{i=1}^{\samplesize}
	   \hat{\heatflux}_h(\vek{x},\event_i) = L\,  \muMCcirc(\stofluxunit),
\end{equation}
\end{linenomath}
where $L = \|\meanMCbold(\stofluxunit)\| \leq 1$ represents resultant Euclidean length 
of the sample mean $\meanMCbold(\stofluxunit)\in\Euclideanspace$, and 
$\muMCcirc(\stofluxunit) = \meanMCbold(\stofluxunit)/\|\meanMCbold(\stofluxunit)\|$ 
is a unit vector defining the sample mean direction of $\stofluxunit$. 
With the help of the former one, the sample circular standard deviation of $\stofluxunit$ 
can further be evaluated as
\begin{linenomath}
\begin{equation}\label{Eq:circvarnormheatflux}
	\stdMCcirc(\stofluxunit) = \sqrt{-2\log L}.
\end{equation}
\end{linenomath}

% where $\stdMCcirc(\stofluxunit)\in[0,\infty]$.

%% file: results.tex
% !TEX root = ../paper_tensors.tex
% !TEX encoding = UTF-8 Unicode

%macros for plot size (matlabtotikz plots)
\newlength\figH
\newlength\figW
\setlength{\figH}{3.9cm}
\setlength{\figW}{3.9cm}

\newlength\figHh
\newlength\figWw
\setlength{\figHh}{6.2cm}
\setlength{\figWw}{7.5cm}

%\newpage
\subsection{Results on 2D proximal femur}\label{Sec:2DResults}
%\begin{wrapfigure}{r}{0.3\textwidth}
%	\centering 
%	\def\svgwidth{0.3\textwidth}  
%	\input{ellipse.eps_tex}
%	\caption{Ellipse}
%	\label{Fig:ellipse}
%\end{wrapfigure}
%
%The equivalent geometric reprsentation of $\genSecTensor$  is a quadric, more specifically, in $d=2$ an ellipse, as shown in \ffig{Fig:ellipse} \cite{malgrange_symmetry_2014}.
%It is a well known proposition that the spectral decomposition of $\genSecTensor$ will always have strictly positive eigenvalues. Hence, assuming $\lambda_1\in\mathbb{R^+}$ and $\lambda_2\in\mathbb{R^+}$ s.t. $\lambda_1<\lambda_2$ as the eigenvalues of $\genSecTensor$ \textcolor{red}{orientation?}. Therefore, the semi-axes $a$ and $b$ of the ellipse are defined as:
%\begin{equation}
%a = \frac{1}{\sqrt{\lambda_1}},\: b = \frac{1}{\sqrt{\lambda_2}}.
%\end{equation} 
\input{results_input_2D-y}

%\clearpage
\input{results_2D}
%\clearpage
% \newpage
\subsection{Results on 3D proximal femur}\label{Sec:3DResults}
\input{results_input_3D}

%\clearpage
\input{results_3D}

%\clearpage

%% file: results_input_2D-y.tex
% !TEX root = ../paper_tensors.tex
% !TEX encoding = UTF-8 Unicode

\begin{figure}[ht]
	\centering
	\includegraphics[width=0.3\linewidth]{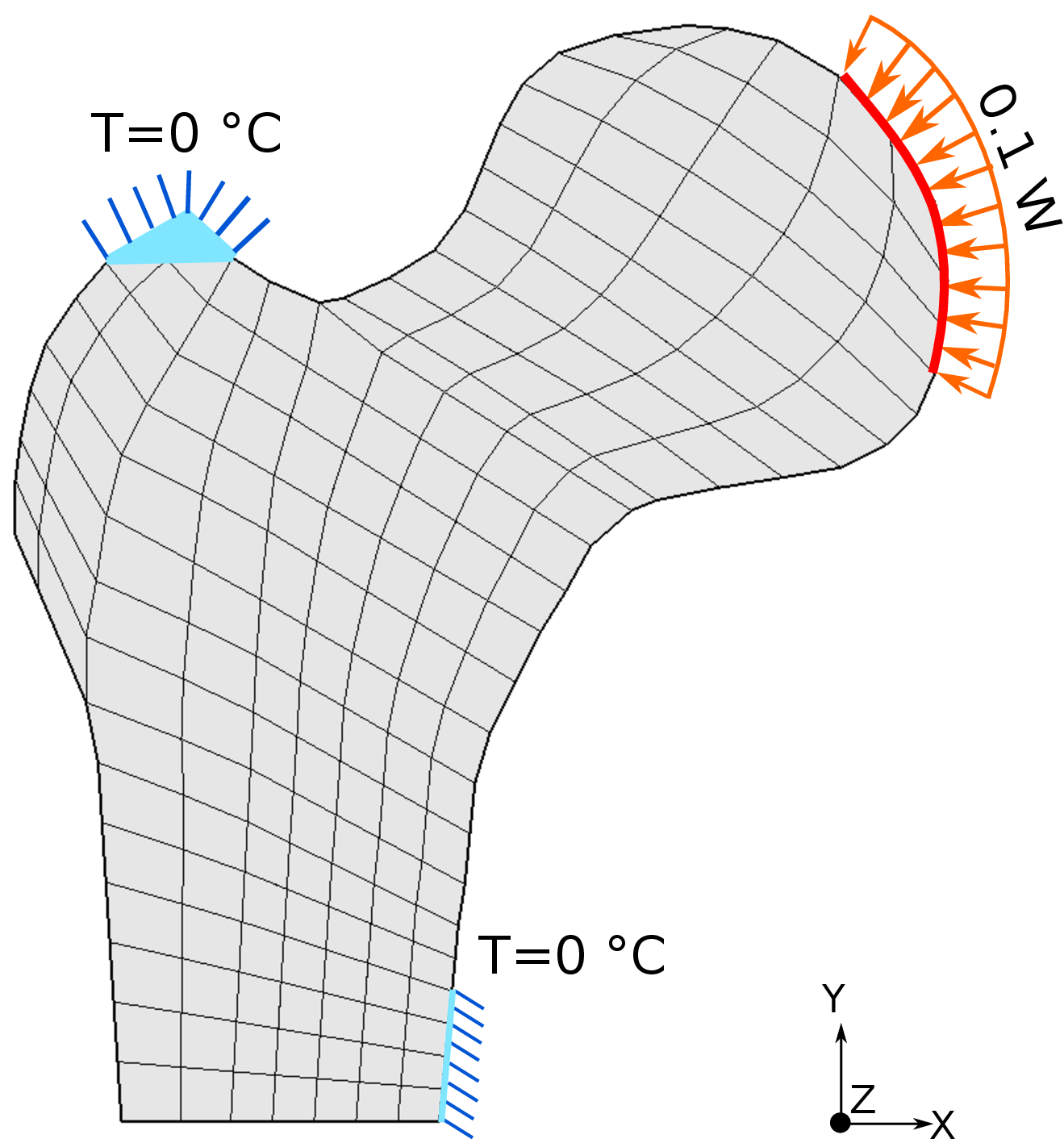}
	\captionof{figure}{Boundary conditions of 2D femur bone}
	\label{Fig:BC}	
\end{figure}
	
A two-dimensional proximal femur with a body of approximately 7 cm in width and 
21.7 cm in height is considered, see \ffig{Fig:BC}.  
A uniformly distributed surface 
heat flux with a resultant heating of 0.1 W and a fixed temperature of $0\ ^\circ C$ 
are applied at locations as shown in \ffig{Fig:BC}. 
The computational domain is discretized by the finite element discretization 
method with 171 four-noded plane-stress elements, %see   
%\ftbl{Tbl:meshspec} for the mesh specifications. 
giving a mesh with 206 nodes and 390 degrees of freedom.
In this study, each deterministic 
simulation obtained by sampling the stochastic parametric space is solved by the 
Calculix solver \cite{Dhondt_CCX_2004}. 
For the \acrshort{mc} simulation, the %Karl 
Pearson coefficient of dispersion $\delta=0.1$ 
(for the scaling uncertainty) and  $\samplesize=10^5$ samples are considered for all the cases 
in this study.
%\begin{table}[ht]
%\centering
%	\begin{tabular}{lll}
%		\specialrule{1pt}{1pt}{1pt}			
%		Elements & Nodes & \begin{tabular}[c]{@{}l@{}}
%			Degrees of \\ freedom 
%		\end{tabular} \\ \midrule
%		171 & 206 & 390 \\  
%		\specialrule{1pt}{1pt}{1pt}
%	\end{tabular}
%	\caption{Mesh specifications of 2D femur}
%	\label{Tbl:meshspec}
%\end{table}
% \vspace{0.6cm}
% \begin{wraptable}[7]{r}{0.39\textwidth}
% 	\begin{tabular}{lll}
% 		\specialrule{1pt}{1pt}{1pt}			
% 		Elements & Nodes & \begin{tabular}[c]{@{}l@{}}
% 			Degrees of \\ freedom 
% 		\end{tabular} \\ \midrule
% 		171 & 206 & 390 \\  
% 		\specialrule{1pt}{1pt}{1pt}
% 	\end{tabular}
% 	\caption{Mesh specifications of 2D femur}
% 	\label{Tbl:meshspec}
% \end{wraptable}

For simplification, we use abbreviations defined in \ftbl{Tbl:abbreviation} which are 
used to refer different modelling scenarios considered in this study.  Each model will 
have a unique name in the form: reference tensor symmetry-realisation symmetry-type of model. 
For instance, ``iso-ortho-scl'' represents a model with random scaling only, such that the 
reference tensor belongs to isotropic symmetry, whereas the realisations are of orthotropic 
symmetry.  Three such possible scenarios are explored in this article---``iso-iso-scl'' 
(isotropy-isotropy-random scaling), ``iso-ortho-scl'' (isotropy-orthotropy-random scaling), 
and ``ortho-ortho-dir'' (orthotropy-orthotropy-random direction).
\begin{table}[h!]
	\centering 
	\begin{tabular}{llll}
	\specialrule{1pt}{1pt}{1pt}
	\multicolumn{2}{l}{\begin{tabular}[c]{@{}l@{}}Mean/realisation  symmetry\end{tabular}} & \multicolumn{2}{l}{Type of model}                                                                                    \\ \midrule
	isotropy                                   & orthotropy                                  & \begin{tabular}[c]{@{}l@{}}random scaling\end{tabular} & \begin{tabular}[c]{@{}l@{}}random direction\end{tabular} \\ \midrule
	\acrshort{iso}                                       & \acrshort{ortho}                                       & \acrshort{scl}                                                      & \acrshort{dir}   \\ \specialrule{1pt}{1pt}{1pt}                                                   
	\end{tabular}
\caption{Model abbreviations\label{Tbl:abbreviation}} 	
\end{table} 
Consequently, using the subscripts defined in \fsec{Ssec:group-inv-detl} for random scaling 
and random orientation, two reference conductivity tensors---
$\meanconducttensorscl\in\PosSymSecRankSpacetwoD$ belonging to the isotropic symmetry class
(denoted by superscript $iso$) and $\meanconducttensorortho\in\PosSymSecRankSpacetwoD$ 
with orthotropic symmetry (with superscript $ortho$)---are considered, whose values are 
tabulated in \ftbl{Tbl:meantensor}.  
% Note that the eigenvalues of the model $\meanconducttensorortho$ are 
%$\meaneigvalscleleone=0.54$ W/mK and $\meaneigvalscleletwo=1$ W/mK.
Furthermore, we simulate the grain orientations of the femur in a spatially homogeneous 
sense. As a result, the eigenvectors of the tensor $\meanconducttensorortho$ corresponding 
to the eigenvalues $\meaneigvalroteleone$ and $\meaneigvalroteletwo$ are oriented at an 
angle of $\pi/4$ in a clockwise direction to the x and y-axis, respectively. 
\begin{table}[h!]
	\centering
	\begin{tabular}{ll}
	\specialrule{1pt}{1pt}{1pt}
	 \begin{tabular}[c]{@{}l@{}}isotropy\\ (W/mK)\end{tabular} & \begin{tabular}[c]{@{}l@{}}orthotropy\\ (W/mK)\end{tabular} \\ 
	 $\meanconducttensorscl$ & $\meanconducttensorortho$\\ \midrule 
	 $
	\begin{bmatrix}
	0.54 & 0  \\ 0 & 0.54
	\end{bmatrix} $  &  $  \begin{bmatrix}
	0.77 & 0.23  \\ 0.23 & 0.77
	\end{bmatrix} $ \\ \midrule
	\multicolumn{2}{c}{Eigenvalues (W/mK)} \\  \midrule
	\begin{tabular}[c]{@{}l@{}}$\meaneigvalscleleone=0.54$,\\ $\meaneigvalscleletwo=0.54$\end{tabular} & \begin{tabular}[c]{@{}l@{}}$\meaneigvalroteleone=0.54$,\\ $\meaneigvalroteletwo=1$\end{tabular} \\
	\specialrule{1pt}{1pt}{1pt}            
	\end{tabular}
\caption{Reference conductivity tensors with corresponding eigenvalues of 2D femur\label{Tbl:meantensor}}
\end{table}
% For simplification, we use abbreviations defined in \ftbl{Tbl:abbreviation} which refer to different modelling scenarios. Each model will have a unique name in the form: reference tensor symmetry-realisation symmetry-type of model.

\begin{figure}[h!]
	\centering
	\begin{subfigure}{\textwidth}
		\centering
		\input{2D_LN_iso_iso}
		\caption{}
		\label{Fig:isoisoLN}
	\end{subfigure}	
	\begin{subfigure}{\textwidth}
		\centering
		\includegraphics[width=0.7\linewidth]{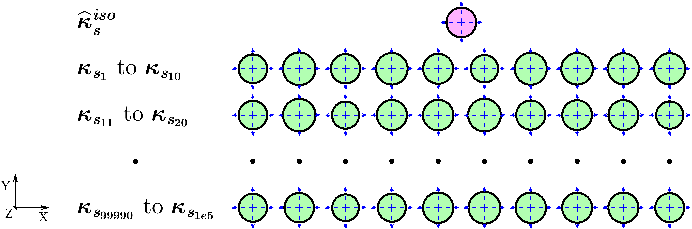}
		\caption{}
		\label{Fig:isoisoquad}
	\end{subfigure}	
	\caption{Stochastic tensor $\conducttensorscl(\event_s)$ with fixed symmetry 
	(\acrshort{iso}-\acrshort{iso}-\acrshort{scl}): \subref{Fig:isoisoLN} log-normal 
	\acrshort{pdf} of identical random scaling values $\eigvalscleleone(\event_s) 
	=\eigvalscleletwo(\event_s)$, \subref{Fig:isoisoquad} geometric visualization of reference 
	tensor $\meanconducttensorscl$ and realisations of random tensor $\conducttensorscl(\event_s)$; 
	the geometries are enhanced by scaling the eigenvalues by a factor of 0.5} 
	\label{Fig:isoisoscl}
\end{figure}

\paragraph{Random scaling with fixed symmetry:}
the first modelling scenario is the one in which 
we perform random scaling only, such that 
the reference tensor $\meanconducttensorscl$, as well as the realisations, belong to 
isotropic symmetry (denoted by \acrshort{iso}-\acrshort{iso}-\acrshort{scl}); 
see \ffig{Fig:isoisoscl} for a schematic overview of the model.
In such a case, %we model 
the random tensor $\conducttensorscl(\event_s)\in\PosSymSecRankSpacetwoD$ 
%the exponential of the model depicted in \feq{Eq:rndscl}   
is modelled---as in \feq{Eq:SpecDec-rep-spl}
---with random scaling only.  Here the uncertain scaling parameters 
$\eigvalscleleone(\event_s)=\eigvalscleletwo(\event_s)$ are modelled as identical 
%and independent 
log-normal random variables. 
% \ffig{Fig:isoisoscl} presents an overview of the stochastic model $\conducttensorscl(\event_s)$, such that fixed isotropic symmetry in the mean and realisations are maintained---denoted as \acrshort{iso}-\acrshort{iso}-\acrshort{scl}. 
Their corresponding \acrshort{pdf}s are shown in \ffig{Fig:isoisoLN}. 
A geometrically equivalent (circular) representation of reference tensor and realisations 
$\{\conducttensorscl(\event_{s_i})\}^\samplesize_{i=1}\equiv\{\conducttensor_{s_i}\}^\samplesize_{i=1}$,
where the radius of circle is defined by ${1}/{\sqrt{\eigvalscleleone(\event_s)}}$ and 
${1}/{\sqrt{\eigvalscleletwo(\event_s)}}$ \cite{malgrange_symmetry_2014}, is presented 
in \ffig{Fig:isoisoquad}. 
% Note that, in this article, all the geometries representing the random tensors are 
% enhanced by scaling the eigenvalues by a factor of 0.5.
As can be seen the preservation of isotropic symmetry is apparent given the constant 
circular shape among the random variable realisations. However, 
the varying size of circles signifies 
a variation in scaling parameters.

\begin{figure}[h!]
	\centering
	\begin{subfigure}{\textwidth}
		\centering
		\input{2D_LN_iso_ortho}
		\caption{}
		\label{Fig:isoorthoLN}
	\end{subfigure}
	\begin{subfigure}{\textwidth}
		\centering
		\includegraphics[width=0.7\linewidth]{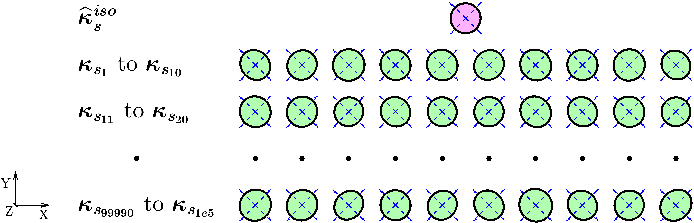}
		\caption{}
		\label{Fig:isoorthoquad}
	\end{subfigure}	
	\caption{Stochastic tensor $\conducttensorscl(\event_s)$ with varying symmetry 
	(\acrshort{iso}-\acrshort{ortho}-\acrshort{scl}): \subref{Fig:isoorthoLN} log-normal 
	\acrshort{pdf} of i.i.d.\ random scaling values $\{\eigvalsclele\}^2_{i=1}$, 
	\subref{Fig:isoorthoquad} geometric visualization of reference tensor 
	$\meanconducttensorscl$ and realisations of random tensor $\conducttensorscl(\event_s)$; 
	the geometries are enhanced by scaling the eigenvalues by a factor of 0.5}
	\label{Fig:isoorthoscl}
\end{figure}

\paragraph{Random scaling with varying symmetry:}
in the next scenario, we consider a similar stochastic tensor 
$\conducttensorscl(\event_s)\in\PosSymSecRankSpacetwoD$ as described previously with a 
similar reference tensor $\meanconducttensorscl$.  The difference being that the eigenvalues 
$\{\eigvalsclele\}^2_{i=1}$ are modelled as independent identically distributed (i.i.d.) 
log-normal random variables.  Due to which the realisations of model 
$\conducttensorscl(\event_s)$ belong to orthotropic symmetry, thus, 
this model is denoted as \acrshort{iso}-\acrshort{ortho}-\acrshort{scl}. 
Furthermore, in orthotropic realisations, the eigenvectors of $\conducttensorscl(\event_s)$ 
related to random eigenvalues $\eigvalscleleone(\event_s)$ and $\eigvalscleletwo(\event_s)$ 
are constrained at an angle of $\pi/4$ in clockwise direction with respect to x and y-axis 
respectively. 
The schematic representation of this model with varying symmetry is displayed in 
\ffig{Fig:isoorthoscl}, in which
the corresponding \acrshort{pdf}s are shown in \ffig{Fig:isoorthoLN}. 
% \ffig{Fig:isoorthoquad} presents the graphical interpretation of the model.
Closer inspection of the \ffig{Fig:isoorthoquad} shows the shift from circular geometry 
of isotropic reference tensor to orthotropic elliptical shape in realisations.  
One may also notice the variation in size of ellipses in the realisations as the 
semi-axes of ellipse are determined by ${1}/{\sqrt{\eigvalscleleone(\event_s)}}$ and 
${1}/{\sqrt{\eigvalscleletwo(\event_s)}}$.  

\begin{figure}[h!]
%	\begin{subfigure}{0.32\textwidth}
%		\centering
%		\input{2D_LN_ortho_ortho}
%		\caption{}
%	\end{subfigure}\hspace{0.45cm}
	\begin{subfigure}{0.49\textwidth}
		\centering
		\includegraphics[width=0.7\linewidth]{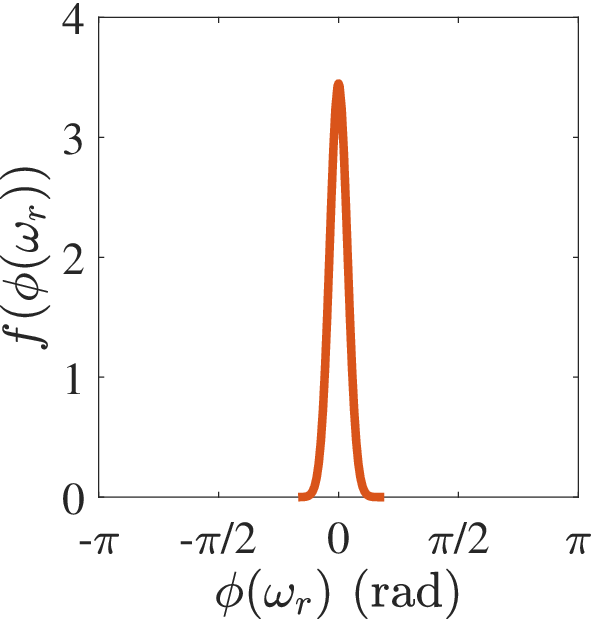} %tex files are too big, hence .eps
		\caption{}
		\label{Fig:orthoorthoVM}
	\end{subfigure}
	\begin{subfigure}{0.49\textwidth}
		\centering
		\includegraphics[width=0.75\linewidth]{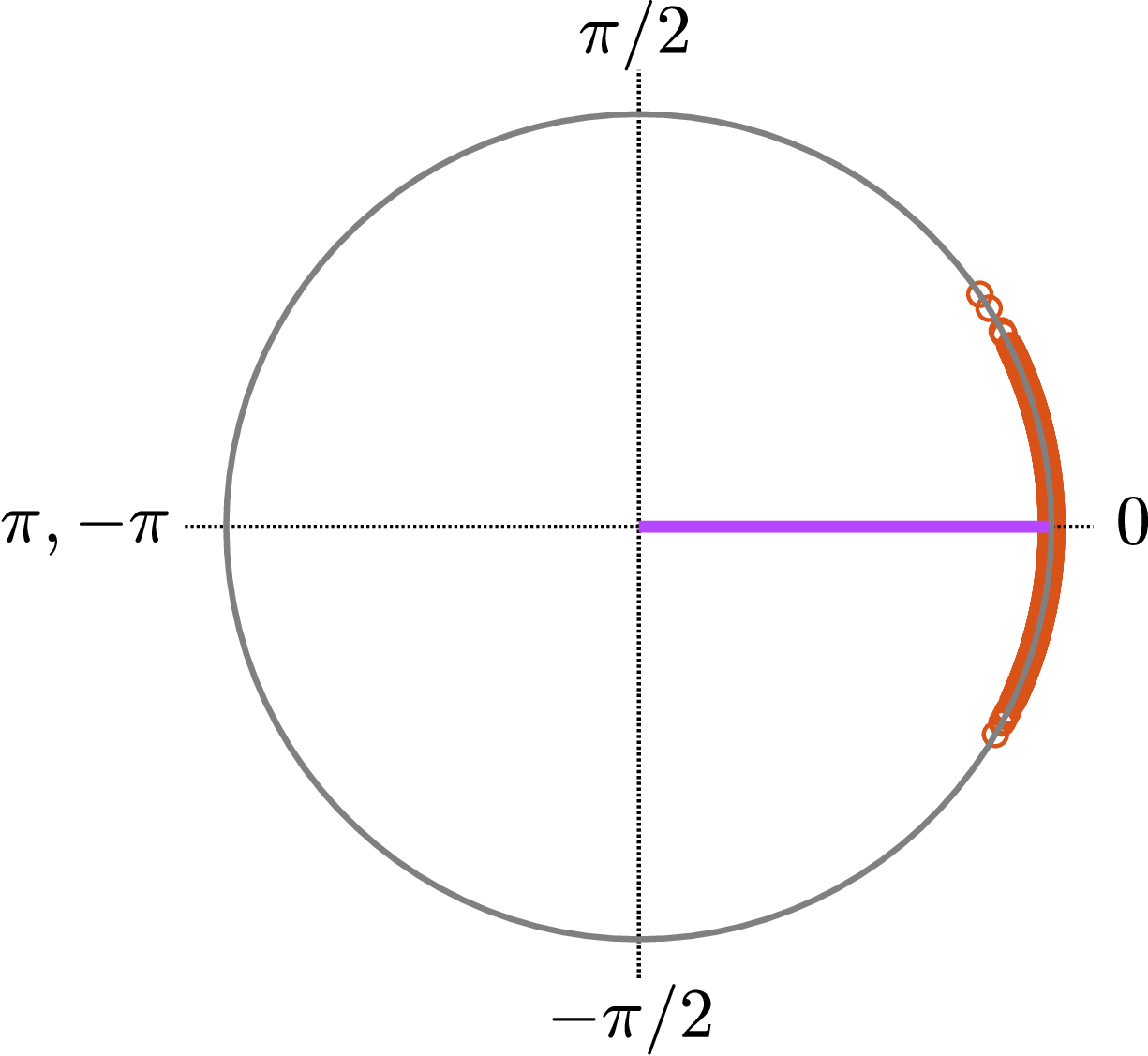}
		\caption{}
		\label{Fig:orthoorthoVMunit}
	\end{subfigure}
%	\caption{(a) Lognormal PDF of eigenvalues $\varLambda_1$ and $\varLambda_2$, (b) von Mises PDF of orientation angle $\rotangleVM$ in Cartesian coordinates, (c) realisations of $\rotangleVM(\event)$ on unit circle with resultant mean vector (red line)}
%	\begin{subfigure}{0.49\textwidth}
%		\centering
%		\input{2D_ortho_ortho_scl}
%		\caption{}
%	\end{subfigure}
	\begin{subfigure}{\textwidth}
		\centering
		\includegraphics[width=0.7\linewidth]{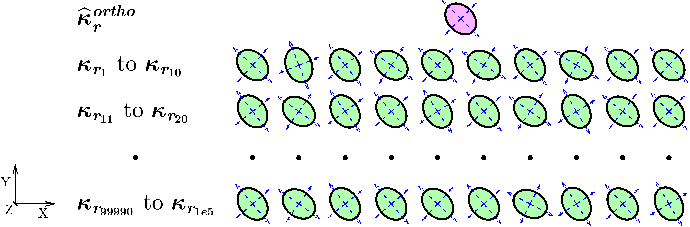}
		\caption{}
		\label{Fig:orthoorthorot}
	\end{subfigure}
%	\begin{subfigure}{\textwidth}
%		\centering
%		\input{2D_ortho_ortho_rotscl}
%		\caption{}
%	\end{subfigure}		
	\caption{Stochastic tensor $\conducttensorrot(\event_r)$ with fixed symmetry (\acrshort{ortho}-\acrshort{ortho}-\acrshort{dir}): \subref{Fig:orthoorthoVM} von Mises \acrshort{pdf} of random rotation angle $\rotangle(\event_r)$, \subref{Fig:orthoorthoVMunit}  realisations of circular random variable $\rotangle(\event_r)$ on unit circle with resultant mean vector (solid straight line), \subref{Fig:orthoorthorot}  geometric visualization of reference tensor $\meanconducttensorortho$ and realisations of random tensor $\conducttensorrot(\event_r)$; the geometries are enhanced by scaling the eigenvalues by a factor of 0.5, and the orientations with respect to the mean are scaled by a factor of 1.5}
\end{figure}

\paragraph{Random direction with fixed symmetry:}
in the last modelling scenario, we build a random tensor $\conducttensorrot(\event_r)\in\PosSymSecRankSpacetwoD$, i.e. 
the exponential mapping of stochastic tensor as shown in
% $\conducttensorrot(\event_r)=\exp(\symtensor_r(\event_r))$ from
\feq{Eq:stomodrot}, with random orientation only, given the reference tensor $\meanconducttensorortho$. Here the model is designated as \acrshort{ortho}-\acrshort{ortho}-\acrshort{dir}, since, a fixed orthotropic symmetry is preserved in mean and realisations. 
The random angle of rotation $\rotangle(\event_r)$ is modelled as a von Mises random variable (see \feq{Eq:VMPDF}) with zero mean $\circmeananly=0$ with respect to x-axis and concentration parameter $\concparam=75$.
% Considering the mean direction of the rotation angle $\rotangle$ which is set at $\circmeananly=0$ with respect to x-axis and concentration parameter $\concparam=75$, the angular parameter $\rotangle(\event_r)$ is modelled as a von Mises random variable (see \feq{Eq:VMPDF}). 
In \ffig{Fig:orthoorthoVM}, the \acrshort{pdf} $f(\rotangle(\event_r))$ is shown on a real line over the domain $(-\pi,\pi]$, whereas the realisations of random variable $\rotangle(\event_r)$ and its sample mean vector (solid straight line) are plotted on a unit circle in \ffig{Fig:orthoorthoVMunit}.
In this study, we model directional uncertainty on the $x-y$ plane with $z$ as the axis of rotation.
Thereby, the corresponding random rotation matrix $\rotmat(\event_r)$ in \feq{Eq:stomodrot} takes the form:
\begin{linenomath}
\begin{equation}\label{Eq:Rmatrix}
\rotmat(\event_r) =
\begin{bmatrix}
\cos\rotangle(\event_r) & -\sin\rotangle(\event_r)  \\ \sin\rotangle(\event_r) & \cos\rotangle(\event_r)
\end{bmatrix}.
\end{equation}
\end{linenomath}
Here we consider a right-handed Cartesian coordinate system, meaning that positive values of random variable $\rotangle(\event_r)$ indicate counterclockwise rotation and vice versa.  
 \ffig{Fig:orthoorthorot} depicts the corresponding elliptical geometric visualization of the model.
% For easier visualization, we scale the random orientations in the realisations of all the random tensor models in this study by a factor of 1.5.
As can be seen that the orientation of elliptical semi-axes fluctuates in the realisations $\{\conducttensorrot(\event_{r_i})\}^\samplesize_{i=1}\equiv\{\conducttensor_{r_i}\}^\samplesize_{i=1}$, whereas its size and shape remain constant in the mean and realisations.

%% file: 2D_LN_iso_iso.tex
\definecolor{mycolor1}{rgb}{0.85000,0.33000,0.10000}%
\begin{tikzpicture}

\begin{axis}[%
width=\figW,
height=\figH,
at={(0\figW,0\figH)},
scale only axis,
xmin=0.333523780208734,
xmax=0.837879846496155,
xlabel style={font=\color{white!15!black}},
xlabel={$\eigvalsclele$ (W/mK)},
ymin=0,
ymax=10,
ylabel style={font=\color{white!15!black}},
ylabel={$f(\eigvalsclele)$},
axis background/.style={fill=white},
legend style={legend cell align=left, align=left, draw=white!15!black}
]
\addplot [color=mycolor1, line width=2.0pt]
  table[row sep=crcr]{%
0.333523780208734	1.06781483367869e-05\\
0.338618285928809	0.000125146641554135\\
0.343712791648884	0.000716457938709621\\
0.348807297368959	0.00217668687921683\\
0.353901803089034	0.00409863064763616\\
0.358996308809109	0.00599991909925284\\
0.364090814529184	0.00851290475255793\\
0.369185320249259	0.0123604328713531\\
0.374279825969334	0.0180406670621828\\
0.379374331689409	0.0285819605209717\\
0.384468837409484	0.0459381598916428\\
0.389563343129558	0.0704544887938538\\
0.394657848849633	0.104238695378309\\
0.399752354569708	0.143056754854802\\
0.404846860289783	0.183653543607387\\
0.409941366009858	0.245699992984782\\
0.415035871729933	0.346856117540483\\
0.420130377450008	0.472334484468659\\
0.425224883170083	0.619521900435522\\
0.430319388890158	0.807926315578874\\
0.435413894610233	1.03153135123811\\
0.440508400330308	1.27675844826076\\
0.445602906050383	1.56822695790261\\
0.450697411770458	1.93466838910659\\
0.455791917490533	2.35171899056274\\
0.460886423210608	2.76026256147619\\
0.465980928930683	3.12830792958184\\
0.471075434650758	3.49395858512132\\
0.476169940370833	3.95004734444102\\
0.481264446090908	4.46104369436053\\
0.486358951810983	4.9352236190815\\
0.491453457531058	5.38499946926963\\
0.496547963251133	5.81812586701607\\
0.501642468971208	6.22495224653765\\
0.506736974691283	6.60116091337403\\
0.511831480411358	6.90551905528756\\
0.516925986131433	7.13197735175847\\
0.522020491851508	7.29783512164279\\
0.527114997571583	7.37677648516861\\
0.532209503291657	7.3848673878748\\
0.537304009011732	7.350757067544\\
0.542398514731807	7.26936917075195\\
0.547493020451882	7.1483706858905\\
0.552587526171957	6.95808570472295\\
0.557682031892032	6.68989057937534\\
0.562776537612107	6.38492147524128\\
0.567871043332182	6.04990745331803\\
0.572965549052257	5.68909071410293\\
0.578060054772332	5.32729741460774\\
0.583154560492407	4.94701146789031\\
0.588249066212482	4.52102252485428\\
0.593343571932557	4.10307521331883\\
0.598438077652632	3.73370039204291\\
0.603532583372707	3.38297386481366\\
0.608627089092782	3.03681445572067\\
0.613721594812857	2.7242268959603\\
0.618816100532932	2.42036939233373\\
0.623910606253007	2.1134364500586\\
0.629005111973082	1.84526115563693\\
0.634099617693157	1.62140457548973\\
0.639194123413232	1.41988422503323\\
0.644288629133307	1.23334926380314\\
0.649383134853382	1.05090955940357\\
0.654477640573457	0.878150068597991\\
0.659572146293532	0.728911514193334\\
0.664666652013607	0.600877950472566\\
0.669761157733682	0.495564348216475\\
0.674855663453756	0.41176498236751\\
0.679950169173831	0.342409220982515\\
0.685044674893906	0.285931135081324\\
0.690139180613981	0.24557089123083\\
0.695233686334056	0.214946033632921\\
0.700328192054131	0.182896625712956\\
0.705422697774206	0.15118682517046\\
0.710517203494281	0.122172211524413\\
0.715611709214356	0.0933209280939708\\
0.720706214934431	0.0700961431863263\\
0.725800720654506	0.0579876649314881\\
0.730895226374581	0.0511757188366701\\
0.735989732094656	0.0416054729495341\\
0.741084237814731	0.0323044645374202\\
0.746178743534806	0.0244529920176591\\
0.751273249254881	0.0172937125149005\\
0.756367754974956	0.013071629841811\\
0.761462260695031	0.0111426333050941\\
0.766556766415106	0.00976113634348158\\
0.771651272135181	0.00770528181986963\\
0.776745777855256	0.00566994494504216\\
0.781840283575331	0.00534675574701623\\
0.786934789295406	0.00551582981872053\\
0.792029295015481	0.00396608538901785\\
0.797123800735556	0.00201612240345166\\
0.802218306455631	0.00132565830598064\\
0.807312812175705	0.00113876363240426\\
0.81240731789578	0.000878806118647164\\
0.817501823615855	0.000796469781892713\\
0.82259632933593	0.000702722306721663\\
0.827690835056005	0.000342932949895322\\
0.83278534077608	7.72219562173645e-05\\
0.837879846496155	7.81257228519023e-06\\
};
\small \addlegendentry{$\eigvalscleleone,\eigvalscleletwo$}

\end{axis}
\end{tikzpicture}%

%% file: 2D_LN_iso_ortho.tex
\definecolor{mycolor1}{rgb}{0.85000,0.33000,0.10000}%
\definecolor{mycolor2}{rgb}{0.72000,0.27000,1.00000}%
\begin{tikzpicture}

\begin{axis}[%
width=\figW,
height=\figH,
at={(0\figW,0\figH)},
scale only axis,
xmin=0.333523780208734,
xmax=0.837879846496155,
xlabel style={font=\color{white!15!black}},
xlabel={$\eigvalsclele$ (W/mK)},
ymin=0,
ymax=10,
ylabel style={font=\color{white!15!black}},
ylabel={$f(\eigvalsclele)$},
axis background/.style={fill=white},
legend style={legend cell align=left, align=left, draw=white!15!black}
]
\addplot [color=mycolor1, line width=2.0pt]
  table[row sep=crcr]{%
0.318034371545172	7.80054195201883e-06\\
0.32335078655169	8.34350540148855e-05\\
0.328667201558207	0.000373519629052774\\
0.333983616564725	0.000723698007361162\\
0.339300031571242	0.000779500522112519\\
0.34461644657776	0.0010570870237011\\
0.349932861584277	0.00211418381032067\\
0.355249276590795	0.00390972075581172\\
0.360565691597312	0.00577602906348665\\
0.36588210660383	0.00764835258218549\\
0.371198521610347	0.0108470031101758\\
0.376514936616865	0.0175815420289205\\
0.381831351623382	0.0305712135709878\\
0.3871477666299	0.0500335133352594\\
0.392464181636417	0.0757616023472796\\
0.397780596642935	0.112843531480652\\
0.403097011649452	0.163674321132789\\
0.40841342665597	0.228347101654996\\
0.413729841662487	0.316157200130122\\
0.419046256669005	0.44213017449827\\
0.424362671675522	0.606761706607081\\
0.42967908668204	0.789116946083822\\
0.434995501688557	0.995544333789892\\
0.440311916695075	1.26269587338433\\
0.445628331701593	1.5847058272021\\
0.45094474670811	1.92974732026118\\
0.456261161714628	2.30979242350764\\
0.461577576721145	2.73750838312084\\
0.466893991727663	3.1832597070655\\
0.47221040673418	3.65380946894086\\
0.477526821740698	4.15804938172588\\
0.482843236747215	4.6782672931643\\
0.488159651753733	5.18687815644666\\
0.49347606676025	5.64659438700246\\
0.498792481766768	6.07898812127084\\
0.504108896773285	6.48665596102453\\
0.509425311779803	6.82681047808237\\
0.51474172678632	7.08218007544207\\
0.520058141792838	7.22309028449937\\
0.525374556799355	7.29401026079918\\
0.530690971805873	7.34459229297453\\
0.53600738681239	7.33894718582828\\
0.541323801818908	7.28962127417048\\
0.546640216825425	7.18537349986148\\
0.551956631831943	6.97364636310761\\
0.55727304683846	6.68761439435511\\
0.562589461844978	6.35297199386754\\
0.567905876851496	5.98389043872555\\
0.573222291858013	5.62581901494169\\
0.57853870686453	5.24809964397854\\
0.583855121871048	4.80446795098422\\
0.589171536877566	4.34523243549715\\
0.594487951884083	3.97029528778148\\
0.599804366890601	3.64976227727199\\
0.605120781897118	3.30931499206233\\
0.610437196903636	2.93634582789143\\
0.615753611910153	2.5672218765336\\
0.621070026916671	2.23911419115869\\
0.626386441923188	1.9508748177403\\
0.631702856929706	1.70032309621681\\
0.637019271936223	1.49006403509107\\
0.642335686942741	1.30469326917044\\
0.647652101949258	1.12915379706356\\
0.652968516955776	0.958181778232407\\
0.658284931962293	0.806702078589301\\
0.663601346968811	0.684051308927109\\
0.668917761975328	0.56772133025819\\
0.674234176981846	0.45878007966469\\
0.679550591988363	0.376596437533859\\
0.684867006994881	0.312395859376061\\
0.690183422001398	0.257037677280281\\
0.695499837007916	0.208949905593403\\
0.700816252014434	0.170604640826205\\
0.706132667020951	0.143284033709517\\
0.711449082027469	0.118446543008565\\
0.716765497033986	0.0951634981846204\\
0.722081912040504	0.075969862637935\\
0.727398327047021	0.0583798948341467\\
0.732714742053539	0.0428355193907454\\
0.738031157060056	0.0296975587681899\\
0.743347572066574	0.0204975974591209\\
0.748663987073091	0.0171115923750802\\
0.753980402079609	0.0160231391081761\\
0.759296817086126	0.0138330054559626\\
0.764613232092644	0.0106270735991171\\
0.769929647099161	0.00828655858374487\\
0.775246062105679	0.00626671760745505\\
0.780562477112196	0.00366000724108886\\
0.785878892118714	0.00213684424981666\\
0.791195307125231	0.00214228213003544\\
0.796511722131749	0.00211987690495529\\
0.801828137138266	0.00145473359094952\\
0.807144552144784	0.000946032197961263\\
0.812460967151301	0.000531501539210157\\
0.817777382157819	0.000311576755575317\\
0.823093797164336	0.00054950640555601\\
0.828410212170854	0.000690000593169896\\
0.833726627177371	0.000371362859552804\\
0.839043042183889	8.3378193808534e-05\\
0.844359457190407	7.79990985309992e-06\\
};
\addlegendentry{$\eigvalscleleone$}

\addplot [color=mycolor2, dashed, line width=2.0pt]
  table[row sep=crcr]{%
0.316843917720996	7.82251808631267e-06\\
0.322391161801415	9.14273277474556e-05\\
0.327938405881835	0.000412042799863609\\
0.333485649962254	0.000751830811541548\\
0.339032894042673	0.000815199743812523\\
0.344580138123093	0.00133528201170125\\
0.350127382203512	0.00239497531340307\\
0.355674626283931	0.00333184633743241\\
0.36122187036435	0.0052813136456002\\
0.36676911444477	0.00879433707549697\\
0.372316358525189	0.0128361885493858\\
0.377863602605608	0.0188728218649595\\
0.383410846686028	0.0299734531879551\\
0.388958090766447	0.049723698801178\\
0.394505334846866	0.0779078028647833\\
0.400052578927285	0.115753948040534\\
0.405599823007705	0.172417547785822\\
0.411147067088124	0.254655486541032\\
0.416694311168543	0.365589455032946\\
0.422241555248963	0.505827370594256\\
0.427788799329382	0.685196069260175\\
0.433336043409801	0.915248542895593\\
0.438883287490221	1.18613863356168\\
0.44443053157064	1.49648221444471\\
0.449977775651059	1.84491965019398\\
0.455525019731478	2.22867096358136\\
0.461072263811898	2.6794124028908\\
0.466619507892317	3.1966717802468\\
0.472166751972736	3.7036211119116\\
0.477713996053156	4.17863454517892\\
0.483261240133575	4.68053994471393\\
0.488808484213994	5.21054904816787\\
0.494355728294414	5.72150317486356\\
0.499902972374833	6.14419208047661\\
0.505450216455252	6.46420219852201\\
0.510997460535671	6.75052638474522\\
0.516544704616091	7.06587191448813\\
0.52209194869651	7.32295006149993\\
0.527639192776929	7.43059409341549\\
0.533186436857349	7.45746720587785\\
0.538733680937768	7.40045715813694\\
0.544280925018187	7.29340421718726\\
0.549828169098606	7.1460083842699\\
0.555375413179026	6.87678432022855\\
0.560922657259445	6.5127312023724\\
0.566469901339864	6.1151097726199\\
0.572017145420284	5.72814927525757\\
0.577564389500703	5.3509340489583\\
0.583111633581122	4.96387705936878\\
0.588658877661542	4.53666273577637\\
0.594206121741961	4.07422657240453\\
0.59975336582238	3.61442956645999\\
0.605300609902799	3.21307374190196\\
0.610847853983219	2.88580677707771\\
0.616395098063638	2.57847567948279\\
0.621942342144057	2.25082218446485\\
0.627489586224477	1.92859105345373\\
0.633036830304896	1.6503618458861\\
0.638584074385315	1.41760050667921\\
0.644131318465734	1.19375524707387\\
0.649678562546154	0.980963151536785\\
0.655225806626573	0.802837146922538\\
0.660773050706992	0.669738864609886\\
0.666320294787412	0.564963839311715\\
0.671867538867831	0.468764183126915\\
0.67741478294825	0.384472485902069\\
0.68296202702867	0.316388879641918\\
0.688509271109089	0.260854493463939\\
0.694056515189508	0.218354048954694\\
0.699603759269928	0.182846850144855\\
0.705151003350347	0.149166615909335\\
0.710698247430766	0.122046753137621\\
0.716245491511185	0.100454598336514\\
0.721792735591605	0.0773635517645659\\
0.727339979672024	0.0551428565042279\\
0.732887223752443	0.0392726755574335\\
0.738434467832863	0.0319940706934147\\
0.743981711913282	0.0285163036068819\\
0.749528955993701	0.0242829461971738\\
0.75507620007412	0.0182836715048358\\
0.76062344415454	0.0125553094966983\\
0.766170688234959	0.00794504818751512\\
0.771717932315378	0.00392039909792059\\
0.777265176395798	0.00138423785530635\\
0.782812420476217	0.000875777799821211\\
0.788359664556636	0.00140382198882671\\
0.793906908637056	0.00252501205540161\\
0.799454152717475	0.00320228610930239\\
0.805001396797894	0.00207155037843331\\
0.810548640878313	0.00060577433894272\\
0.816095884958733	7.81769439581497e-05\\
0.821643129039152	5.47272727778243e-05\\
0.827190373119571	0.000298162621082136\\
0.832737617199991	0.000677450306780269\\
0.83828486128041	0.000688330282213736\\
0.843832105360829	0.000650757925349482\\
0.849379349441248	0.00072642939775905\\
0.854926593521668	0.000410324253044\\
0.860473837602087	9.13599116136099e-05\\
0.866021081682506	7.82189948836931e-06\\
};
\addlegendentry{$\eigvalscleletwo$}

\end{axis}
\end{tikzpicture}%

%% file: results_2D.tex
% !TEX root = ../paper_tensors.tex
% !TEX encoding = UTF-8 Unicode

\newcommand{\widthtwoD}{1.0} %width size for 2D images in results section

\begin{figure}[!h]
	\begin{subfigure}[b]{0.32\textwidth}
		\centering
		\includegraphics[width=\widthtwoD\linewidth]{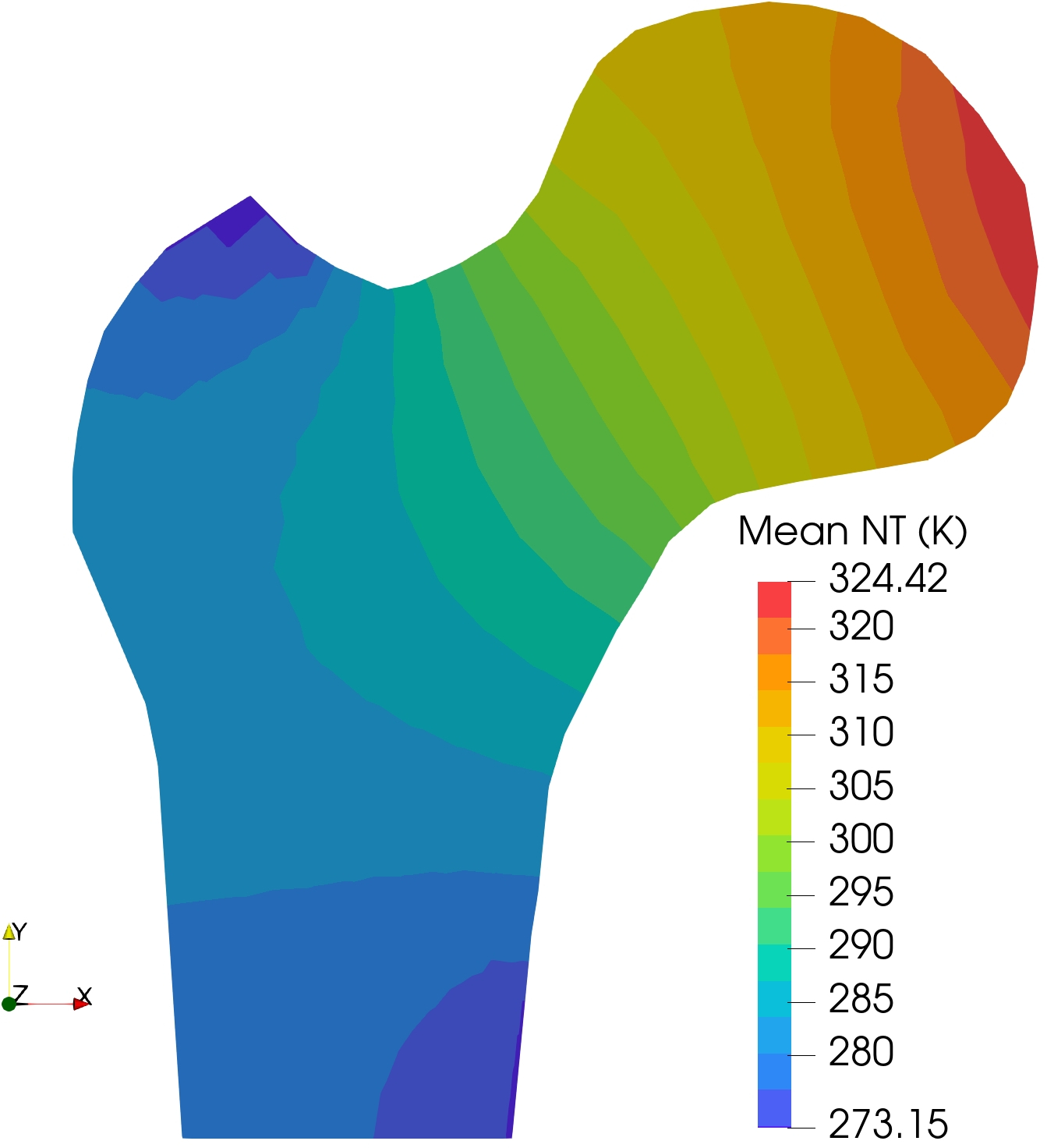}
		\caption{\acrshort{iso}-\acrshort{iso}-\acrshort{scl}}
		\label{Fig:meanNTa}
	\end{subfigure}
	\begin{subfigure}[b]{0.32\textwidth}
		\centering
		\includegraphics[width=\widthtwoD\linewidth]{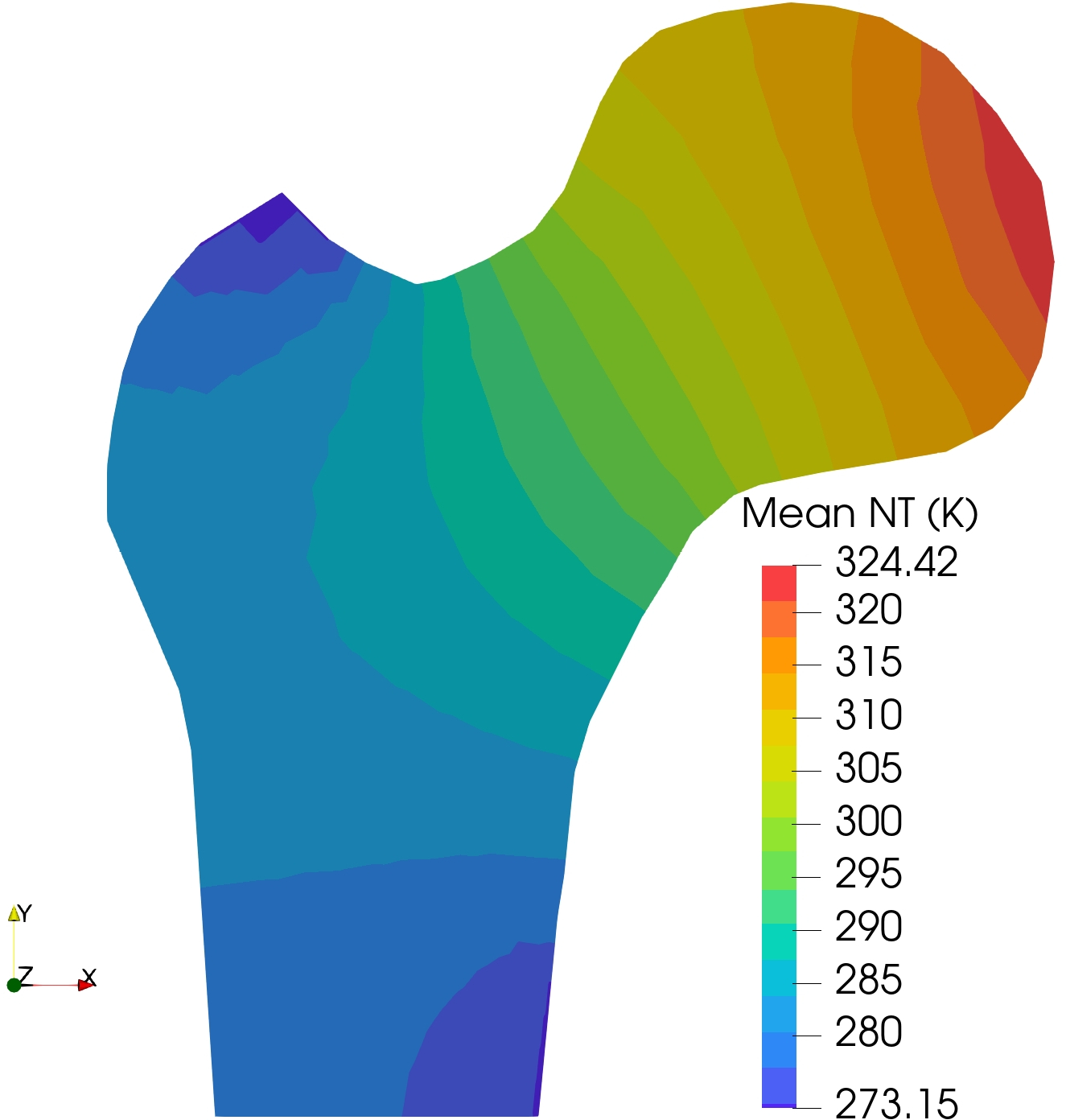}
		\caption{\acrshort{iso}-\acrshort{ortho}-\acrshort{scl}}
		\label{Fig:meanNTb}
	\end{subfigure}
%	\begin{subfigure}{0.49\textwidth}
%		\centering
%		\includegraphics[width=\widthtwoD\linewidth]{2D_scl_ortho_ortho_CCX_sideflux_1e5_mean_NT}
%		\caption{ortho-ortho-scl}
%	\end{subfigure}
	\begin{subfigure}[b]{0.32\textwidth}
		\centering
		\includegraphics[width=\widthtwoD\linewidth]{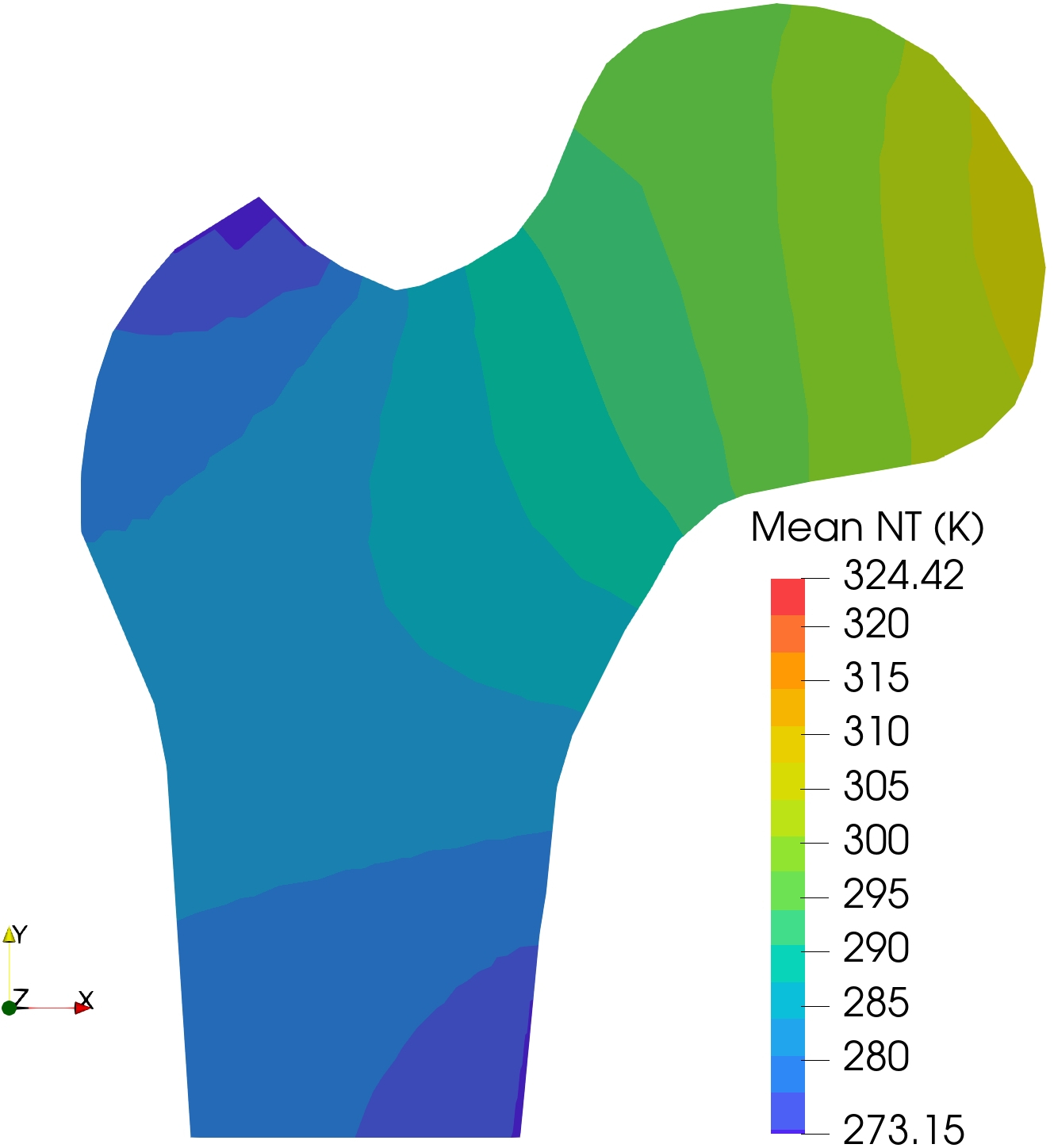}
		\caption{\acrshort{ortho}-\acrshort{ortho}-\acrshort{dir}}
		\label{Fig:meanNTc}
	\end{subfigure}
%	\begin{subfigure}{0.49\textwidth}
%		\centering
%		\includegraphics[width=\widthtwoD\linewidth]{2D_rotscl_ortho_ortho_CCX_sideflux_1e5_mean_NT}
%		\caption{ortho-ortho-rotscl}
%	\end{subfigure}
%	\begin{subfigure}{0.49\textwidth}
%		\centering
%		\includegraphics[width=\widthtwoD\linewidth]{2D_rotscl_iso_ortho_CCX_sideflux_1e5_mean_NT}
%		\caption{iso-ortho-rotscl}
%	\end{subfigure}
	\caption{\acrshort{mc} mean estimate of nodal temperature}
	\label{Fig:meanNT}
\end{figure}

\paragraph{Uncertainty quantification:}
\ffig{Fig:meanNT} shows the \acrshort{mc} mean estimate of the \acrfull{nt} $\stotemp$ for the three described stochastic models.  For easier interpretation, in this study, the results are displayed on a uniform scale. Hence, the values range from zero to a maximum value (corresponding to a given estimate of the desired output quantity). Clearly, the results in \ffigs{Fig:meanNTa}{Fig:meanNTb} are similar, as the considered reference tensor $\meanconducttensorscl$ is identical in both the cases. Due to the assumption of an orthotropic mean tensor $\meanconducttensorortho$,
in the third figure \ffig{Fig:meanNTc}, we see a significant difference in the magnitude/contour pattern of the temperature.

%------------------------------------------------------------------mean NT---------------------------------------------------------
\begin{figure}[h!]
	\begin{subfigure}[b]{0.32\textwidth}
		\centering
		\includegraphics[width=\widthtwoD\linewidth]{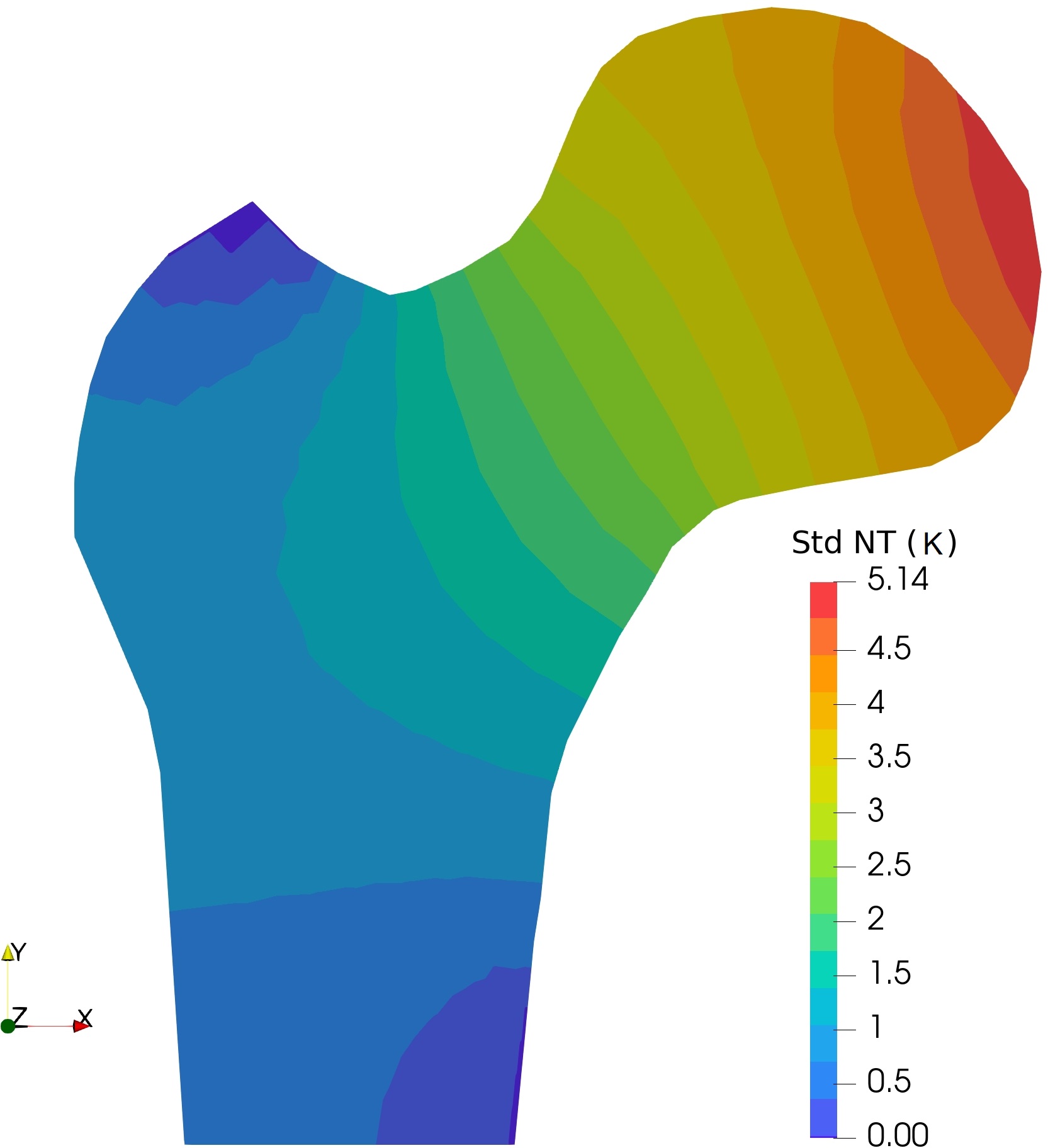}
		\caption{\acrshort{iso}-\acrshort{iso}-\acrshort{scl}}
		\label{Fig:varNTa}
	\end{subfigure}
	\begin{subfigure}[b]{0.32\textwidth}
		\centering
		\includegraphics[width=\widthtwoD\linewidth]{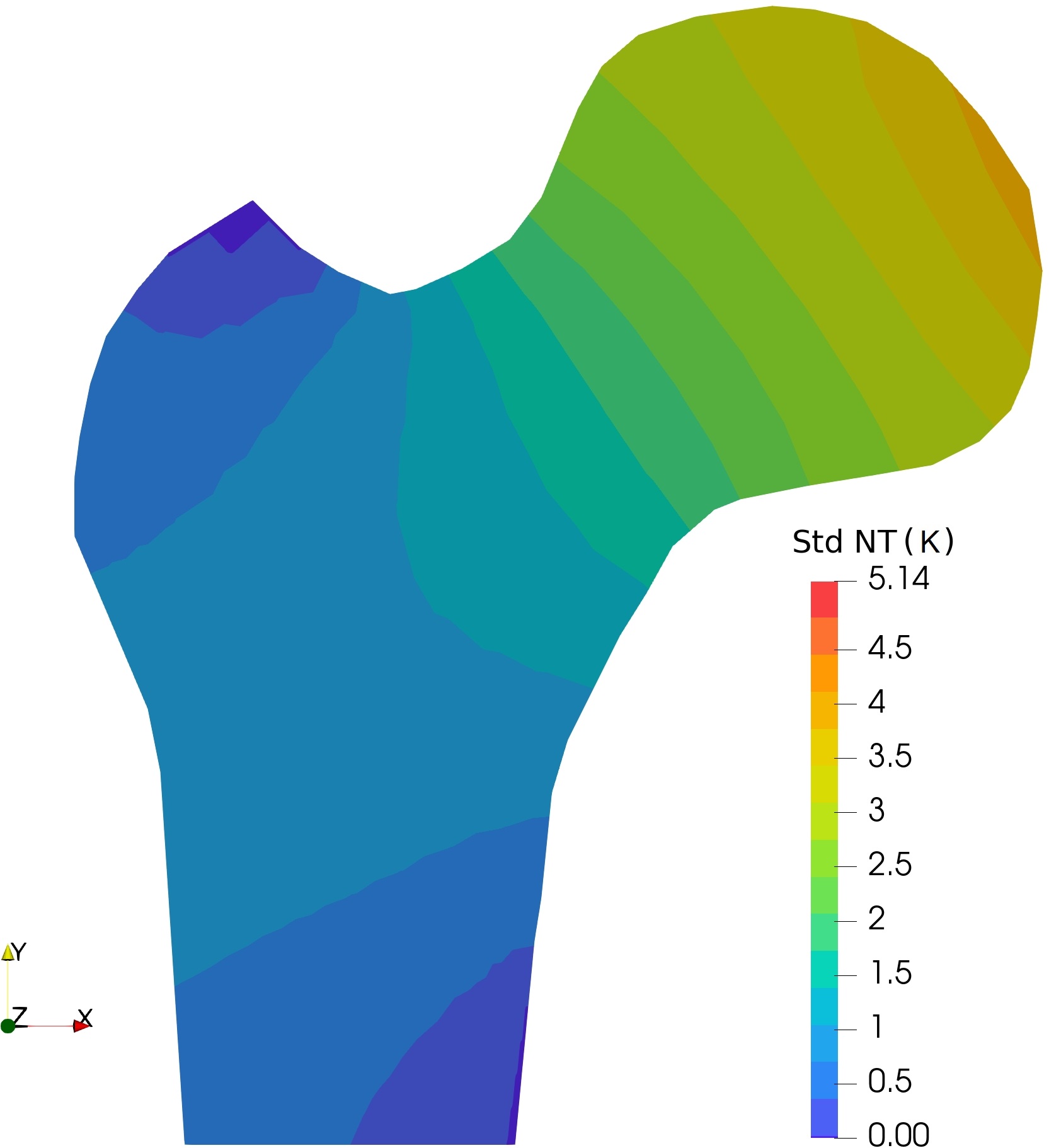}
		\caption{\acrshort{iso}-\acrshort{ortho}-\acrshort{scl}}
		\label{Fig:varNTb}
	\end{subfigure}
%	\begin{subfigure}{0.49\textwidth}
%		\centering
%		\includegraphics[width=\widthtwoD\linewidth]{2D_scl_ortho_ortho_CCX_sideflux_1e5_var_NT}
%		\caption{ortho-ortho-scl}
%	\end{subfigure}
	\begin{subfigure}[b]{0.32\textwidth}
		\centering
		\includegraphics[width=\widthtwoD\linewidth]{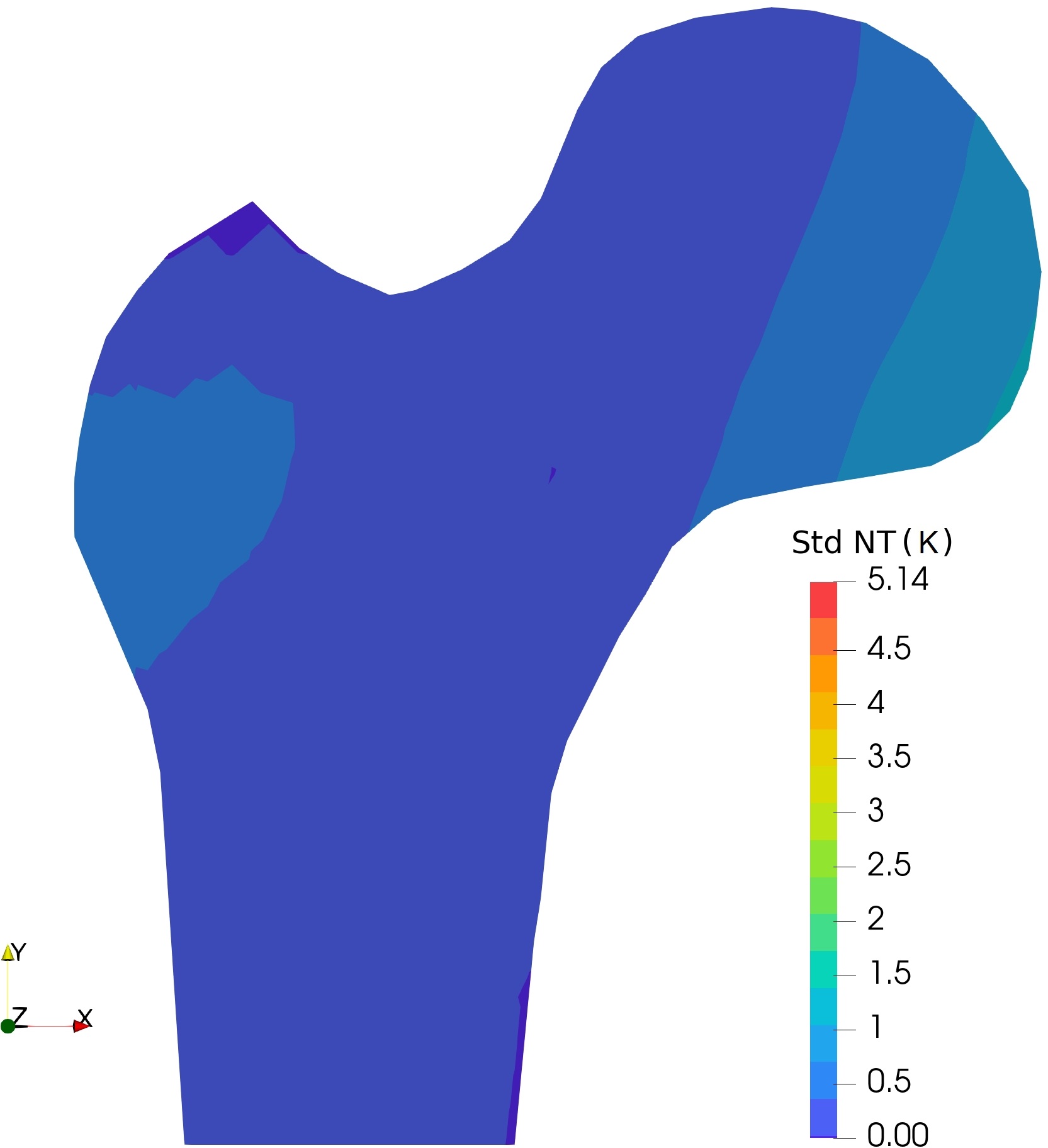}
		\caption{\acrshort{ortho}-\acrshort{ortho}-\acrshort{dir}}
		\label{Fig:varNTc}
	\end{subfigure}
%	\begin{subfigure}{0.49\textwidth}
%		\centering
%		\includegraphics[width=\widthtwoD\linewidth]{2D_rotscl_ortho_ortho_CCX_sideflux_1e5_var_NT}
%		\caption{ortho-ortho-rotscl}
%	\end{subfigure}
%	\begin{subfigure}{0.49\textwidth}
%		\centering
%		\includegraphics[width=\widthtwoD\linewidth]{2D_rotscl_iso_ortho_CCX_sideflux_1e5_var_NT}
%		\caption{iso-ortho-rotscl}
%	\end{subfigure}
	\caption{\acrshort{mc} standard deviation estimate of nodal temperature}
	\label{Fig:varNT}
\end{figure}

%----------------------------------------------------------var NT--------------------------------------------------------------------
\ffig{Fig:varNT} further depicts the standard deviation of temperature $\stotemp$. Due to different modelling assumptions, it is apparent that all three results in this figure are different from each other. But \ffig{Fig:varNTc} stands out the most, where the stochastic influence on nodal temperature is much lower as compared to the other two cases.

In the first two scenarios, only the scaling parameter is uncertain, whereas in the third example we consider only directional uncertainty.
We know that under the assumption of deterministic boundary conditions, the random temperature field $\temperature(x,\event)$ varies inversely to the stochastic coefficeint $\conducttensorrot(\event_r)$.  As the scaling parameters of the model $\conducttensorrot(\event_r)$ remain constant, the impact on the variation of temperature field $\temperature(x,\event)$ due to directional randomness is small. 
On the other hand, when the scaling parameters are modelled as varying, the standard deviation estimate of temperature $\stotemp$ becomes more sensitive as evident
%The higher sensitivity of standard deviation estimate of $\temperature(x,\event)$ on the scaling uncertainty in input model is certainly evident 
in \ffigs{Fig:varNTa}{Fig:varNTb}. Interestingly, when \ffig{Fig:varNTb} is compared to \ffig{Fig:varNTa}, it is clear that, varying the material symmetry from higher (isotropy) to lower (orthotropy) order in the model \acrshort{iso}-\acrshort{ortho}-\acrshort{scl} results in lower stochastic influence on temperature $\stotemp$.
%-------------------------------------------------------------mean HFL--------------------------------------------------------------
\begin{figure}[h!]
	\begin{subfigure}[b]{0.32\textwidth}
		\centering
		\includegraphics[width=\widthtwoD\linewidth]{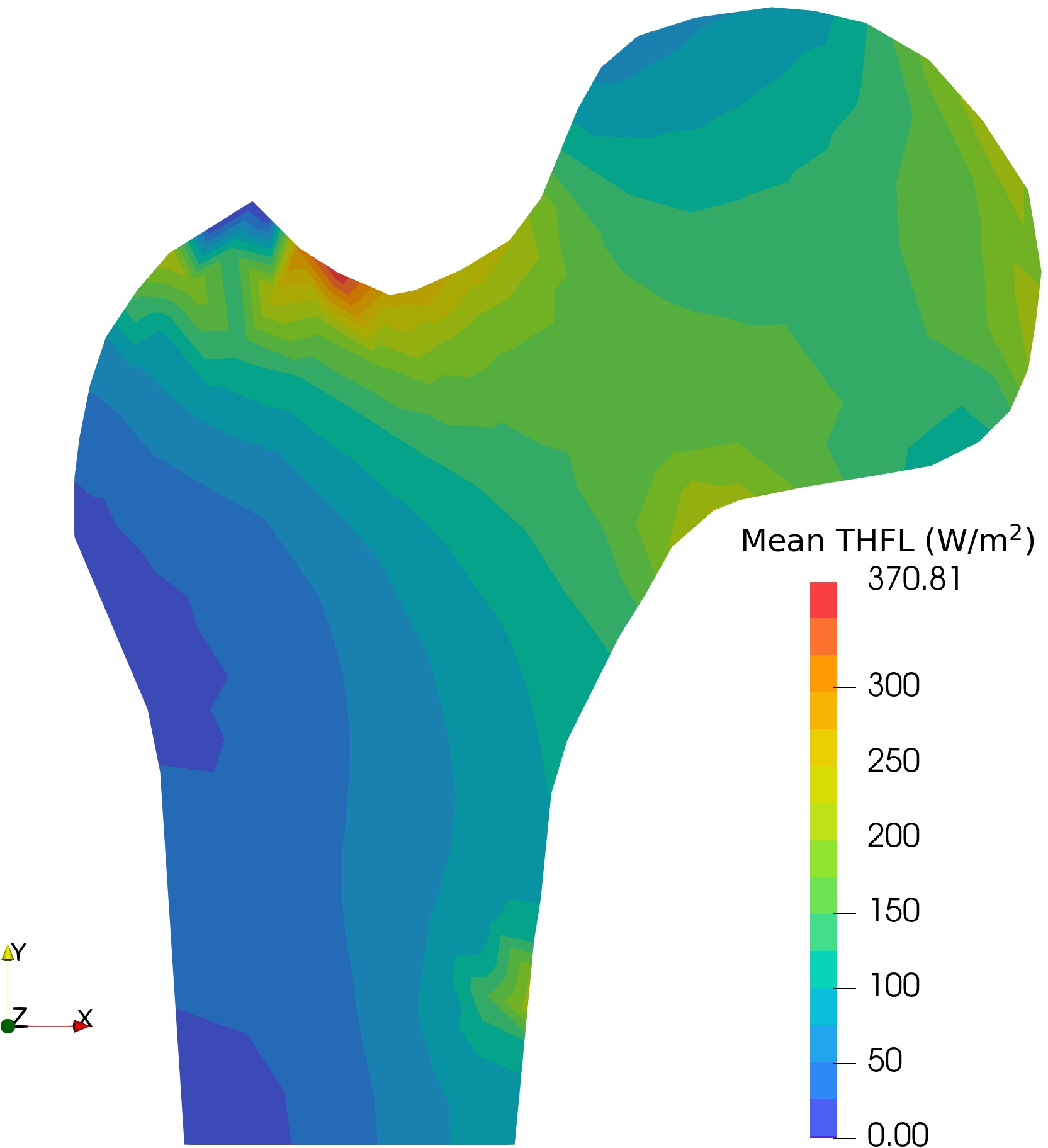}
		\caption{\acrshort{iso}-\acrshort{iso}-\acrshort{scl}}
		\label{Fig:meanHFLa}
	\end{subfigure}
	\begin{subfigure}[b]{0.32\textwidth}
		\centering
		\includegraphics[width=\widthtwoD\linewidth]{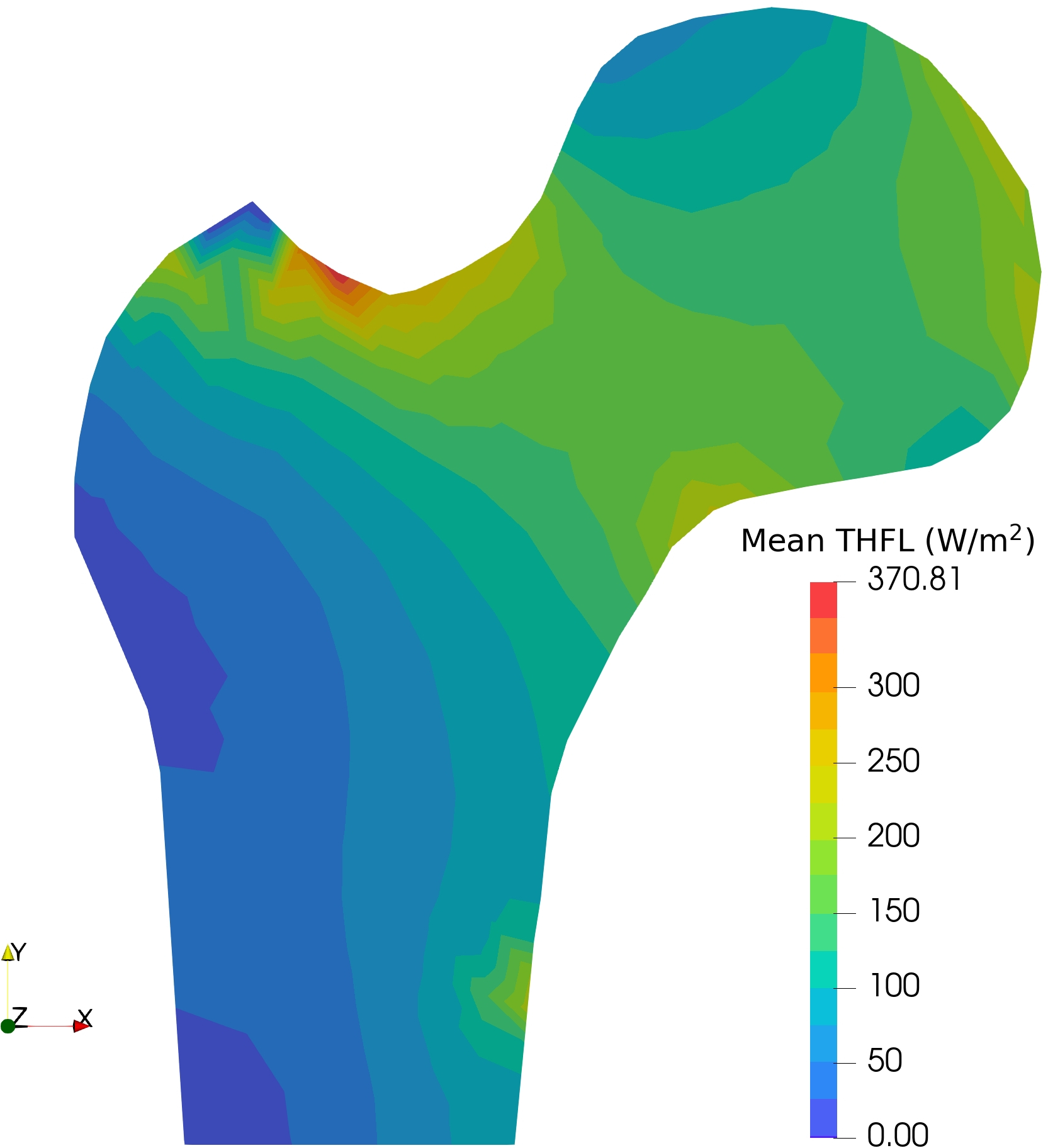}
		\caption{\acrshort{iso}-\acrshort{ortho}-\acrshort{scl}}
		\label{Fig:meanHFLb}
	\end{subfigure}
%	\begin{subfigure}{0.49\textwidth}
%		\centering
%		\includegraphics[width=\widthtwoD\linewidth]{2D_scl_ortho_ortho_CCX_sideflux_1e5_mean_HFL}
%		\caption{ortho-ortho-scl}
%	\end{subfigure}
	\begin{subfigure}[b]{0.32\textwidth}
		\centering
		\includegraphics[width=\widthtwoD\linewidth]{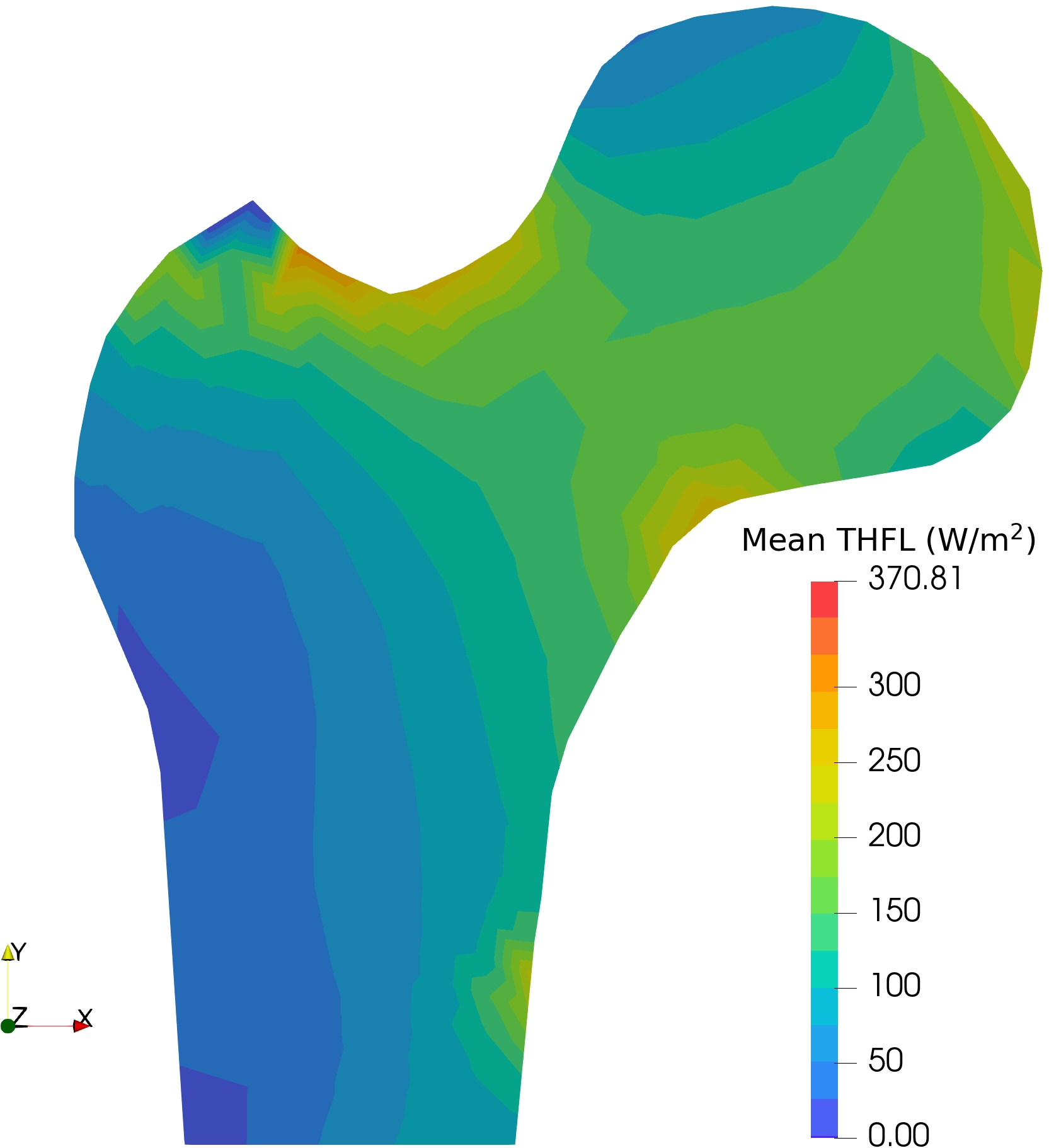}
		\caption{\acrshort{ortho}-\acrshort{ortho}-\acrshort{dir}}
		\label{Fig:meanHFLc}
	\end{subfigure}
%	\begin{subfigure}{0.49\textwidth}
%		\centering
%		\includegraphics[width=\widthtwoD\linewidth]{2D_rotscl_ortho_ortho_CCX_sideflux_1e5_mean_HFL}
%		\caption{ortho-ortho-rotscl}
%	\end{subfigure}
%	\begin{subfigure}{0.49\textwidth}
%		\centering
%		\includegraphics[width=\widthtwoD\linewidth]{2D_rotscl_iso_ortho_CCX_sideflux_1e5_mean_HFL}
%		\caption{iso-ortho-rotscl}
%	\end{subfigure}
	\caption{\acrshort{mc} mean estimate of total heat flux (Euclidean norm)}
	\label{Fig:meanHFL}
\end{figure}

\begin{figure}[h!]
	\begin{subfigure}[b]{0.32\textwidth}
		\centering
		\includegraphics[width=\widthtwoD\linewidth]{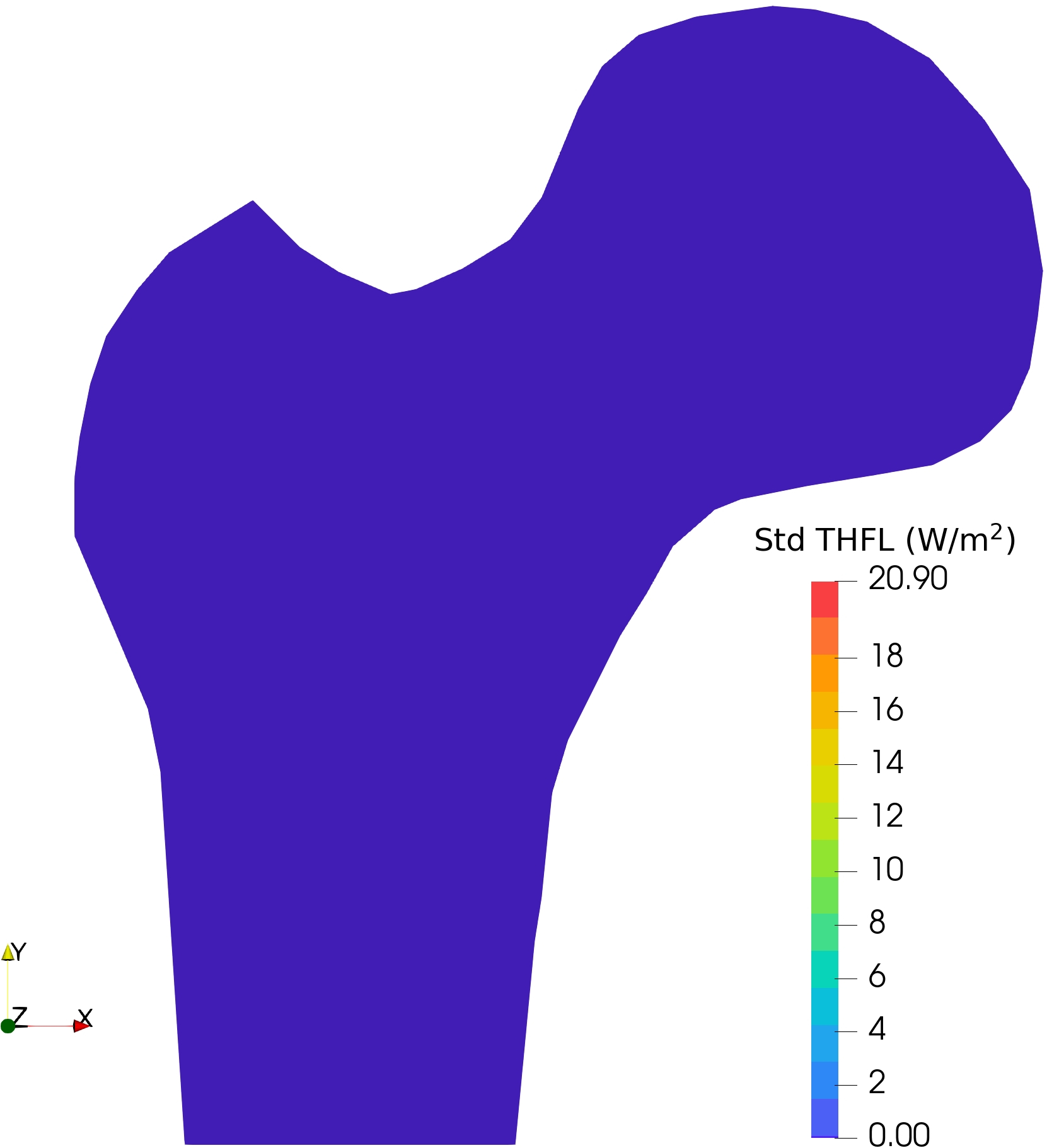}
		\caption{\acrshort{iso}-\acrshort{iso}-\acrshort{scl}}
		\label{Fig:varHFLa}
	\end{subfigure}
	\begin{subfigure}[b]{0.32\textwidth}
		\centering
		\includegraphics[width=\widthtwoD\linewidth]{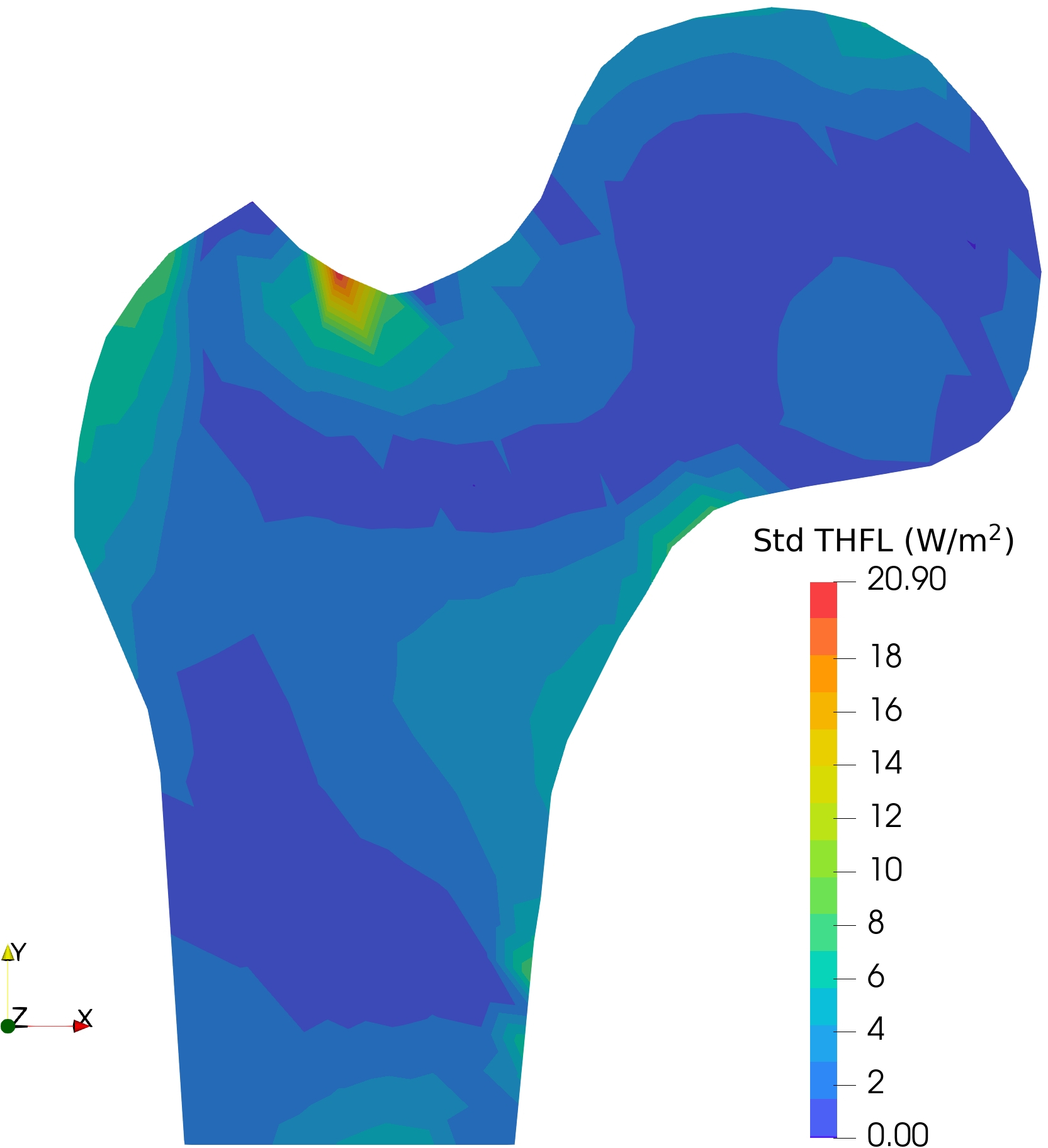}
		\caption{\acrshort{iso}-\acrshort{ortho}-\acrshort{scl}}
		\label{Fig:varHFLb}
	\end{subfigure}
%	\begin{subfigure}{0.49\textwidth}
%		\centering
%		\includegraphics[width=\widthtwoD\linewidth]{2D_scl_ortho_ortho_CCX_sideflux_1e5_var_HFL}
%		\caption{ortho-ortho-scl}
%	\end{subfigure}
	\begin{subfigure}[b]{0.32\textwidth}
		\centering
		\includegraphics[width=\widthtwoD\linewidth]{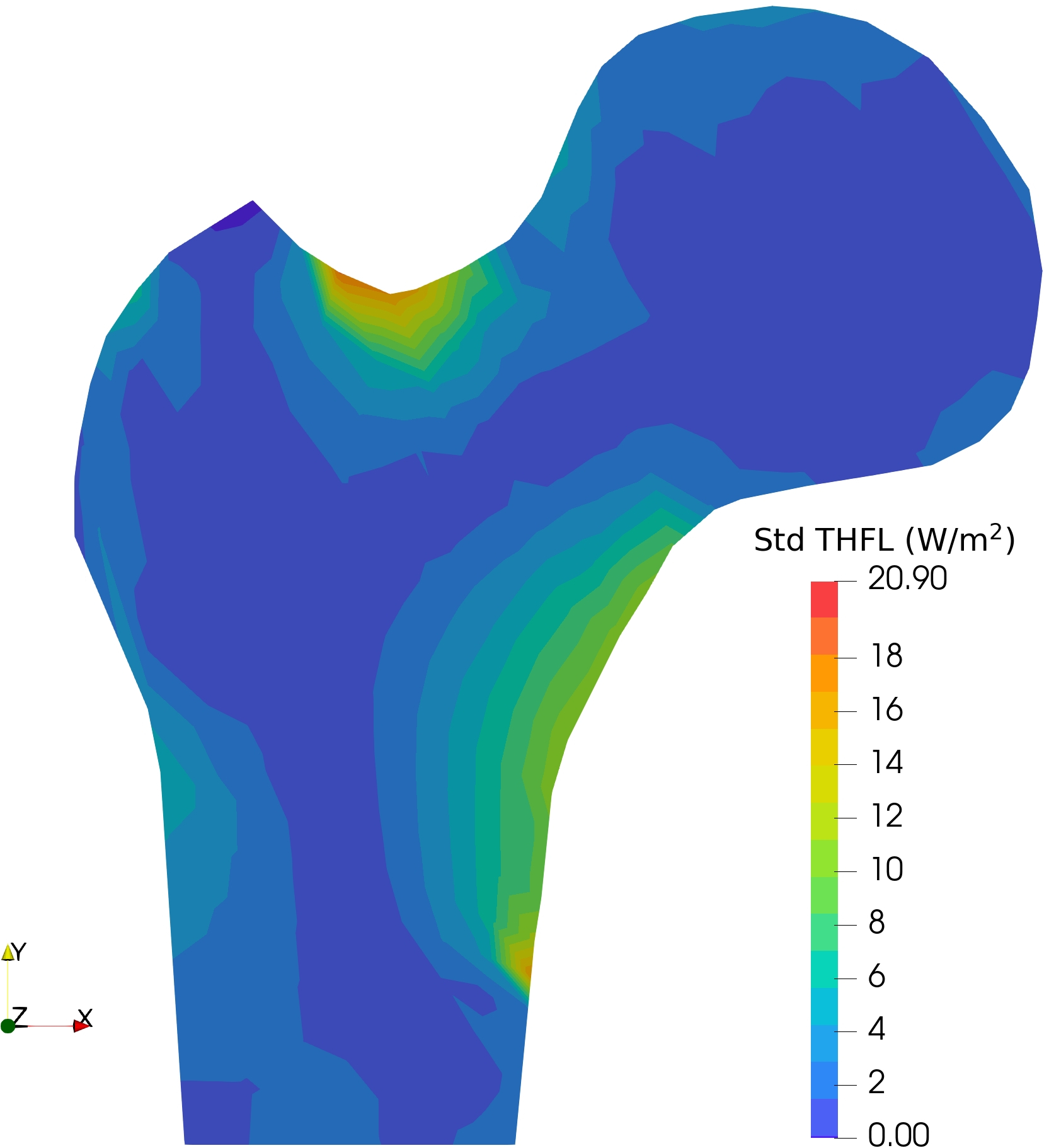}
		\caption{\acrshort{ortho}-\acrshort{ortho}-\acrshort{dir}}
		\label{Fig:varHFLc}
	\end{subfigure}
%	\begin{subfigure}{0.49\textwidth}
%		\centering
%		\includegraphics[width=\widthtwoD\linewidth]{2D_rotscl_ortho_ortho_CCX_sideflux_1e5_var_HFL}
%		\caption{ortho-ortho-rotscl}
%	\end{subfigure}
%	\begin{subfigure}{0.49\textwidth}
%		\centering
%		\includegraphics[width=\widthtwoD\linewidth]{2D_rotscl_iso_ortho_CCX_sideflux_1e5_var_HFL}
%		\caption{iso-ortho-rotscl}
%	\end{subfigure}
	\caption{\acrshort{mc} standard deviation estimate of total heat flux (Euclidean norm)}
	\label{Fig:varHFL}
\end{figure}

Furthermore, the sample mean estimates of \acrfull{thfl} $\stototalflux$ are presented in \ffig{Fig:meanHFL}, where, \ffigs{Fig:meanHFLa}{Fig:meanHFLb} showcase similar results as described previously (see \ffigs{Fig:meanNTa}{Fig:meanNTb}). However, in comparison to these two figures, a visible difference in maximum total heat flux and contour pattern is noticed in \ffig{Fig:meanHFLc}.  In \ffig{Fig:varHFL}, the corresponding estimated standard deviation of THFL $\stototalflux$ is plotted. Here, the most interesting aspect is seen in \ffig{Fig:varHFLa}, where, the standard deviation is close to zero. It turns out that, in the model \acrshort{iso}-\acrshort{iso}-\acrshort{scl}, the stochastic coefficient $\conducttensorscl(\event_s)$ has almost perfect inverse correlation with temperature gradient field $\nabla\temperature(\vek{x},\event)$. 
Thus, the randomness in the model $\conducttensorscl(\event_s)$ has almost no stochastic impact on THFL $\stototalflux$. 
However, on the contrary, the estimate in \ffigs{Fig:varHFLb}{Fig:varHFLc} is visible, signifying the stochastic influence of the considered input models on THFL $\stototalflux$.
One may conclude that directivity has more impact on the heat flux than the scaling parameter. 
%-------------------------------------------------------------var HFL----------------------------------------------------------------
%-----------------------------------------------------------mean var dir-------------------------------------------------------------

\begin{figure}[h!]
	\begin{subfigure}[b]{0.32\textwidth}
		\centering
		\includegraphics[width=\widthtwoD\linewidth]{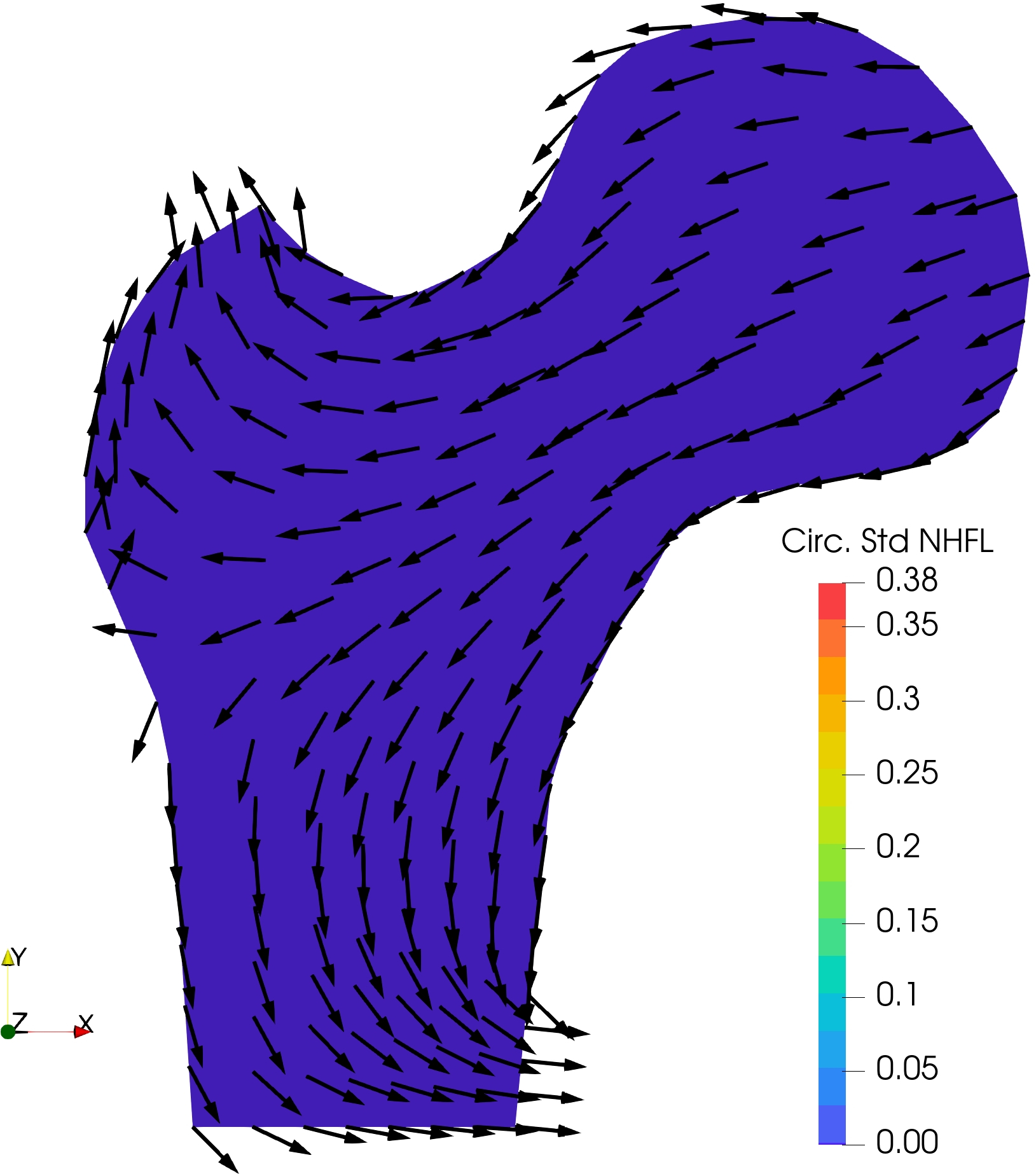}
		\caption{\acrshort{iso}-\acrshort{iso}-\acrshort{scl}}
		\label{Fig:circHFLa}
	\end{subfigure}
	\begin{subfigure}[b]{0.32\textwidth}
		\centering
		\includegraphics[width=\widthtwoD\linewidth]{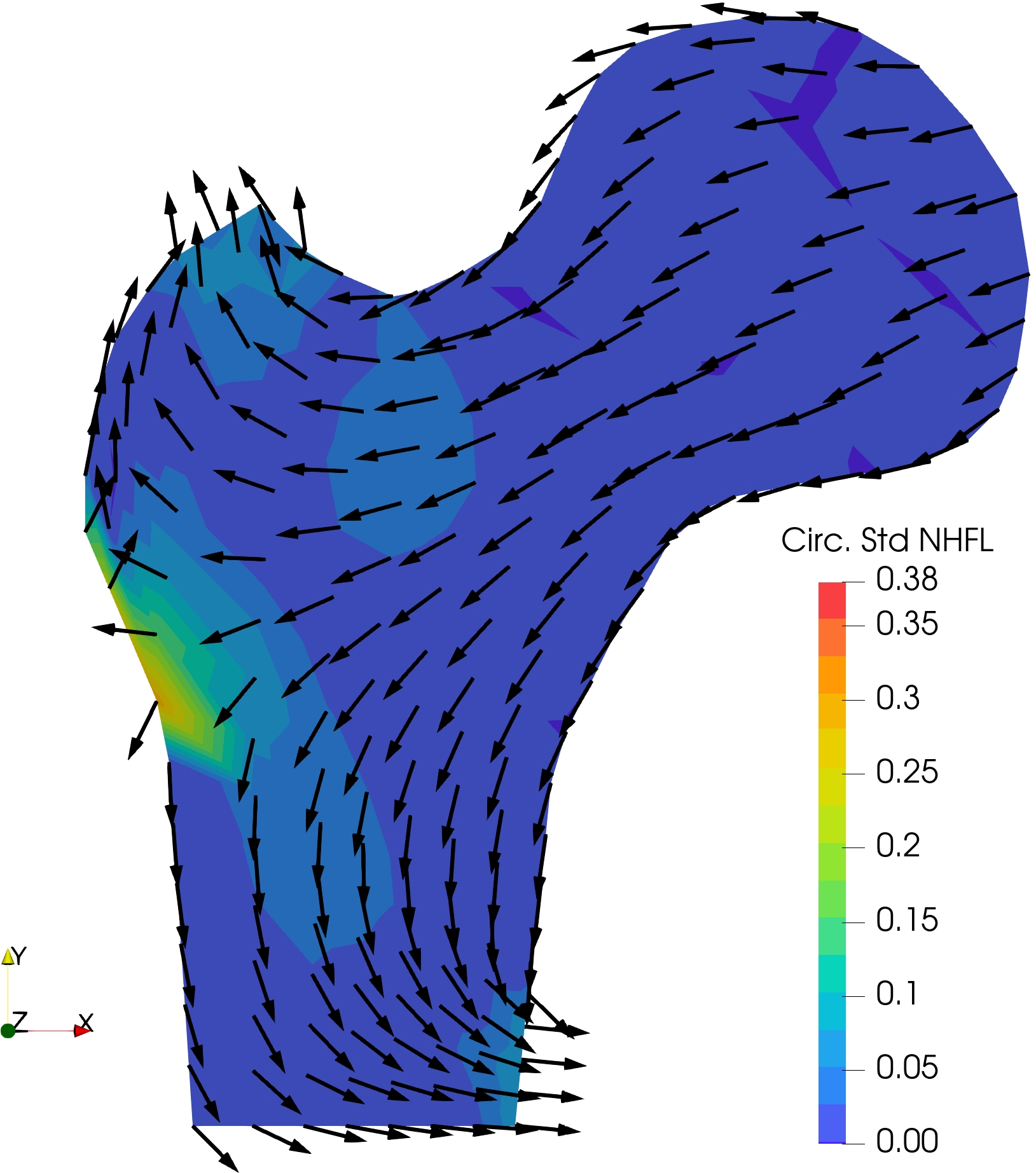}
		\caption{\acrshort{iso}-\acrshort{ortho}-\acrshort{scl}}
		\label{Fig:circHFLb}
	\end{subfigure}
%	\begin{subfigure}{0.49\textwidth}
%		\centering
%		\includegraphics[width=\widthtwoD\linewidth]{2D_scl_ortho_ortho_CCX_sideflux_1e5_mean_var_dir}
%		\caption{ortho-ortho-scl}
%	\end{subfigure}
	\begin{subfigure}[b]{0.32\textwidth}
		\centering
		\includegraphics[width=\widthtwoD\linewidth]{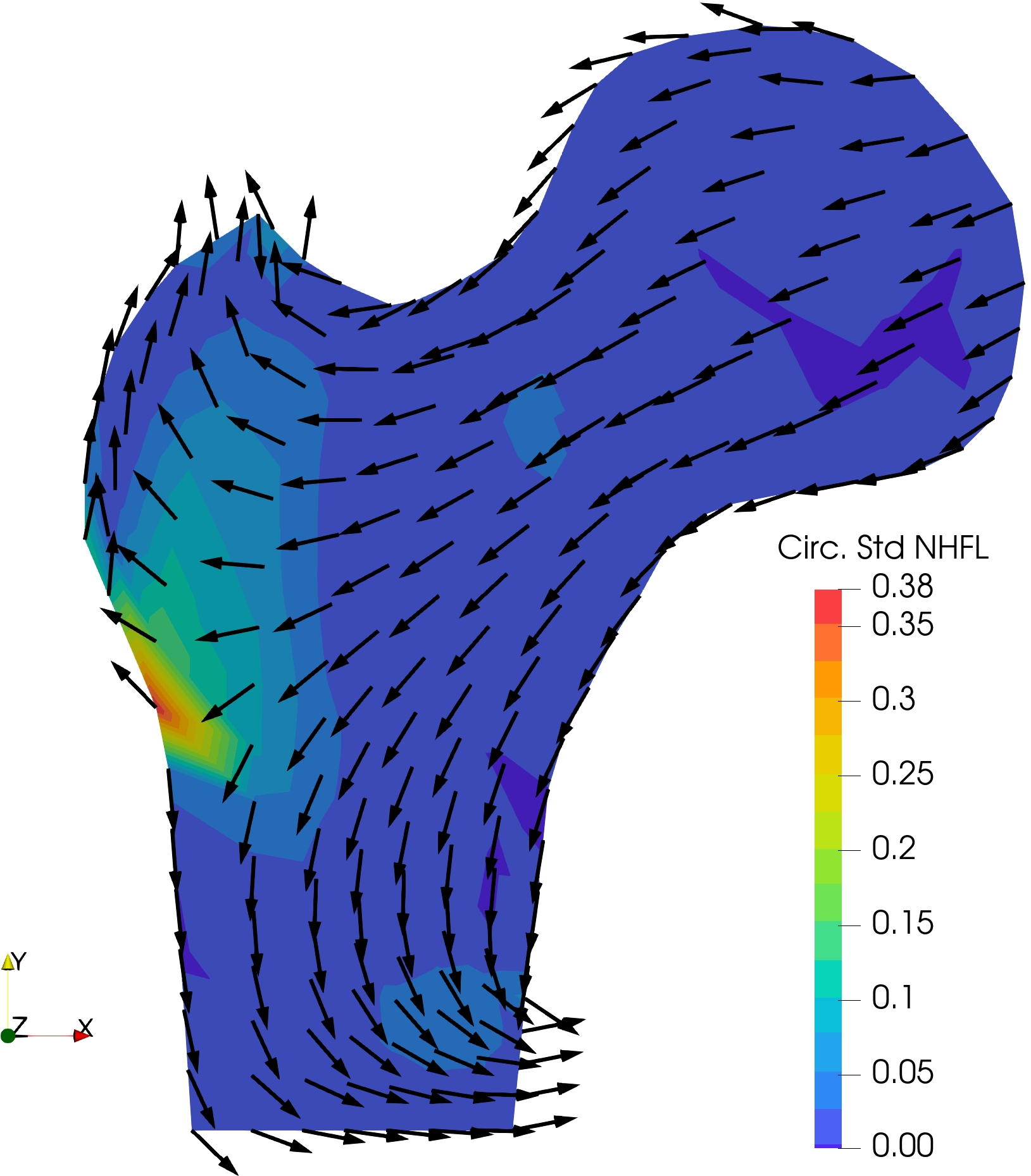}
		\caption{\acrshort{ortho}-\acrshort{ortho}-\acrshort{dir}}
		\label{Fig:circHFLc}
	\end{subfigure}
%	\begin{subfigure}{0.49\textwidth}
%		\centering
%		\includegraphics[width=\widthtwoD\linewidth]{2D_rotscl_ortho_ortho_CCX_sideflux_1e5_mean_var_dir}
%		\caption{ortho-ortho-rotscl}
%	\end{subfigure}
%	\begin{subfigure}{0.49\textwidth}
%		\centering
%		\includegraphics[width=\widthtwoD\linewidth]{2D_rotscl_iso_ortho_CCX_sideflux_1e5_mean_var_dir}
%		\caption{iso-ortho-rotscl}
%	\end{subfigure}
	\caption{\acrshort{mc} circular mean (vector representation) and circular standard deviation estimates (dimensionless quantity) of normalized heat flux}
	\label{Fig:circHFL}
\end{figure}

Additionally, the directional mean $\mucirc^{MC}(\stofluxunit)$ and standard deviation $\stdMCcirc(\stofluxunit)$ of \acrfull{nhfl} $\stofluxunit$---from 
%\feqs{Eq:normheatflux}{Eq:meannormheatflux}---are 
\feq{Eq:meannormheatflux}---are 
displayed in \ffig{Fig:circHFL}.
The mean orientation (vector representation) in
\ffigs{Fig:circHFLa}{Fig:circHFLb} look similar, however, with a closer inspection of \ffig{Fig:circHFLc},
the difference is apparent. Similar to the results in \ffig{Fig:varHFL}, the circular standard deviation estimate in \ffig{Fig:circHFLa} is also close to zero, showing once again the insensitivity of directional attribute of NHFL $\stofluxunit$ to scaling uncertainty present in the model $\conducttensorscl(\event_s)$. 
Also, in comparison, the standard deviation estimate in \ffig{Fig:circHFLc} is more significant than in \ffig{Fig:circHFLb}.

%% file: results_input_3D.tex
% !TEX root = ../paper_tensors.tex
% !TEX encoding = UTF-8 Unicode

% \begin{wrapfigure}[16]{r}{0.39\textwidth}
% 	\centering
% 	\includegraphics[width=0.7\linewidth]{3D_femur_BC_thermal_sideflux_1.eps}
% 	\captionof{figure}{Boundary conditions of 3D femur bone}
% 	\label{Fig:BC3D}	
% \end{wrapfigure}

\begin{figure}[ht]
	\centering
	\includegraphics[width=0.25\linewidth]{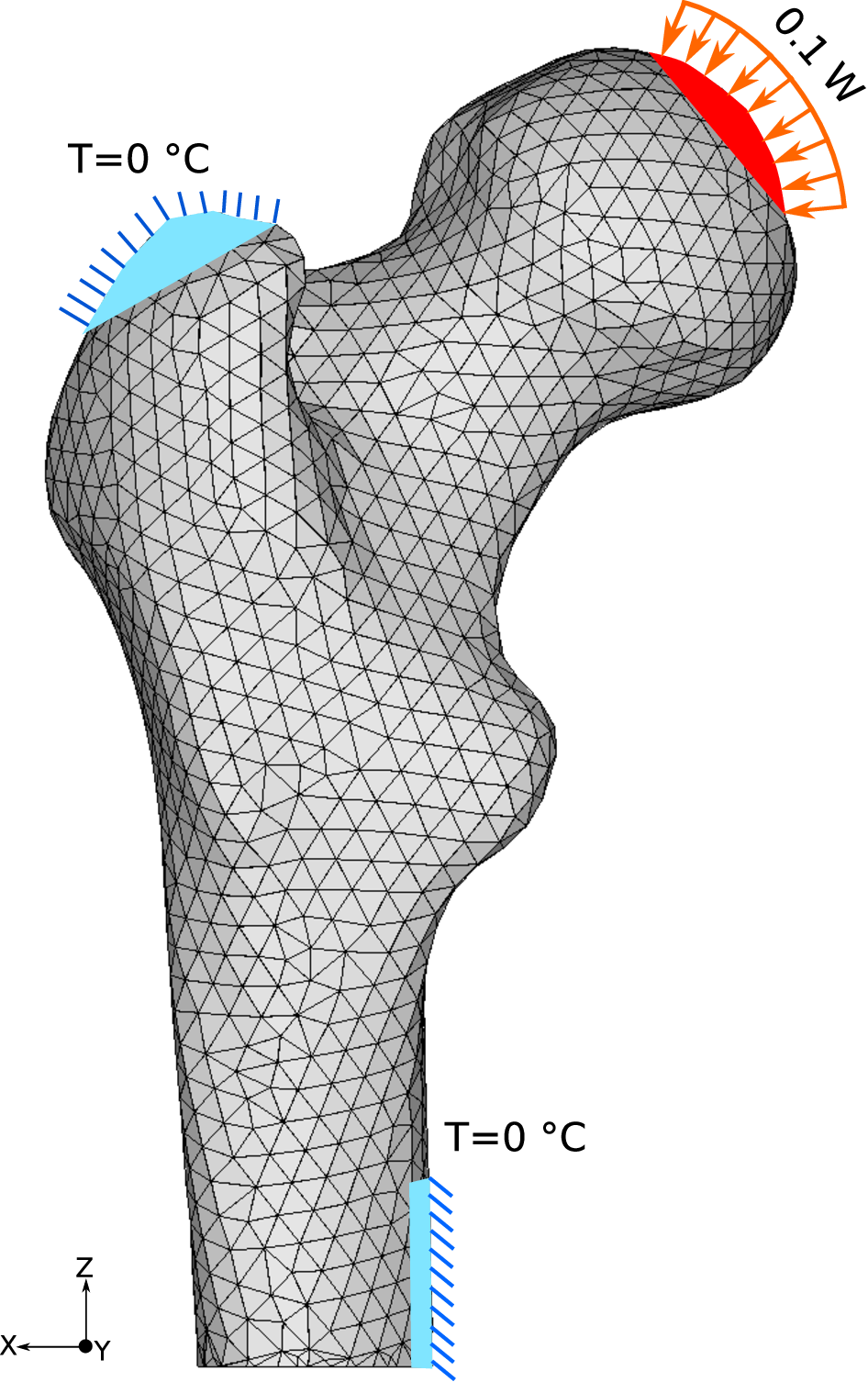}
	\caption{Boundary conditions of 3D femur bone}
	\label{Fig:BC3D}	
\end{figure}

A three-dimensional proximal femur configuration of dimensions 45 mm in width and 154 mm in height is considered. The boundary conditions with identical values as used in the 2D example are imposed, shown in \ffig{Fig:BC3D}. The computational domain is discretized by a four-noded tetrahedral finite element mesh, comprising 12504 nodes, 3166 elements and 35115 degrees of freedom.
% , the deterministic simulation is performed by the Calculix solver \cite{Dhondt_CCX_2004}. 
We use similar modelling scenarios as in 2D case i.e. iso-iso-scl, iso-ortho-scl and ortho-ortho-dir (see \ftbl{Tbl:abbreviation} for model abbreviations). 
The considered isotropic $\meanconducttensorscl\in\PosSymSecRankSpacethreeD$ and orthotropic $\meanconducttensorortho\in\PosSymSecRankSpacethreeD$ reference conductivity tensors are tabulated in \ftbl{Tbl:meantensor3D}. 
% The values of thermal conductivity of the model $\meanconducttensorortho$ in x, y and z-direction are $\meaneigvalroteleone=0.54$ W/mK, $\meaneigvalroteletwo=0.75$ W/mK and $\meaneigvalrotelethree=1$ W/mK, respectively. 
Furthermore, by fixing the negative y-axis as the rotational axis, the directional vectors of the model $\meanconducttensorortho$ corresponding to the eigenvalues $\meaneigvalroteleone$ and $\meaneigvalrotelethree$ are positioned at an angle of $\pi/4$ in an anti-clockwise direction with respect to the $x$ and $z$-axis, respectively.
% \begin{table}[!h]
% \centering
% 	\begin{tabular}{lll}
% 		\specialrule{1pt}{1pt}{1pt}			
% 		Elements & Nodes & \begin{tabular}[c]{@{}l@{}}
% 			Degrees of \\ freedom 
% 		\end{tabular} \\ \midrule
% 		3166 & 12504 &  35115 \\  
% 		\specialrule{1pt}{1pt}{1pt}
% 	\end{tabular}
% 	\caption{Mesh specifications}
% 	\label{Tbl:meshspec3D}
% \end{table}
\begin{table}[h!]
	\centering
	\begin{tabular}{ll}
	\specialrule{1pt}{1pt}{1pt}
	\begin{tabular}[c]{@{}l@{}}isotropy\\ (W/mK)\end{tabular} & \begin{tabular}[c]{@{}l@{}}orthotropy\\ (W/mK)\end{tabular} \\ 
	$\meanconducttensorscl$ & $\meanconducttensorortho$\\ \midrule
	$
		\begin{bmatrix}
		0.54 & 0 & 0  \\ 0 & 0.54 & 0 \\ 0 & 0 & 0.54
		\end{bmatrix} $  &  $  		\begin{bmatrix}
		0.77 & 0 & 0.23  \\ 0 & 0.75 & 0 \\ 0.23 & 0 & 0.77
		\end{bmatrix} $ \\ \midrule
	\multicolumn{2}{c}{Eigenvalues (W/mK)} \\  \midrule
	\begin{tabular}[c]{@{}l@{}}$\meaneigvalscleleone=0.54$, $\meaneigvalscleletwo=0.54$\\ $\meaneigvalsclelethree=0.54$\end{tabular} & \begin{tabular}[c]{@{}l@{}}$\meaneigvalroteleone=0.54$, $\meaneigvalroteletwo=0.75$\\ $\meaneigvalrotelethree=1$\end{tabular} \\
	\specialrule{1pt}{1pt}{1pt}            
	\end{tabular}
\caption{Reference conductivity tensors with corresponding eigenvalues of 3D femur\label{Tbl:meantensor3D}}
\end{table}

%=========================================================================================================================
\begin{figure}[!h]
	\begin{subfigure}{\textwidth}
		\centering
		\input{3D_LN_iso_iso}
		\caption{}
		\label{Fig:isoisoLN3D}
	\end{subfigure}	
	\begin{subfigure}{\textwidth}
		\centering
		\includegraphics[width=0.85\linewidth]{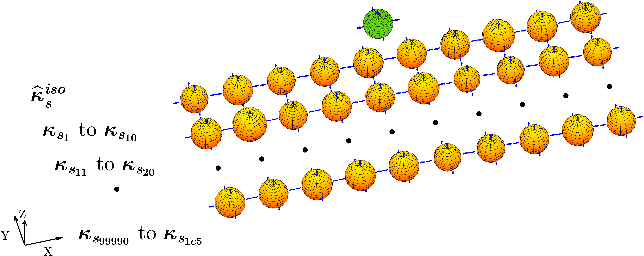}
		\caption{}
		\label{Fig:isoisoquad3D}
	\end{subfigure}	
	\caption{Stochastic tensor $\conducttensorscl(\event_s)$ with fixed symmetry (\acrshort{iso}-\acrshort{iso}-\acrshort{scl}): \subref{Fig:isoisoLN} log-normal \acrshort{pdf} of identical and dependent random scaling values $\eigvalscleleone(\event_s)=\eigvalscleletwo(\event_s)=\eigvalsclelethree(\event_s)$, \subref{Fig:isoisoquad} geometric visualization of reference tensor $\meanconducttensorscl$ and realisations of random tensor $\conducttensorscl(\event_s)$; the geometries are enhanced by scaling the eigenvalues by a factor of 0.5} 
	\label{Fig:isoisoscl3D}
\end{figure}

\paragraph{Random scaling with fixed symmetry:}
A stochastic tensor $\conducttensorscl(\event_s)\in\PosSymSecRankSpacethreeD$
% , which is equivalent to the exponential of $\log\conducttensorscl(\event_s)$ in \feq{Eq:rndscl} 
with fixed symmetry (\acrshort{iso}-\acrshort{iso}-\acrshort{scl}), where the random scaling elements 
$\eigvalscleleone(\event_s)=\eigvalscleletwo(\event_s)=\eigvalsclelethree(\event_s)$ are 
modelled by an identical log-normal random variable, is constructed.
%(see \feq{Eq:rndscl} for the exponential mapping). 
The schematic representation of the model is shown in \ffig{Fig:isoisoscl3D}, where 
\ffig{Fig:isoisoLN3D} presents the log-normal PDF of scaling parameters 
$\{\eigvalsclele\}^3_{i=1}$, whereas the geometric visualization---in the form of equivalent 
spheres---of the reference tensor $\meanconducttensorscl$ and realisations 
$\{\conducttensor_{s_i}\}^\samplesize_{i=1}$ is displayed in \ffig{Fig:isoisoquad3D}. 
Here the radius of the sphere is determined by ${1}/{\sqrt{\eigvalscleleone(\event_s)}}$, 
${1}/{\sqrt{\eigvalscleletwo(\event_s)}}$ and ${1}/{\sqrt{\eigvalsclelethree(\event_s)}}$. 
As can be seen that the radius of the spheres in the realisations vary and the spherical 
shape of the mean is maintained in the realisations, owing to the preservation of 
isotropic material symmetry in the mean and realisations.

%=========================================================================================================================
\begin{figure}[!h]
	\begin{subfigure}{\textwidth}
		\centering
		\input{3D_LN_iso_ortho}
		\caption{}
		\label{Fig:isoorthoLN3D}
		\end{subfigure}
	\begin{subfigure}{\textwidth}
		\centering
		\includegraphics[width=0.85\linewidth]{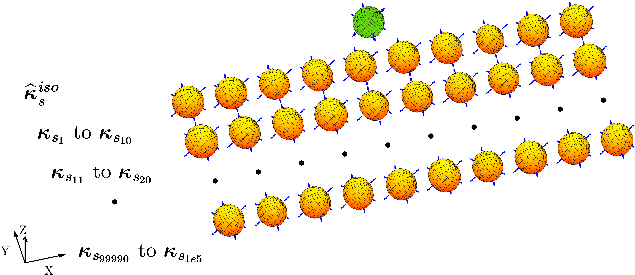} %tex files are too big, hence .eps
		\caption{}
		\label{Fig:isoorthoquad3D}
	\end{subfigure}	
	\caption{Stochastic tensor $\conducttensorscl(\event_s)$ with varying symmetry (\acrshort{iso}-\acrshort{ortho}-\acrshort{scl}): \subref{Fig:isoorthoLN} log-normal \acrshort{pdf} of i.i.d random scaling values $\{\eigvalsclele\}^3_{i=1}$, \subref{Fig:isoorthoquad} geometric visualization of reference tensor $\meanconducttensorscl$ and realisations of random tensor $\conducttensorscl(\event_s)$; the geometries are enhanced by scaling the eigenvalues by a factor of 0.5} 
	\label{Fig:isoorthoscl3D}
\end{figure}

\paragraph{Random scaling with varying symmetry:}
a random tensor model $\conducttensorscl(\event_s)\in\PosSymSecRankSpacethreeD$ similar to the previous scenario, however with a varying material symmetry (\acrshort{iso}-\acrshort{ortho}-\acrshort{scl}) is showcased in \ffig{Fig:isoorthoscl3D}. In \ffig{Fig:isoorthoLN3D}, the PDF of i.i.d log-normal scaling parameters $\{\eigvalsclele\}^3_{i=1}$ can be seen, whereas 
\ffig{Fig:isoorthoquad3D} shows the geometric interpretation of the mean $\meanconducttensorscl$ and realisations $\{\conducttensor_{s_i}\}^\samplesize_{i=1}$. 
Looking from the negative y-axis, the eigenvectors of the model $\conducttensorscl(\event_s)$ conforming to the eigenvalues $\meaneigvalscleleone$ and $\meaneigvalsclelethree$ are at an angle of $\pi/4$ in an anti-clockwise direction corresponding to x and z-axis, respectively.
A change in symmetry from isotropic in the mean (spherical shape) to orthotropic (ellipsoidal shape) in the realisations is evident. Also, the fluctuation in lengths of semi-major axes---calculated by ${1}/{\sqrt{\eigvalscleleone(\event_s)}}$, ${1}/{\sqrt{\eigvalscleletwo(\event_s)}}$ and ${1}/{\sqrt{\eigvalsclelethree(\event_s)}}$---of the ellipsoids represents the aspect of varying scaling parameters. 

% with mean $\meanconducttensorortho$
\begin{figure}[!h]
	\begin{subfigure}{\textwidth}
		\centering
		\includegraphics[width=0.37\linewidth]{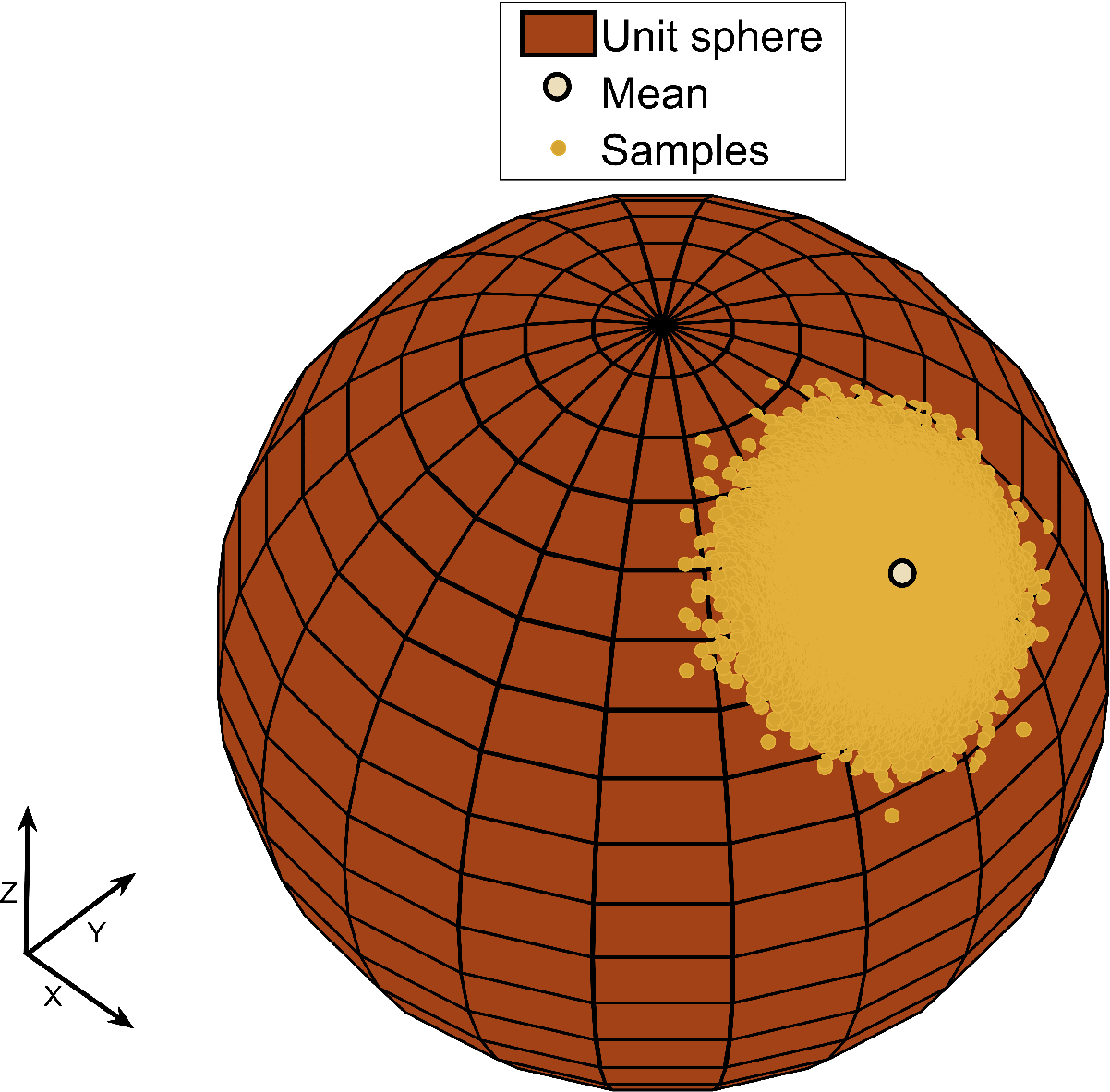}
		\caption{}
		\label{Fig:orthoorthoVMunit3D}
	\end{subfigure}
	\begin{subfigure}{\textwidth}
		\centering
		\includegraphics[width=0.85\linewidth]{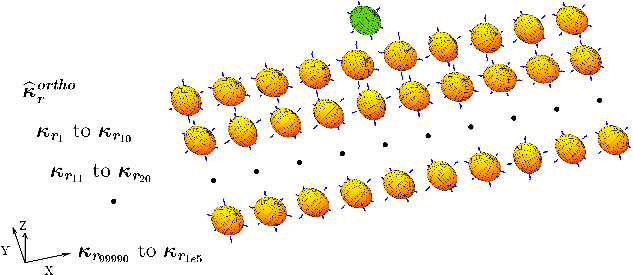} %tex files are too big, hence .eps
		\caption{}
		\label{Fig:orthoorthorot3D}
	\end{subfigure}
	\caption{Stochastic tensor $\conducttensorrot(\event_r)$ with fixed symmetry (\acrshort{ortho}-\acrshort{ortho}-\acrshort{dir}): \subref{Fig:orthoorthoVMunit3D}  realisations of von Mises Fisher random variable $\rotvect(\event_r)$ on unit sphere with mean $\mucirc$, \subref{Fig:orthoorthorot3D}  geometric visualization of reference tensor $\meanconducttensorortho$ and realisations of random tensor $\conducttensorrot(\event_r)$; the geometries are enhanced by scaling the eigenvalues by a factor of 0.5, and the orientations with respect to the mean are scaled by a factor of 1.5}
	\label{Fig:orthoortho3D}
\end{figure}

\textbf{Random orientation with fixed symmetry}:
to account for directional uncertainty only, a random tensor with fixed symmetry (\acrshort{ortho}-\acrshort{ortho}-\acrshort{dir}) $\conducttensorrot(\event_r)\in\PosSymSecRankSpacethreeD$, which is the exponential of the model described in \feq{Eq:stomodrot}, is considered. 
The eigenvector $\rotvect(\event_r)$ of random tensor $\conducttensorrot(\event_r)$ related to the eigenvalue $\meaneigvalsclelethree$ is modelled as a von Mises-Fisher distribution (see \feq{Eq:VMFPDF}) with
the mean direction $\mucirc$ (the eigenvector of the reference tensor $\meanconducttensorortho$ which corresponds to the eigenvalue $\meaneigvalsclelethree$), and the concentration parameter $\concparam=75$. 
% then the eigenvector of random tensor $\conducttensorrot(\event_r)$ related to the eigenvalue $\meaneigvalsclelethree$ is modelled as a von Mises-Fisher distribution (see \feq{Eq:VMFPDF})--let the random eigenvector be denoted by $\rotvect(\event_r)$. 
The random rotational angle $\rotangle(\event_r)$ is determined using the dot product operation, $\rotangle(\event_r)=\textrm{arccos}(\mucirc\cdot\rotvect(\event_r))$, whereas 
the unit random rotational axis $\axis(\event_r) / \|\axis(\event_r)\|$ is evaluated by the cross-product operation, $\axis(\event_r) / \|\axis(\event_r)\|=\mucirc\times\rotvect(\event_r)/\|\mucirc\times\rotvect \|$.
% By determining the random rotational angle $\rotangle(\event_r)$ and the random rotational axis $\axis(\event_r)$ using the cross-product operation, 
The corresponding random rotation matrix $\rotmat(\event_r)$ is then modelled as per the Rodrigues rotation formaula, described in \feq{Eq:Rodrigues}.
\ffig{Fig:orthoortho3D} summarizes the model $\conducttensorrot(\event_r)$. Here \ffig{Fig:orthoorthoVMunit3D} displays the mean and realisations of von Mises Fisher random variable on unit sphere; the equivalent ellipsoidal representation of stochastic model $\conducttensorrot(\event_r)$ is portrayed in \ffig{Fig:orthoorthorot3D}. It is evident that the shape and size of the ellipsoid in the mean and realisations remain constant; only the orientation of semi-axes of ellipsoids in the realisations fluctuate.

%% file: 3D_LN_iso_iso.tex
\definecolor{mycolor1}{rgb}{0.85000,0.33000,0.10000}%
\begin{tikzpicture}

\begin{axis}[%
width=\figW,
height=\figH,
at={(0\figW,0\figH)},
scale only axis,
xmin=0.2,
xmax=1,
xlabel style={font=\color{white!15!black}},
xlabel={$\eigvalsclele$ (W/mK)},
ymin=0,
ymax=10,
ylabel style={font=\color{white!15!black}},
ylabel={$f(\eigvalsclele)$},
axis background/.style={fill=white},
legend style={legend cell align=left, align=left, draw=white!15!black}
]
\addplot [color=mycolor1, line width=2.0pt]
  table[row sep=crcr]{%
0.317918242041234	7.83090365907542e-06\\
0.323459115978396	9.14664725867283e-05\\
0.328999989915558	0.000409772510664598\\
0.33454086385272	0.000708206276232978\\
0.340081737789882	0.000534535413881742\\
0.345622611727044	0.000575313016694337\\
0.351163485664206	0.00143974603678239\\
0.356704359601367	0.00269117078052519\\
0.362245233538529	0.00433202655396806\\
0.367786107475691	0.00758689700477421\\
0.373326981412853	0.0135201471305045\\
0.378867855350015	0.0243205381074356\\
0.384408729287177	0.04364798739184\\
0.389949603224339	0.0714221144231814\\
0.395490477161501	0.107524272228574\\
0.401031351098663	0.150624572948692\\
0.406572225035825	0.203021045085569\\
0.412113098972987	0.279220260303385\\
0.417653972910149	0.389742525510713\\
0.423194846847311	0.544069026994695\\
0.428735720784473	0.734152711906598\\
0.434276594721635	0.954098695424456\\
0.439817468658797	1.22211245143433\\
0.445358342595959	1.53888331100637\\
0.450899216533121	1.90580418336284\\
0.456440090470283	2.31319254093504\\
0.461980964407445	2.76544047961043\\
0.467521838344607	3.26726779260584\\
0.473062712281769	3.76568992325912\\
0.478603586218931	4.24449781991841\\
0.484144460156093	4.74073406633398\\
0.489685334093255	5.26479602899659\\
0.495226208030417	5.77122389748578\\
0.500767081967579	6.20734736464781\\
0.506307955904741	6.54614338556565\\
0.511848829841903	6.80072074482809\\
0.517389703779065	7.04678149067951\\
0.522930577716227	7.27528367966814\\
0.528471451653389	7.41847479056535\\
0.534012325590551	7.45793686314046\\
0.539553199527713	7.3913148552751\\
0.545094073464875	7.21575736219334\\
0.550634947402037	7.00356433620907\\
0.556175821339199	6.82661598136723\\
0.561716695276361	6.60182338373832\\
0.567257569213523	6.22319725499515\\
0.572798443150685	5.72522499162462\\
0.578339317087847	5.23810321874365\\
0.583880191025009	4.82278364187415\\
0.589421064962171	4.41967319320697\\
0.594961938899333	4.00537495553546\\
0.600502812836495	3.60581499934089\\
0.606043686773656	3.22397912657491\\
0.611584560710819	2.85415812443322\\
0.617125434647981	2.49589090497551\\
0.622666308585142	2.1725499670693\\
0.628207182522304	1.88384163337496\\
0.633748056459466	1.61350380149326\\
0.639288930396628	1.35901647137094\\
0.64482980433379	1.14334918829336\\
0.650370678270952	0.974625580685602\\
0.655911552208114	0.844037767478477\\
0.661452426145276	0.717879002580712\\
0.666993300082438	0.589851395995294\\
0.6725341740196	0.47934169139225\\
0.678075047956762	0.385587932938629\\
0.683615921893924	0.310758413001969\\
0.689156795831086	0.256156407288708\\
0.694697669768248	0.215228920968094\\
0.70023854370541	0.17652224713128\\
0.705779417642572	0.13748591913111\\
0.711320291579734	0.106397207496757\\
0.716861165516896	0.0858325193587758\\
0.722402039454058	0.0700293682873176\\
0.72794291339122	0.0528778057093145\\
0.733483787328382	0.038320620787391\\
0.739024661265544	0.0283020220634039\\
0.744565535202706	0.0207391639910899\\
0.750106409139868	0.0157860346261977\\
0.75564728307703	0.0116796867813348\\
0.761188157014192	0.0084645548205112\\
0.766729030951354	0.00614223239838463\\
0.772269904888516	0.00423124401836766\\
0.777810778825678	0.00367878857262198\\
0.78335165276284	0.00402912965968851\\
0.788892526700002	0.00366545062418795\\
0.794433400637164	0.00204985916763785\\
0.799974274574326	0.000550640154112587\\
0.805515148511488	6.20182455848949e-05\\
0.81105602244865	4.55737172445624e-06\\
0.816596896385812	3.29266018949391e-05\\
0.822137770322974	0.000230114052784536\\
0.827678644260136	0.000617056593525549\\
0.833219518197298	0.000644960657313844\\
0.83876039213446	0.000366655247868873\\
0.844301266071622	0.000501017191908049\\
0.849842140008784	0.000705709570264503\\
0.855383013945946	0.000409648651198648\\
0.860923887883108	9.14651646268963e-05\\
0.86646476182027	7.83089836279116e-06\\
};
\small \addlegendentry{$\eigvalscleleone,\eigvalscleletwo,\eigvalsclelethree$}

\end{axis}
\end{tikzpicture}%

%% file: 3D_LN_iso_ortho.tex
\definecolor{mycolor1}{rgb}{0.85000,0.33000,0.10000}%
\definecolor{mycolor2}{rgb}{0.72000,0.27000,1.00000}%
\definecolor{mycolor3}{rgb}{0.39000,0.83000,0.07000}%
\begin{tikzpicture}

\begin{axis}[%
width=\figW,
height=\figH,
at={(0\figW,0\figH)},
scale only axis,
xmin=0.328692615403957,
xmax=0.82334178060821,
xlabel style={font=\color{white!15!black}},
xlabel={$\eigvalsclele$ (W/mK)},
ymin=0,
ymax=8,
ylabel style={font=\color{white!15!black}},
ylabel={$f(\eigvalsclele)$},
axis background/.style={fill=white},
legend style={at={(1.03,1)}, anchor=north west, legend cell align=left, align=left, draw=white!15!black}
]
\addplot [color=mycolor1, line width=2.0pt]
  table[row sep=crcr]{%
0.328692615403957	7.83949680833059e-06\\
0.333689071618142	7.62895143237402e-05\\
0.338685527832326	0.000360072726883199\\
0.34368198404651	0.000964096484881909\\
0.348678440260695	0.00201646218690532\\
0.353674896474879	0.00369961484288961\\
0.358671352689063	0.0051734668462498\\
0.363667808903248	0.00611454295074284\\
0.368664265117432	0.00836782687472623\\
0.373660721331617	0.0147176661489005\\
0.378657177545801	0.0264980643680087\\
0.383653633759985	0.0413554300932365\\
0.38865008997417	0.0585569286082138\\
0.393646546188354	0.0828766492104202\\
0.398643002402538	0.118406354991588\\
0.403639458616723	0.166455155081374\\
0.408635914830907	0.235321508664423\\
0.413632371045092	0.333090454508992\\
0.418628827259276	0.457722009421113\\
0.42362528347346	0.599819898429183\\
0.428621739687645	0.751352377739134\\
0.433618195901829	0.923473392789292\\
0.438614652116013	1.14416549830251\\
0.443611108330198	1.42002444608891\\
0.448607564544382	1.73384173862733\\
0.453604020758567	2.07931578163138\\
0.458600476972751	2.46578684023272\\
0.463596933186935	2.8934447973109\\
0.46859338940112	3.32071408520995\\
0.473589845615304	3.73616976345403\\
0.478586301829488	4.18892272440373\\
0.483582758043673	4.67595227386792\\
0.488579214257857	5.1643034305314\\
0.493575670472042	5.61284537169262\\
0.498572126686226	5.99761148481012\\
0.50356858290041	6.33038519171247\\
0.508565039114595	6.6170973713731\\
0.513561495328779	6.8968860804757\\
0.518557951542963	7.17186394093725\\
0.523554407757148	7.39394065242762\\
0.528550863971332	7.52758267399566\\
0.533547320185517	7.56253711169674\\
0.538543776399701	7.50520216517596\\
0.543540232613885	7.35417940136771\\
0.54853668882807	7.10100104858339\\
0.553533145042254	6.80603621407592\\
0.558529601256438	6.54607223715993\\
0.563526057470623	6.30098739162379\\
0.568522513684807	6.00139124295453\\
0.573518969898992	5.64840046639359\\
0.578515426113176	5.28544809240874\\
0.58351188232736	4.9153701198387\\
0.588508338541545	4.53152245651234\\
0.593504794755729	4.13353440014273\\
0.598501250969913	3.74439094476256\\
0.603497707184098	3.39624624309344\\
0.608494163398282	3.06894179643761\\
0.613490619612467	2.74328640965401\\
0.618487075826651	2.43144663511484\\
0.623483532040835	2.1367897523227\\
0.62847998825502	1.86040614142116\\
0.633476444469204	1.61785089849885\\
0.638472900683388	1.41673014017243\\
0.643469356897573	1.24202191689672\\
0.648465813111757	1.07361000306636\\
0.653462269325942	0.911392132561069\\
0.658458725540126	0.771632367125628\\
0.66345518175431	0.66050854639095\\
0.668451637968495	0.564760760319964\\
0.673448094182679	0.469459371966871\\
0.678444550396863	0.381837695867582\\
0.683441006611048	0.313693375214894\\
0.688437462825232	0.259147285586185\\
0.693433919039417	0.211716769411788\\
0.698430375253601	0.17450878720293\\
0.703426831467785	0.146115811707075\\
0.70842328768197	0.122956869514763\\
0.713419743896154	0.105766949301517\\
0.718416200110338	0.0884209159634591\\
0.723412656324523	0.066060966850306\\
0.728409112538707	0.0463149141632995\\
0.733405568752892	0.0367498841765538\\
0.738402024967076	0.0337016531124491\\
0.74339848118126	0.0294672074198342\\
0.748394937395445	0.0242007487809013\\
0.753391393609629	0.0201401147055918\\
0.758387849823813	0.016276449692189\\
0.763384306037998	0.0128210845726041\\
0.768380762252182	0.00955127489528071\\
0.773377218466367	0.00597810211668584\\
0.778373674680551	0.00391209209408971\\
0.783370130894735	0.00344547841296212\\
0.78836658710892	0.00295499197868688\\
0.793363043323104	0.00198738842870724\\
0.798359499537288	0.00109844560421539\\
0.803355955751473	0.000816939228680363\\
0.808352411965657	0.000686247643071834\\
0.813348868179842	0.000328983331872484\\
0.818345324394026	7.46554833713754e-05\\
0.82334178060821	7.81987878148505e-06\\
};
\addlegendentry{$\eigvalscleleone$}

\addplot [color=mycolor2, dashed, line width=2.0pt]
  table[row sep=crcr]{%
0.328692615403957	7.83949680833059e-06\\
0.333689071618142	7.62895143237402e-05\\
0.338685527832326	0.000360072726883199\\
0.34368198404651	0.000964096484881909\\
0.348678440260695	0.00201646218690532\\
0.353674896474879	0.00369961484288961\\
0.358671352689063	0.0051734668462498\\
0.363667808903248	0.00611454295074284\\
0.368664265117432	0.00836782687472623\\
0.373660721331617	0.0147176661489005\\
0.378657177545801	0.0264980643680087\\
0.383653633759985	0.0413554300932365\\
0.38865008997417	0.0585569286082138\\
0.393646546188354	0.0828766492104202\\
0.398643002402538	0.118406354991588\\
0.403639458616723	0.166455155081374\\
0.408635914830907	0.235321508664423\\
0.413632371045092	0.333090454508992\\
0.418628827259276	0.457722009421113\\
0.42362528347346	0.599819898429183\\
0.428621739687645	0.751352377739134\\
0.433618195901829	0.923473392789292\\
0.438614652116013	1.14416549830251\\
0.443611108330198	1.42002444608891\\
0.448607564544382	1.73384173862733\\
0.453604020758567	2.07931578163138\\
0.458600476972751	2.46578684023272\\
0.463596933186935	2.8934447973109\\
0.46859338940112	3.32071408520995\\
0.473589845615304	3.73616976345403\\
0.478586301829488	4.18892272440373\\
0.483582758043673	4.67595227386792\\
0.488579214257857	5.1643034305314\\
0.493575670472042	5.61284537169262\\
0.498572126686226	5.99761148481012\\
0.50356858290041	6.33038519171247\\
0.508565039114595	6.6170973713731\\
0.513561495328779	6.8968860804757\\
0.518557951542963	7.17186394093725\\
0.523554407757148	7.39394065242762\\
0.528550863971332	7.52758267399566\\
0.533547320185517	7.56253711169674\\
0.538543776399701	7.50520216517596\\
0.543540232613885	7.35417940136771\\
0.54853668882807	7.10100104858339\\
0.553533145042254	6.80603621407592\\
0.558529601256438	6.54607223715993\\
0.563526057470623	6.30098739162379\\
0.568522513684807	6.00139124295453\\
0.573518969898992	5.64840046639359\\
0.578515426113176	5.28544809240874\\
0.58351188232736	4.9153701198387\\
0.588508338541545	4.53152245651234\\
0.593504794755729	4.13353440014273\\
0.598501250969913	3.74439094476256\\
0.603497707184098	3.39624624309344\\
0.608494163398282	3.06894179643761\\
0.613490619612467	2.74328640965401\\
0.618487075826651	2.43144663511484\\
0.623483532040835	2.1367897523227\\
0.62847998825502	1.86040614142116\\
0.633476444469204	1.61785089849885\\
0.638472900683388	1.41673014017243\\
0.643469356897573	1.24202191689672\\
0.648465813111757	1.07361000306636\\
0.653462269325942	0.911392132561069\\
0.658458725540126	0.771632367125628\\
0.66345518175431	0.66050854639095\\
0.668451637968495	0.564760760319964\\
0.673448094182679	0.469459371966871\\
0.678444550396863	0.381837695867582\\
0.683441006611048	0.313693375214894\\
0.688437462825232	0.259147285586185\\
0.693433919039417	0.211716769411788\\
0.698430375253601	0.17450878720293\\
0.703426831467785	0.146115811707075\\
0.70842328768197	0.122956869514763\\
0.713419743896154	0.105766949301517\\
0.718416200110338	0.0884209159634591\\
0.723412656324523	0.066060966850306\\
0.728409112538707	0.0463149141632995\\
0.733405568752892	0.0367498841765538\\
0.738402024967076	0.0337016531124491\\
0.74339848118126	0.0294672074198342\\
0.748394937395445	0.0242007487809013\\
0.753391393609629	0.0201401147055918\\
0.758387849823813	0.016276449692189\\
0.763384306037998	0.0128210845726041\\
0.768380762252182	0.00955127489528071\\
0.773377218466367	0.00597810211668584\\
0.778373674680551	0.00391209209408971\\
0.783370130894735	0.00344547841296212\\
0.78836658710892	0.00295499197868688\\
0.793363043323104	0.00198738842870724\\
0.798359499537288	0.00109844560421539\\
0.803355955751473	0.000816939228680363\\
0.808352411965657	0.000686247643071834\\
0.813348868179842	0.000328983331872484\\
0.818345324394026	7.46554833713754e-05\\
0.82334178060821	7.81987878148505e-06\\
};
\addlegendentry{$\eigvalscleletwo$}

\addplot [color=mycolor3, dashed, line width=2.0pt]
  table[row sep=crcr]{%
0.328692615403957	7.83949680833059e-06\\
0.333689071618142	7.62895143237402e-05\\
0.338685527832326	0.000360072726883199\\
0.34368198404651	0.000964096484881909\\
0.348678440260695	0.00201646218690532\\
0.353674896474879	0.00369961484288961\\
0.358671352689063	0.0051734668462498\\
0.363667808903248	0.00611454295074284\\
0.368664265117432	0.00836782687472623\\
0.373660721331617	0.0147176661489005\\
0.378657177545801	0.0264980643680087\\
0.383653633759985	0.0413554300932365\\
0.38865008997417	0.0585569286082138\\
0.393646546188354	0.0828766492104202\\
0.398643002402538	0.118406354991588\\
0.403639458616723	0.166455155081374\\
0.408635914830907	0.235321508664423\\
0.413632371045092	0.333090454508992\\
0.418628827259276	0.457722009421113\\
0.42362528347346	0.599819898429183\\
0.428621739687645	0.751352377739134\\
0.433618195901829	0.923473392789292\\
0.438614652116013	1.14416549830251\\
0.443611108330198	1.42002444608891\\
0.448607564544382	1.73384173862733\\
0.453604020758567	2.07931578163138\\
0.458600476972751	2.46578684023272\\
0.463596933186935	2.8934447973109\\
0.46859338940112	3.32071408520995\\
0.473589845615304	3.73616976345403\\
0.478586301829488	4.18892272440373\\
0.483582758043673	4.67595227386792\\
0.488579214257857	5.1643034305314\\
0.493575670472042	5.61284537169262\\
0.498572126686226	5.99761148481012\\
0.50356858290041	6.33038519171247\\
0.508565039114595	6.6170973713731\\
0.513561495328779	6.8968860804757\\
0.518557951542963	7.17186394093725\\
0.523554407757148	7.39394065242762\\
0.528550863971332	7.52758267399566\\
0.533547320185517	7.56253711169674\\
0.538543776399701	7.50520216517596\\
0.543540232613885	7.35417940136771\\
0.54853668882807	7.10100104858339\\
0.553533145042254	6.80603621407592\\
0.558529601256438	6.54607223715993\\
0.563526057470623	6.30098739162379\\
0.568522513684807	6.00139124295453\\
0.573518969898992	5.64840046639359\\
0.578515426113176	5.28544809240874\\
0.58351188232736	4.9153701198387\\
0.588508338541545	4.53152245651234\\
0.593504794755729	4.13353440014273\\
0.598501250969913	3.74439094476256\\
0.603497707184098	3.39624624309344\\
0.608494163398282	3.06894179643761\\
0.613490619612467	2.74328640965401\\
0.618487075826651	2.43144663511484\\
0.623483532040835	2.1367897523227\\
0.62847998825502	1.86040614142116\\
0.633476444469204	1.61785089849885\\
0.638472900683388	1.41673014017243\\
0.643469356897573	1.24202191689672\\
0.648465813111757	1.07361000306636\\
0.653462269325942	0.911392132561069\\
0.658458725540126	0.771632367125628\\
0.66345518175431	0.66050854639095\\
0.668451637968495	0.564760760319964\\
0.673448094182679	0.469459371966871\\
0.678444550396863	0.381837695867582\\
0.683441006611048	0.313693375214894\\
0.688437462825232	0.259147285586185\\
0.693433919039417	0.211716769411788\\
0.698430375253601	0.17450878720293\\
0.703426831467785	0.146115811707075\\
0.70842328768197	0.122956869514763\\
0.713419743896154	0.105766949301517\\
0.718416200110338	0.0884209159634591\\
0.723412656324523	0.066060966850306\\
0.728409112538707	0.0463149141632995\\
0.733405568752892	0.0367498841765538\\
0.738402024967076	0.0337016531124491\\
0.74339848118126	0.0294672074198342\\
0.748394937395445	0.0242007487809013\\
0.753391393609629	0.0201401147055918\\
0.758387849823813	0.016276449692189\\
0.763384306037998	0.0128210845726041\\
0.768380762252182	0.00955127489528071\\
0.773377218466367	0.00597810211668584\\
0.778373674680551	0.00391209209408971\\
0.783370130894735	0.00344547841296212\\
0.78836658710892	0.00295499197868688\\
0.793363043323104	0.00198738842870724\\
0.798359499537288	0.00109844560421539\\
0.803355955751473	0.000816939228680363\\
0.808352411965657	0.000686247643071834\\
0.813348868179842	0.000328983331872484\\
0.818345324394026	7.46554833713754e-05\\
0.82334178060821	7.81987878148505e-06\\
};
\addlegendentry{$\eigvalsclelethree$}

\end{axis}
\end{tikzpicture}%

%% file: results_3D.tex
% !TEX root = ../paper_tensors.tex
% !TEX encoding = UTF-8 Unicode
\newcommand{\widththreeD}{0.85} %width size for 3D images in results section
\newcommand{\widththreeDdir}{0.95} %width size for 3D images in results section
%------------------------------------------------------------------mean NT---------------------------------------------------------

\begin{figure}[!h]
	\begin{subfigure}{0.32\textwidth}
		\centering
		\includegraphics[width=\widththreeD\linewidth]{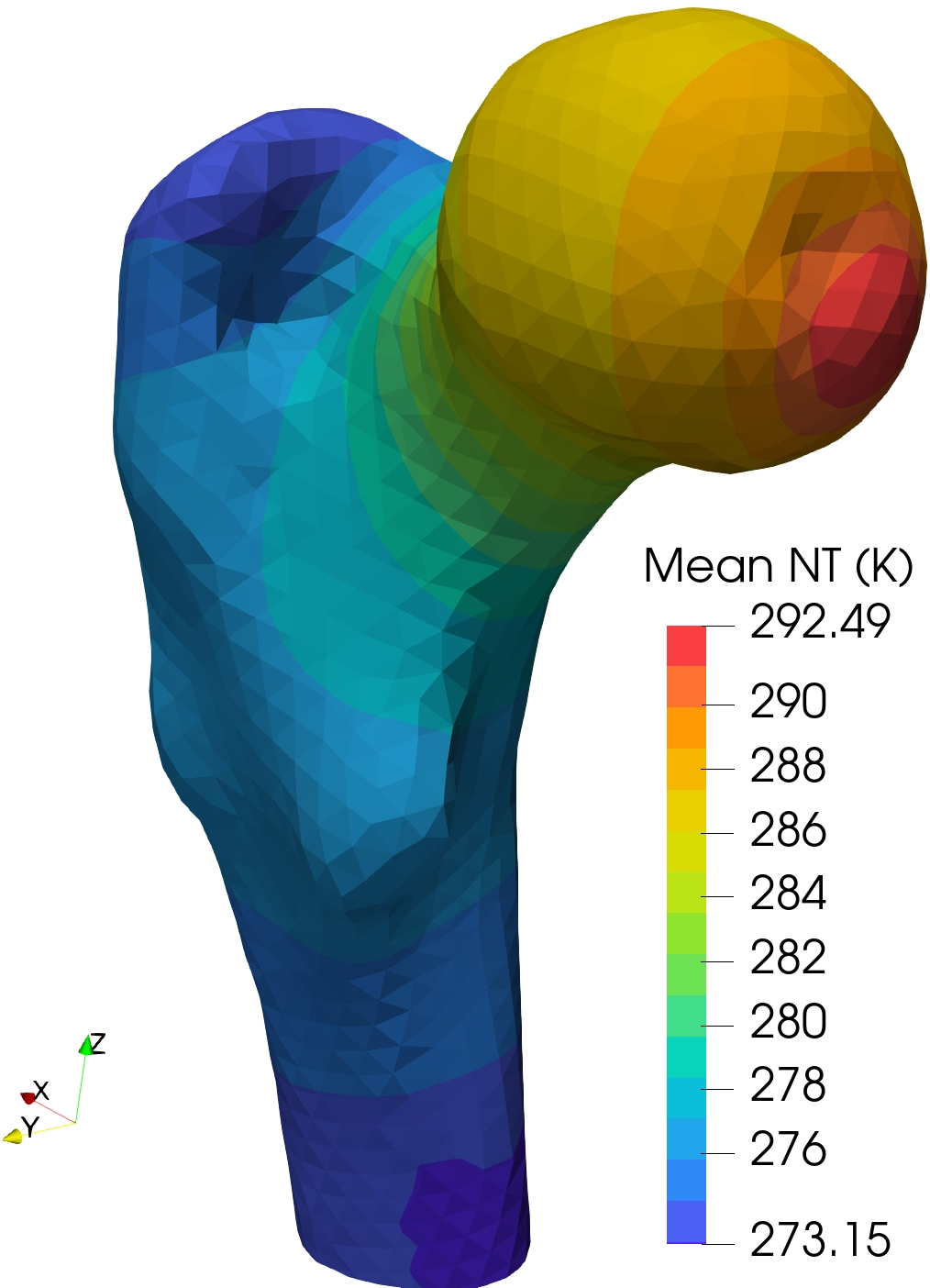}
		\caption{iso-iso-scl}
	\end{subfigure}
	\begin{subfigure}{0.32\textwidth}
		\centering
		\includegraphics[width=\widththreeD\linewidth]{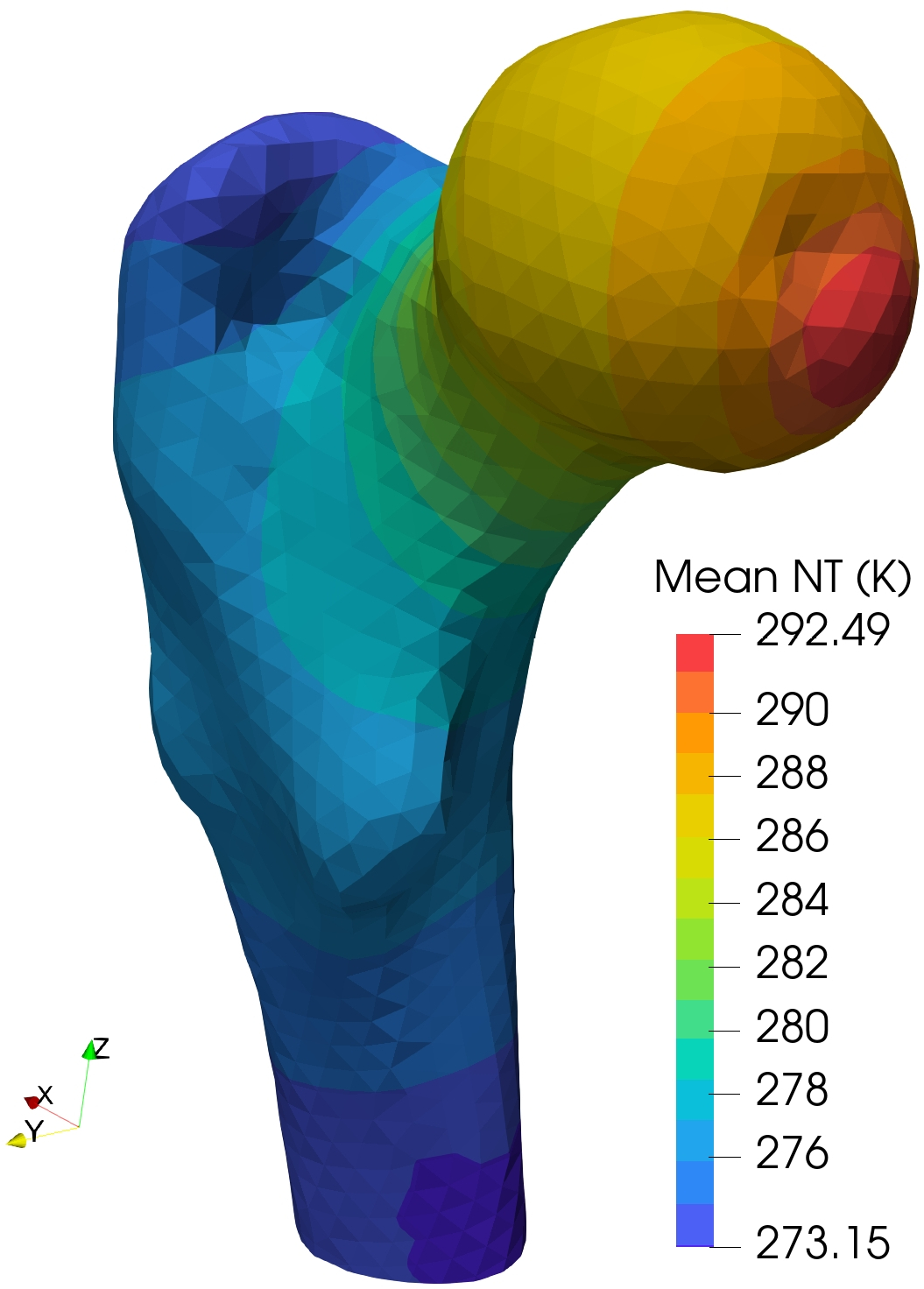}
		\caption{iso-ortho-scl}
	\end{subfigure}
% 	\begin{subfigure}{0.32\textwidth}
% 		\centering
% 		\includegraphics[width=\widththreeD\linewidth]{3D_scl_ortho_ortho_CCX_sideflux_1e5_mean_NT}
% 		\caption{ortho-ortho-scl}
% 	\end{subfigure}
	\begin{subfigure}{0.32\textwidth}
		\centering
		\includegraphics[width=\widththreeD\linewidth]{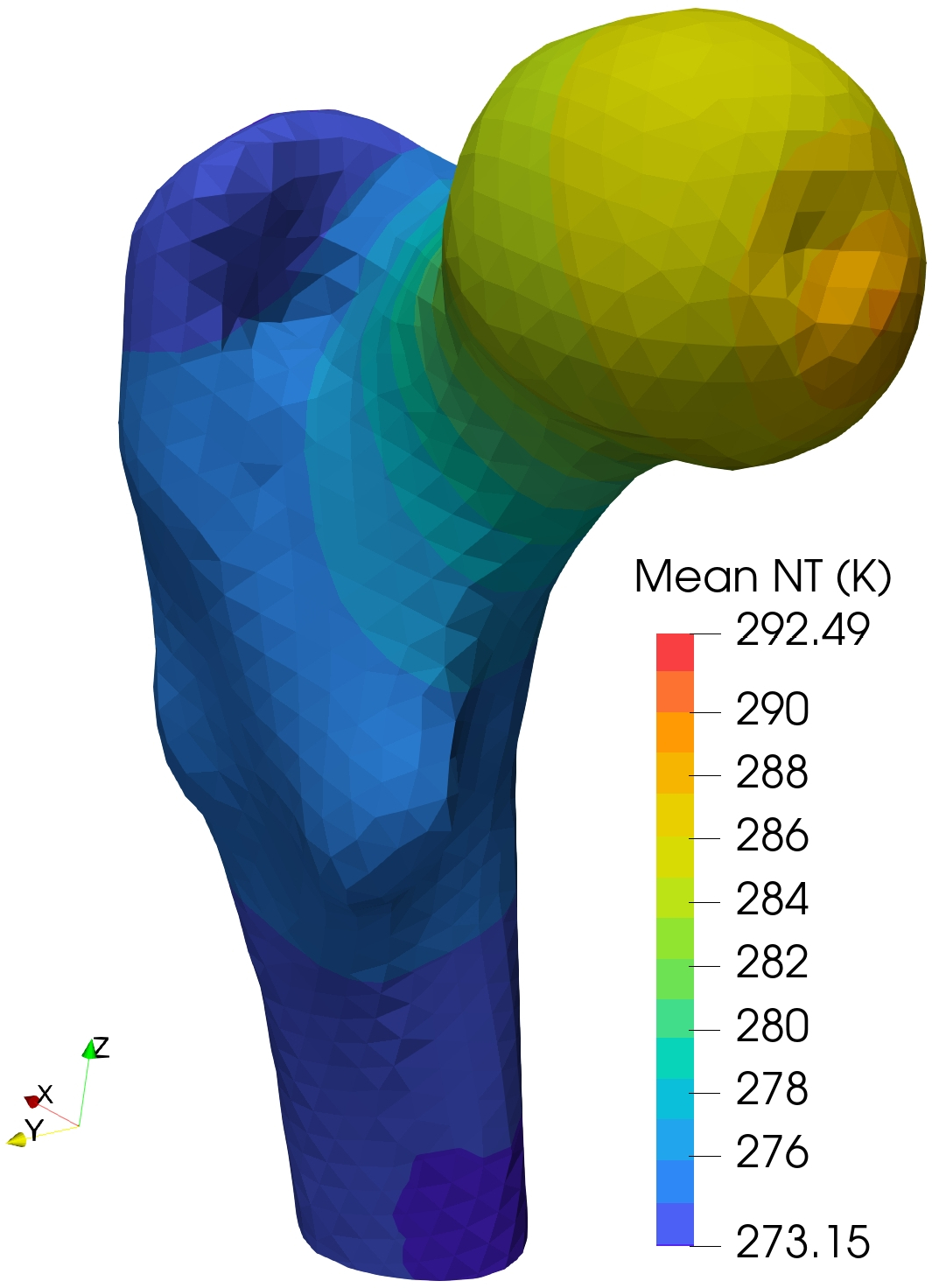}
		\caption{\acrshort{ortho}-\acrshort{ortho}-\acrshort{dir}}
	\end{subfigure}
% 	\begin{subfigure}{0.32\textwidth}
% 		\centering
% 		\includegraphics[width=\widththreeD\linewidth]{3D_rotscl_ortho_ortho_CCX_sideflux_1e5_mean_NT}
% 		\caption{ortho-ortho-rotscl}
% 	\end{subfigure}
% 	\begin{subfigure}{0.32\textwidth}
% 		\centering
% 		\includegraphics[width=\widththreeD\linewidth]{3D_rotscl_iso_ortho_CCX_sideflux_1e5_mean_NT}
% 		\caption{iso-ortho-rotscl}
% 	\end{subfigure}
	\caption{\acrshort{mc} mean estimate of nodal temperature}
	\label{Fig:meanNT3D}
\end{figure}

\begin{figure}[!h]
	\begin{subfigure}{0.32\textwidth}
		\centering
		\includegraphics[width=\widththreeD\linewidth]{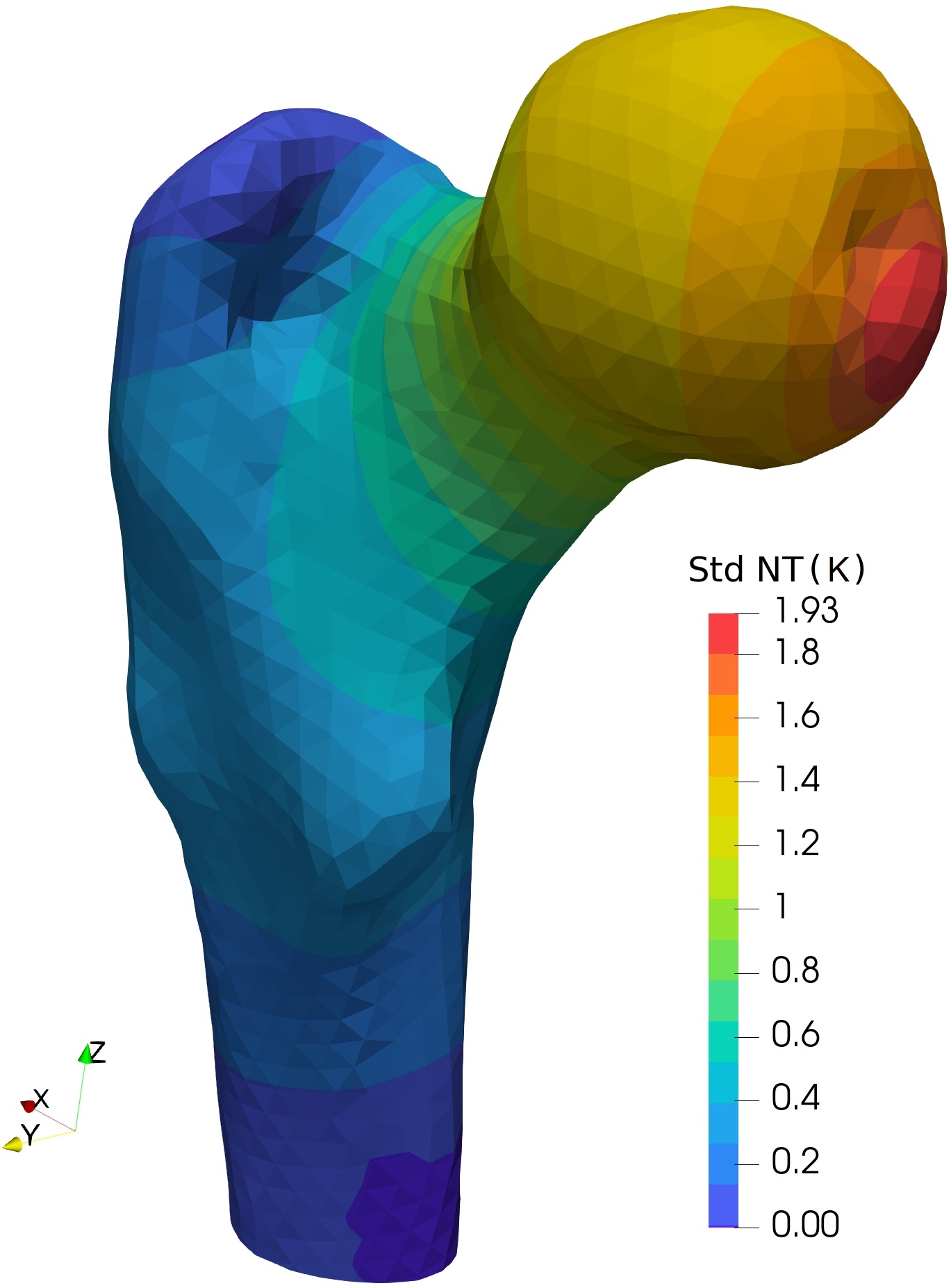}
		\caption{iso-iso-scl}
	\end{subfigure}
	\begin{subfigure}{0.32\textwidth}
		\centering
		\includegraphics[width=\widththreeD\linewidth]{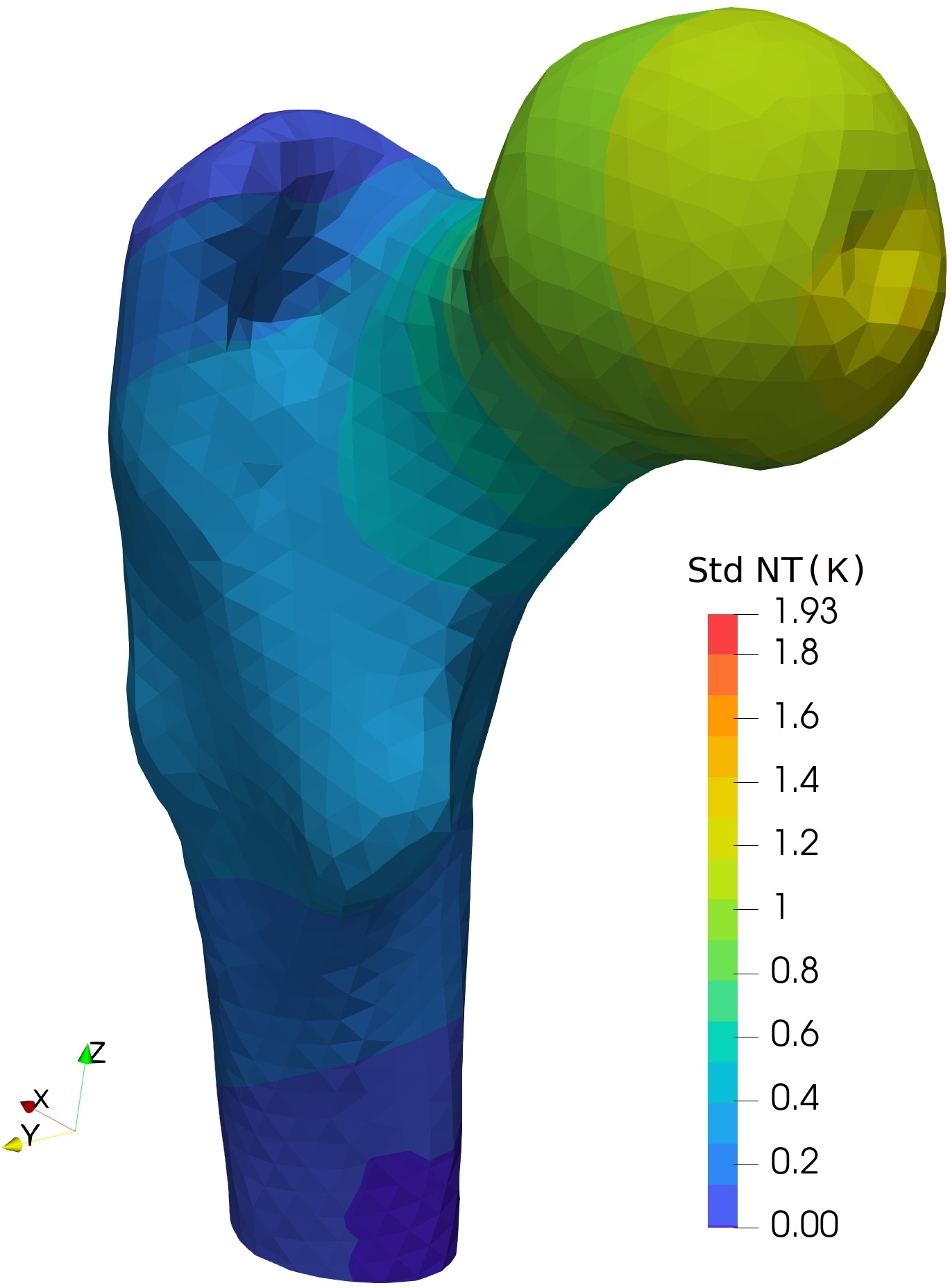}
		\caption{iso-ortho-scl}
	\end{subfigure}
% 	\begin{subfigure}{0.32\textwidth}
% 		\centering
% 		\includegraphics[width=\widththreeD\linewidth]{3D_scl_ortho_ortho_CCX_sideflux_1e5_var_NT}
% 		\caption{ortho-ortho-scl}
% 	\end{subfigure}
	\begin{subfigure}{0.32\textwidth}
		\centering
		\includegraphics[width=\widththreeD\linewidth]{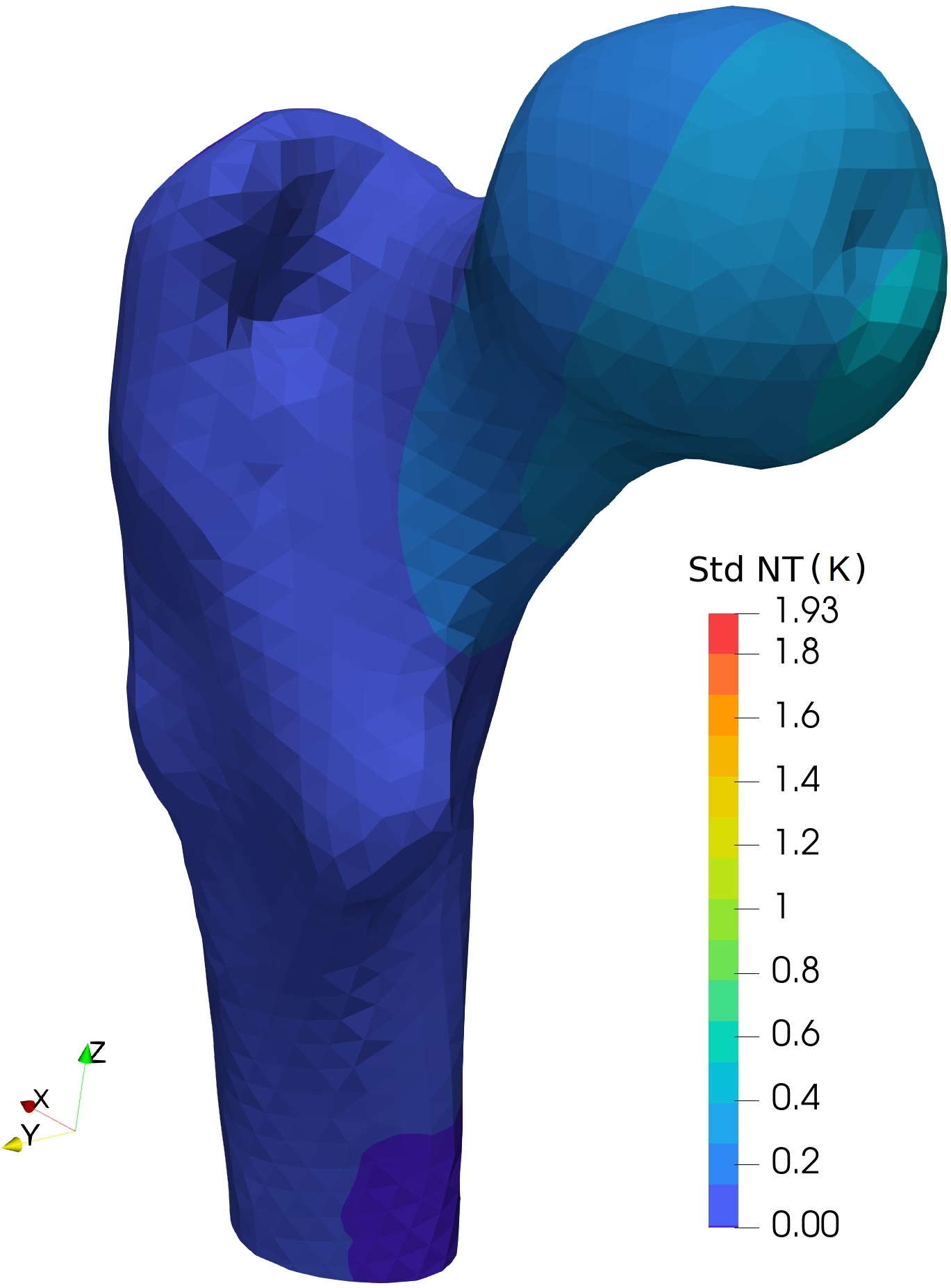}
		\caption{\acrshort{ortho}-\acrshort{ortho}-\acrshort{dir}}
	\end{subfigure}
% 	\begin{subfigure}{0.32\textwidth}
% 		\centering
% 		\includegraphics[width=\widththreeD\linewidth]{3D_rotscl_ortho_ortho_CCX_sideflux_1e5_var_NT}
% 		\caption{ortho-ortho-rotscl}
% 	\end{subfigure}
% 	\begin{subfigure}{0.32\textwidth}
% 		\centering
% 		\includegraphics[width=\widththreeD\linewidth]{3D_rotscl_iso_ortho_CCX_sideflux_1e5_var_NT}
% 		\caption{iso-ortho-rotscl}
% 	\end{subfigure}
	\caption{\acrshort{mc} standard deviation estimate of nodal temperature}
	\label{Fig:stdNT3D}
\end{figure}

\begin{figure}[!h]
	\begin{subfigure}{0.32\textwidth}
		\centering
		\includegraphics[width=\widththreeDdir\linewidth]{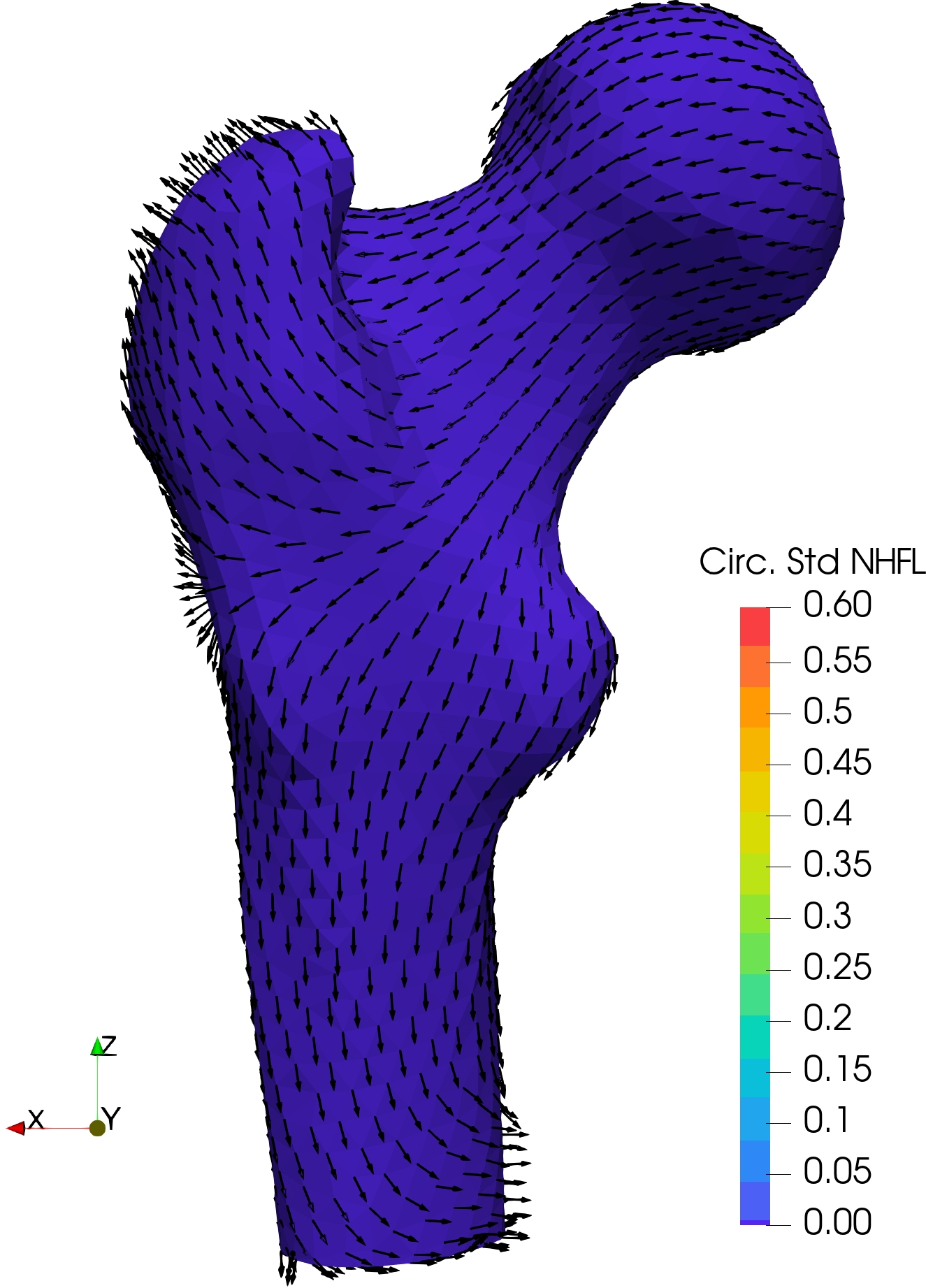}
		\caption{iso-iso-scl}
	\end{subfigure}
	\begin{subfigure}{0.32\textwidth}
		\centering
		\includegraphics[width=\widththreeDdir\linewidth]{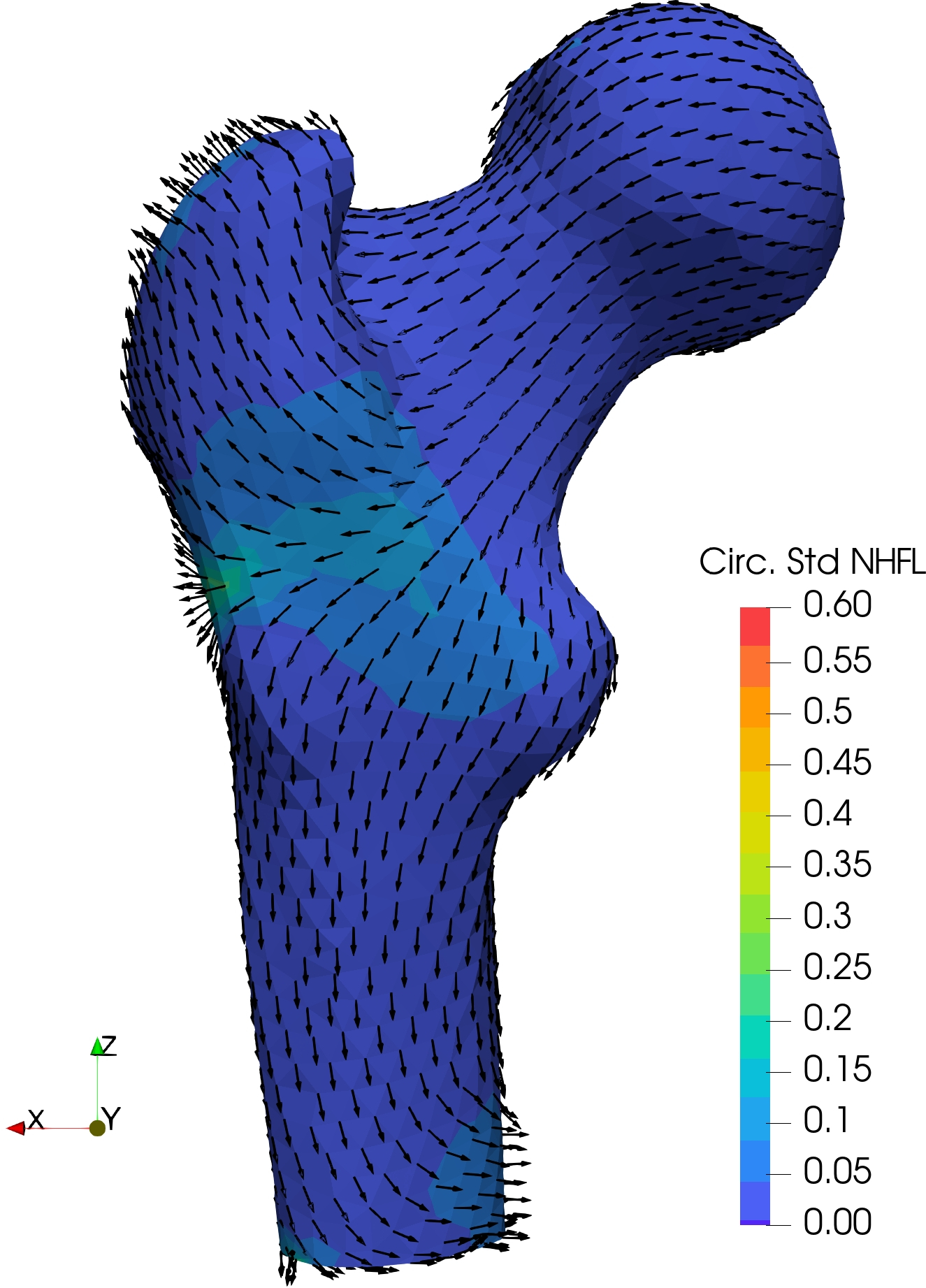}
		\caption{iso-ortho-scl}
	\end{subfigure}
% 	\begin{subfigure}{0.32\textwidth}
% 		\centering
% 		\includegraphics[width=\widththreeDdir\linewidth]{3D_scl_ortho_ortho_CCX_sideflux_1e5_mean_var_dir}
% 		\caption{ortho-ortho-scl}
% 	\end{subfigure}
	\begin{subfigure}{0.32\textwidth}
		\centering
		\includegraphics[width=\widththreeDdir\linewidth]{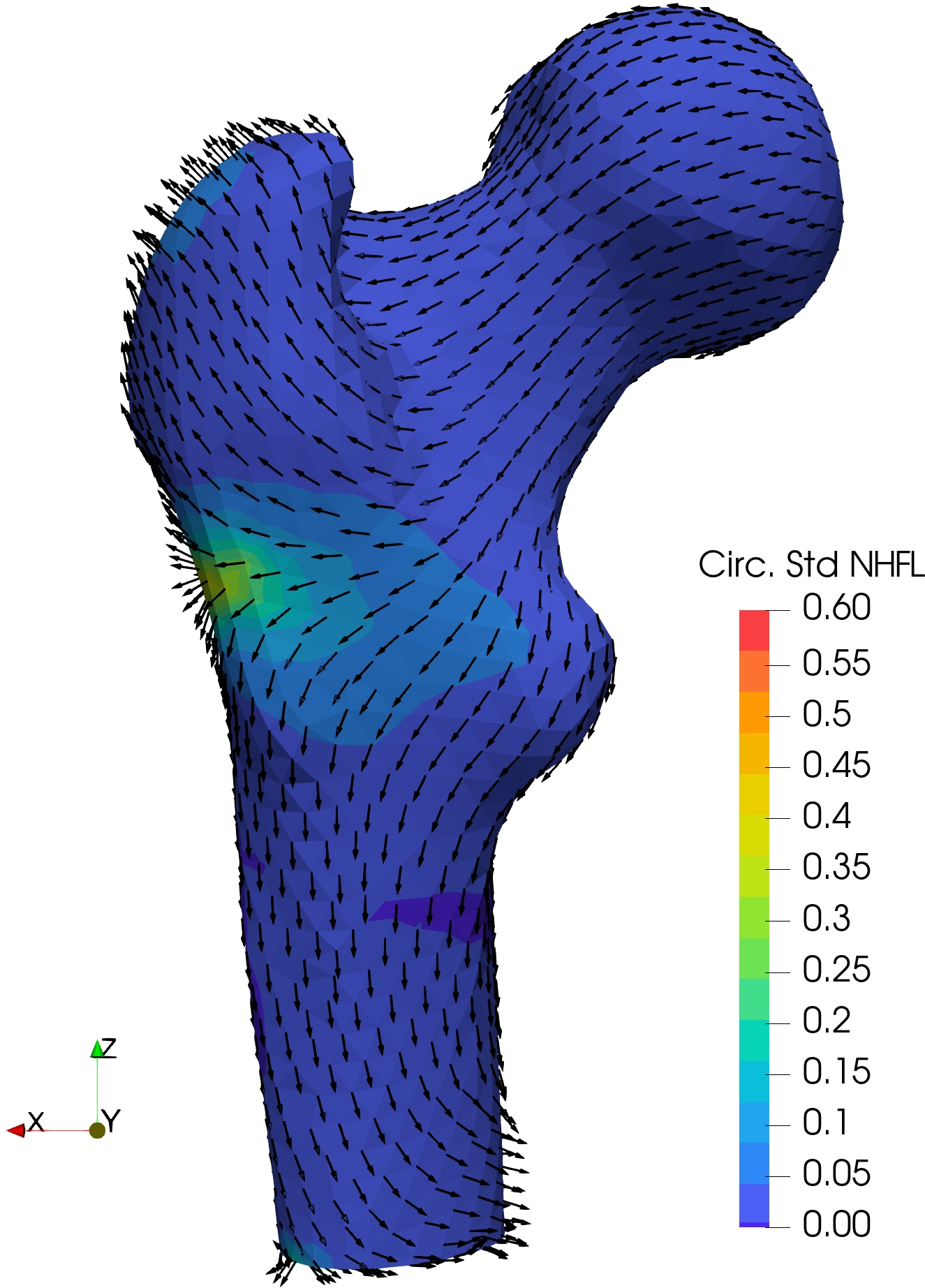}
		\caption{\acrshort{ortho}-\acrshort{ortho}-\acrshort{dir}}
	\end{subfigure}
% 	\begin{subfigure}{0.32\textwidth}
% 		\centering
% 		\includegraphics[width=\widththreeDdir\linewidth]{3D_rotscl_ortho_ortho_CCX_sideflux_1e5_mean_var_dir_1}
% 		\caption{ortho-ortho-rotscl}
% 	\end{subfigure}
% 	\begin{subfigure}{0.32\textwidth}
% 		\centering
% 		\includegraphics[width=\widththreeDdir\linewidth]{3D_rotscl_iso_ortho_CCX_sideflux_1e5_mean_var_dir_1}
% 		\caption{iso-ortho-rotscl}
% 	\end{subfigure}
	\caption{\acrshort{mc} circular mean (vector representation) and circular standard deviation estimates (dimensionless quantity) of normalized heat flux}
	\label{Fig:meanstdNHFL3D}
\end{figure}

\paragraph{Uncertainty quantification}:
\ffiigs{Fig:meanNT3D}{Fig:stdTHFL3D} display the MC mean and standard deviation estimates of NT and THFL, whereas \ffig{Fig:meanstdNHFL3D} shows the circular MC mean and circular standard deviation estimate of NHFL. 
Concerning the sensitivity of different stochastic models on the desired quantities, we make similar observations, and hence, draw similar conclusions on the 3D example as in the previous 2D case (from \fsec{Sec:2DResults}). That is, the findings in  \ffiigs{Fig:meanNT}{Fig:circHFL} correspond approximately to the results in \ffiigs{Fig:meanNT3D}{Fig:meanstdNHFL3D}.
% However, the magnitude of MC estimates of quantities of interest of a particular type of stochastic model does not correspond to the results in the 2D scenario.
% We make similar observations and hence draw similar conclusions as explained in the 2D case. 

The numerical results of both 2D and 3D proximal femur signify the prominence of incorporating different material uncertainties---scaling, orientation and material symmetry---independently into the constitutive model;  the results showcase the distinct impact of stochastic models on the desired quantities of interest, such as nodal temperature and heat flux. 
% the stochastic impact on the solutions is different for different stochastic model..

%----------------------------------------------------------var NT--------------------------------------------------------------------

%-------------------------------------------------------------mean HFL--------------------------------------------------------------

\begin{figure}[!h]
	\begin{subfigure}{0.32\textwidth}
		\centering
		\includegraphics[width=\widththreeD\linewidth]{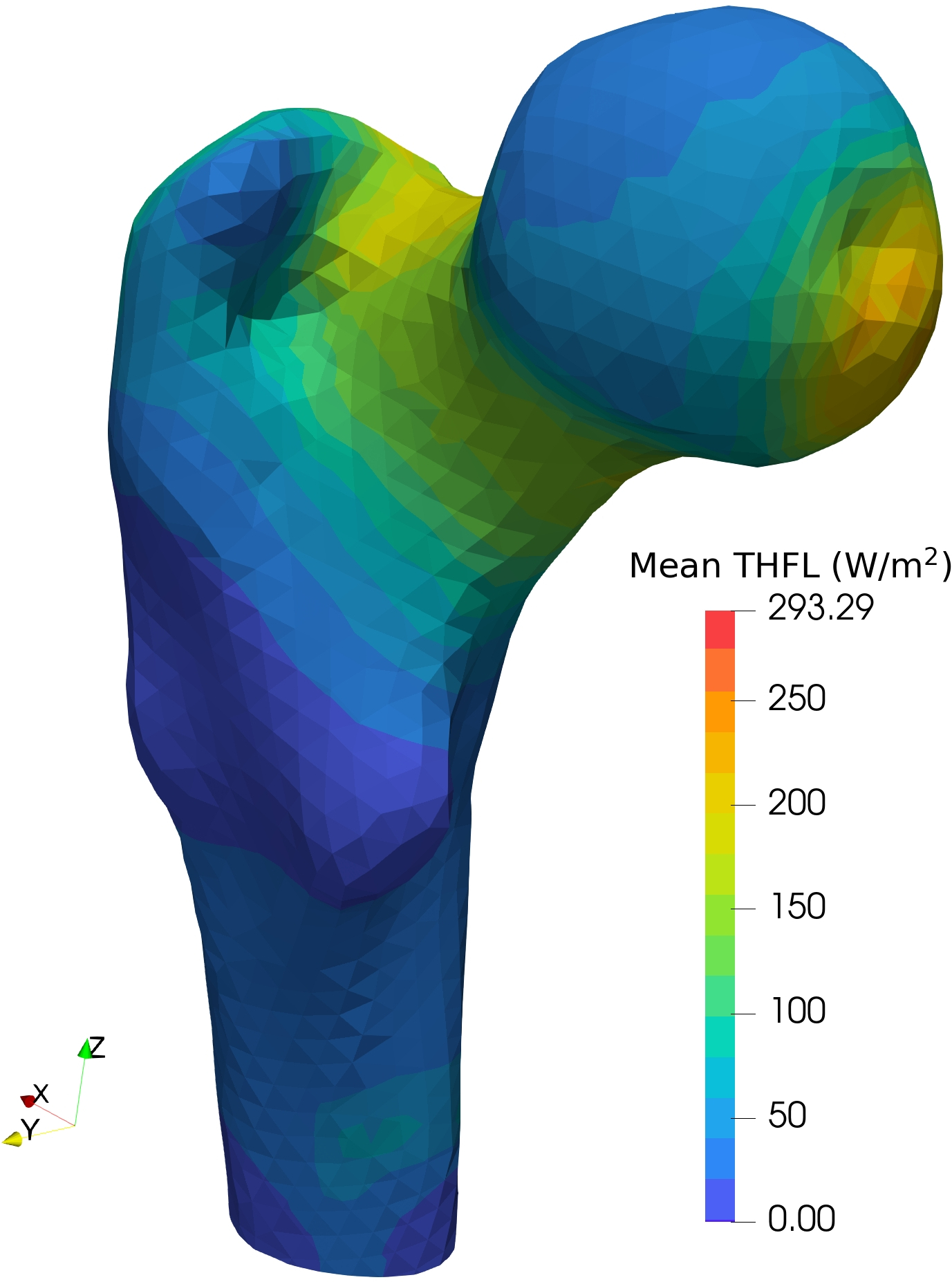}
		\caption{iso-iso-scl}
	\end{subfigure}
	\begin{subfigure}{0.32\textwidth}
		\centering
		\includegraphics[width=\widththreeD\linewidth]{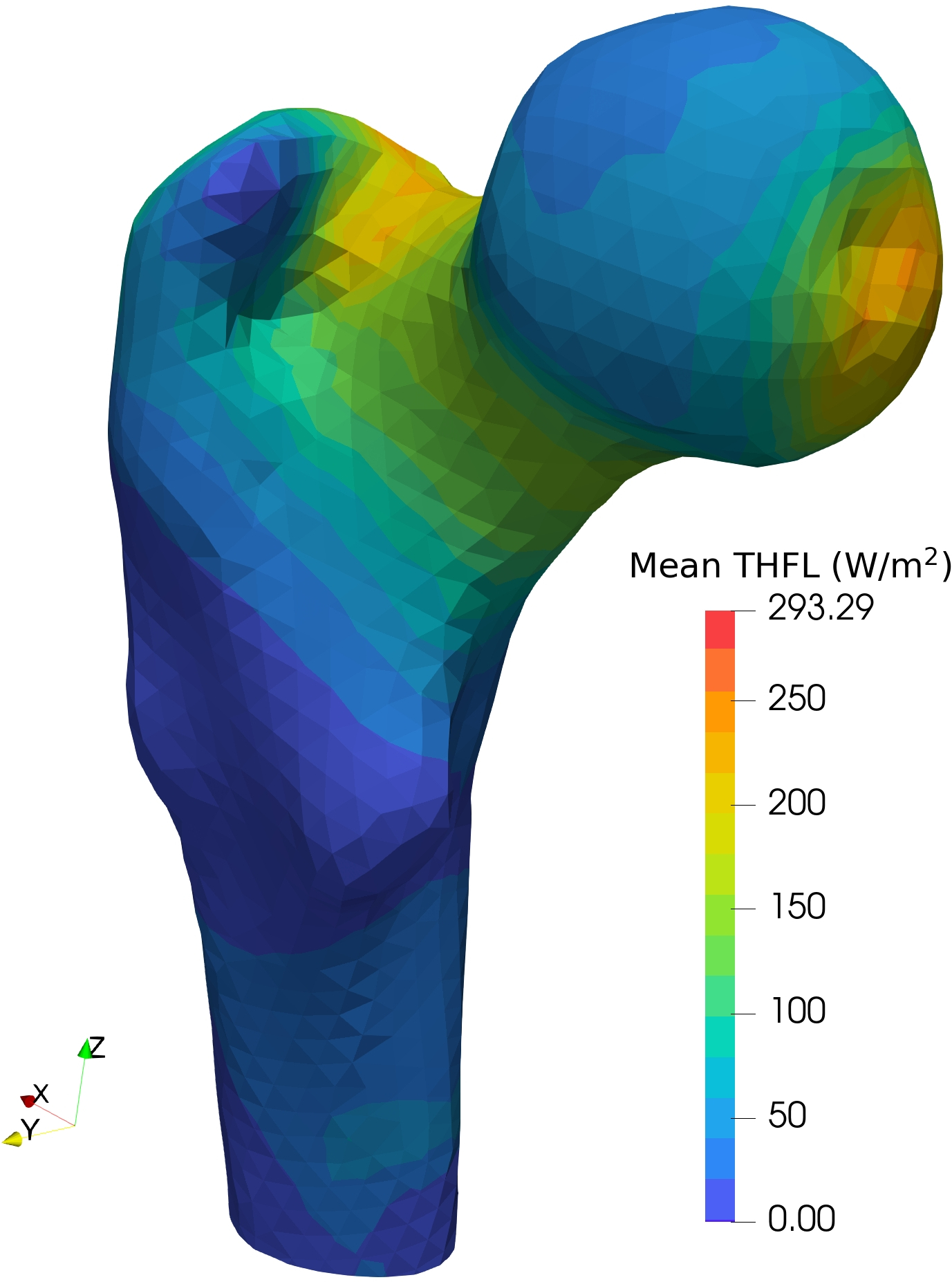}
		\caption{iso-ortho-scl}
	\end{subfigure}
% 	\begin{subfigure}{0.32\textwidth}
% 		\centering
% 		\includegraphics[width=\widththreeD\linewidth]{3D_scl_ortho_ortho_CCX_sideflux_1e5_mean_HFL}
% 		\caption{ortho-ortho-scl}
% 	\end{subfigure}
	\begin{subfigure}{0.32\textwidth}
		\centering
		\includegraphics[width=\widththreeD\linewidth]{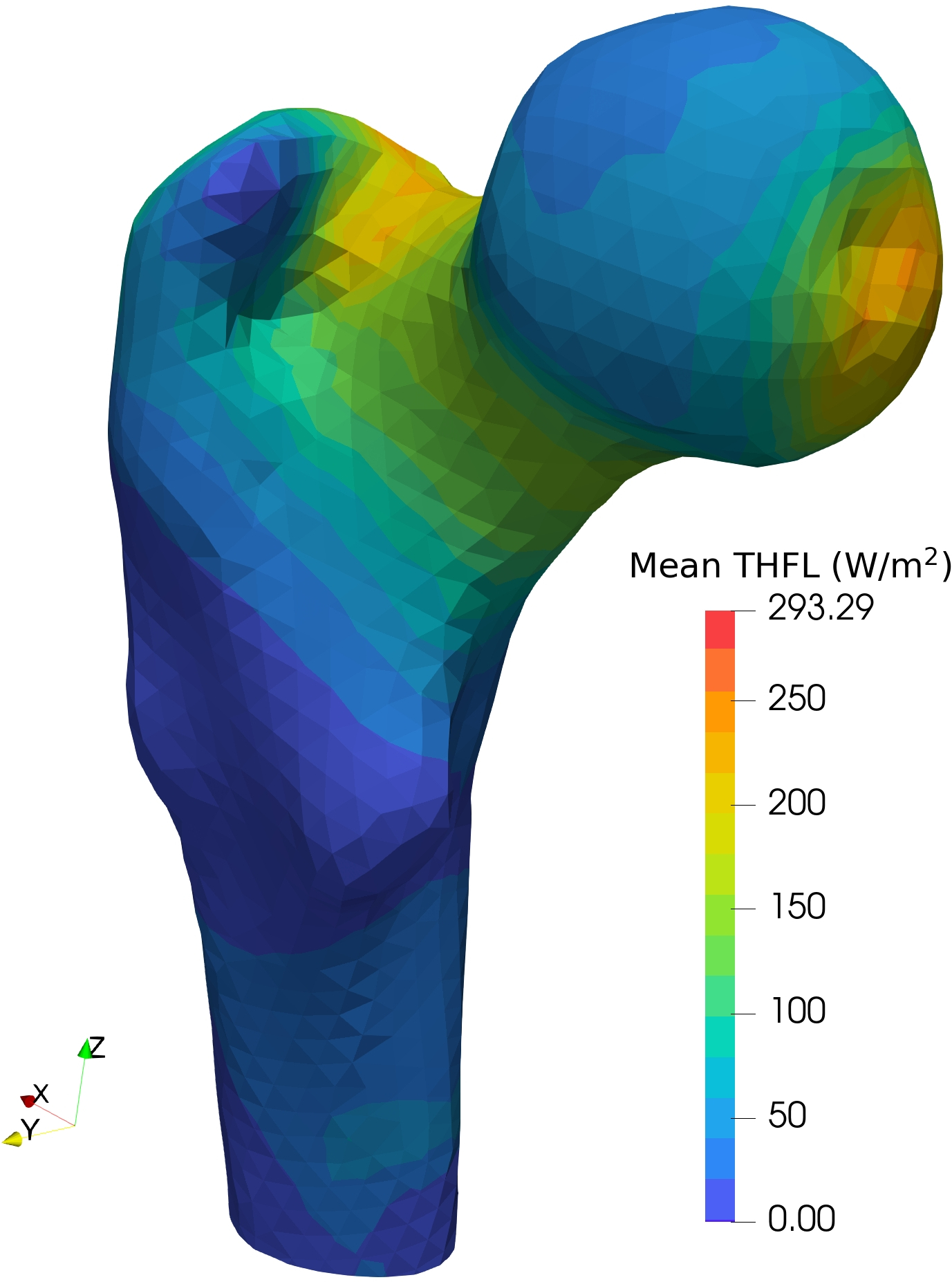}
		\caption{\acrshort{ortho}-\acrshort{ortho}-\acrshort{dir}}
	\end{subfigure}
% 	\begin{subfigure}{0.32\textwidth}
% 		\centering
% 		\includegraphics[width=\widththreeD\linewidth]{3D_rotscl_ortho_ortho_CCX_sideflux_1e5_mean_HFL}
% 		\caption{ortho-ortho-rotscl}
% 	\end{subfigure}
% 	\begin{subfigure}{0.32\textwidth}
% 		\centering
% 		\includegraphics[width=\widththreeD\linewidth]{3D_rotscl_iso_ortho_CCX_sideflux_1e5_mean_HFL}
% 		\caption{iso-ortho-rotscl}
% 	\end{subfigure}
	\caption{\acrshort{mc} mean estimate of total heat flux (Euclidean norm)}
	\label{Fig:meanTHFL3D}
\end{figure}

%-------------------------------------------------------------var HFL----------------------------------------------------------------

\begin{figure}[!h]
	\begin{subfigure}{0.32\textwidth}
		\centering
		\includegraphics[width=\widththreeD\linewidth]{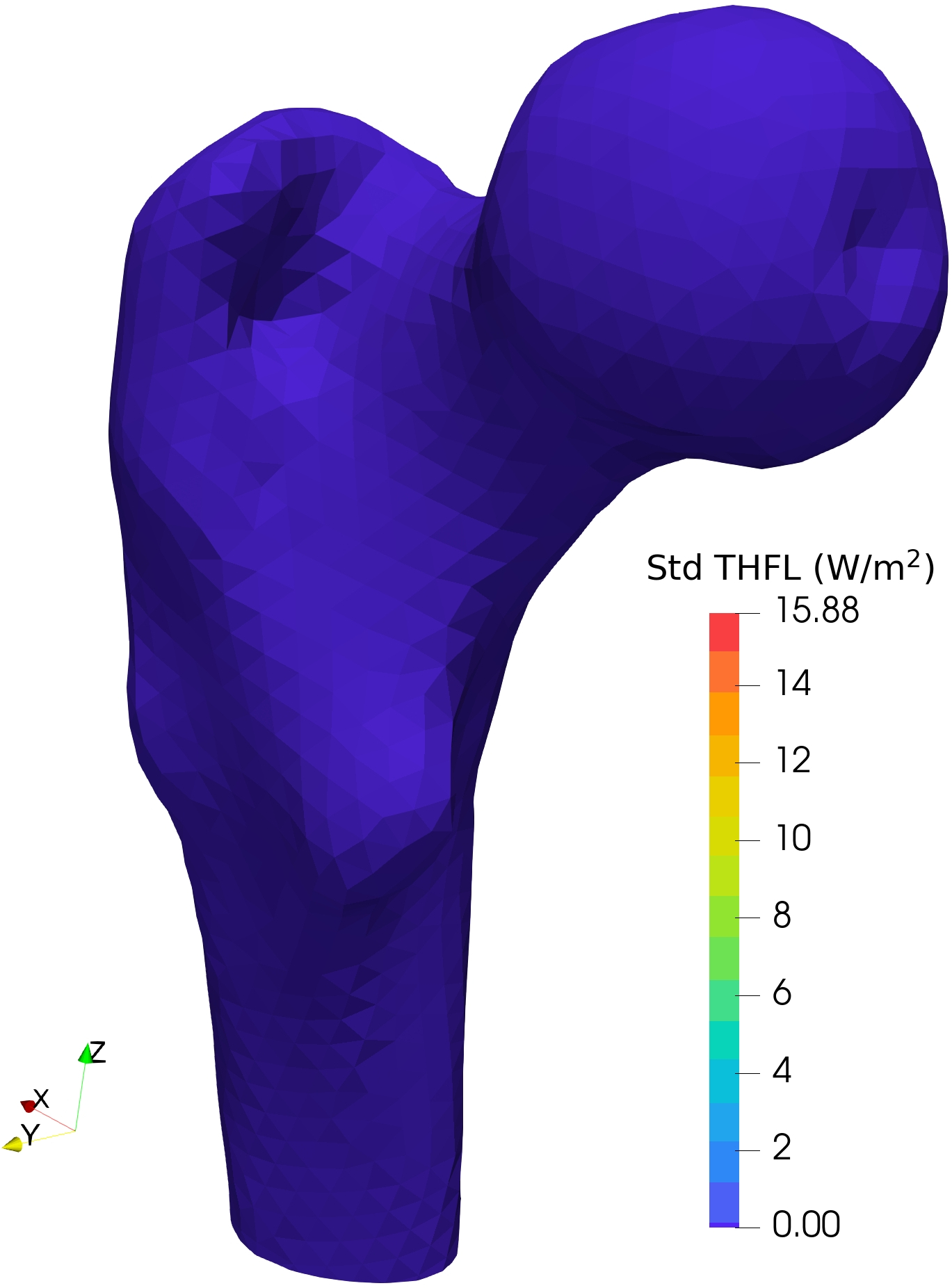}
		\caption{iso-iso-scl}
	\end{subfigure}
	\begin{subfigure}{0.32\textwidth}
		\centering
		\includegraphics[width=\widththreeD\linewidth]{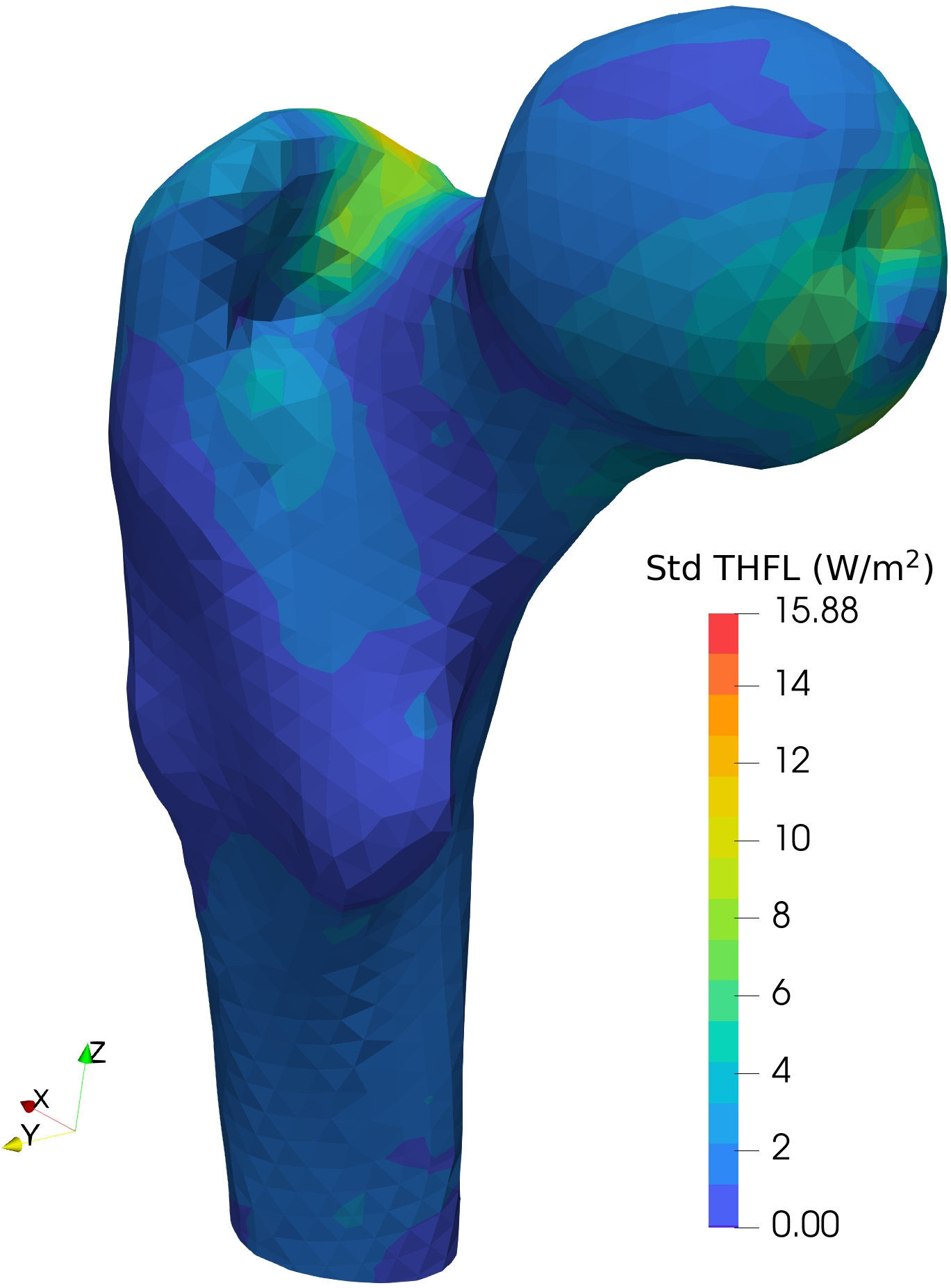}
		\caption{iso-ortho-scl}
	\end{subfigure}
% 	\begin{subfigure}{0.32\textwidth}
% 		\centering
% 		\includegraphics[width=\widththreeD\linewidth]{3D_scl_ortho_ortho_CCX_sideflux_1e5_var_HFL}
% 		\caption{ortho-ortho-scl}
% 	\end{subfigure}
	\begin{subfigure}{0.32\textwidth}
		\centering
		\includegraphics[width=\widththreeD\linewidth]{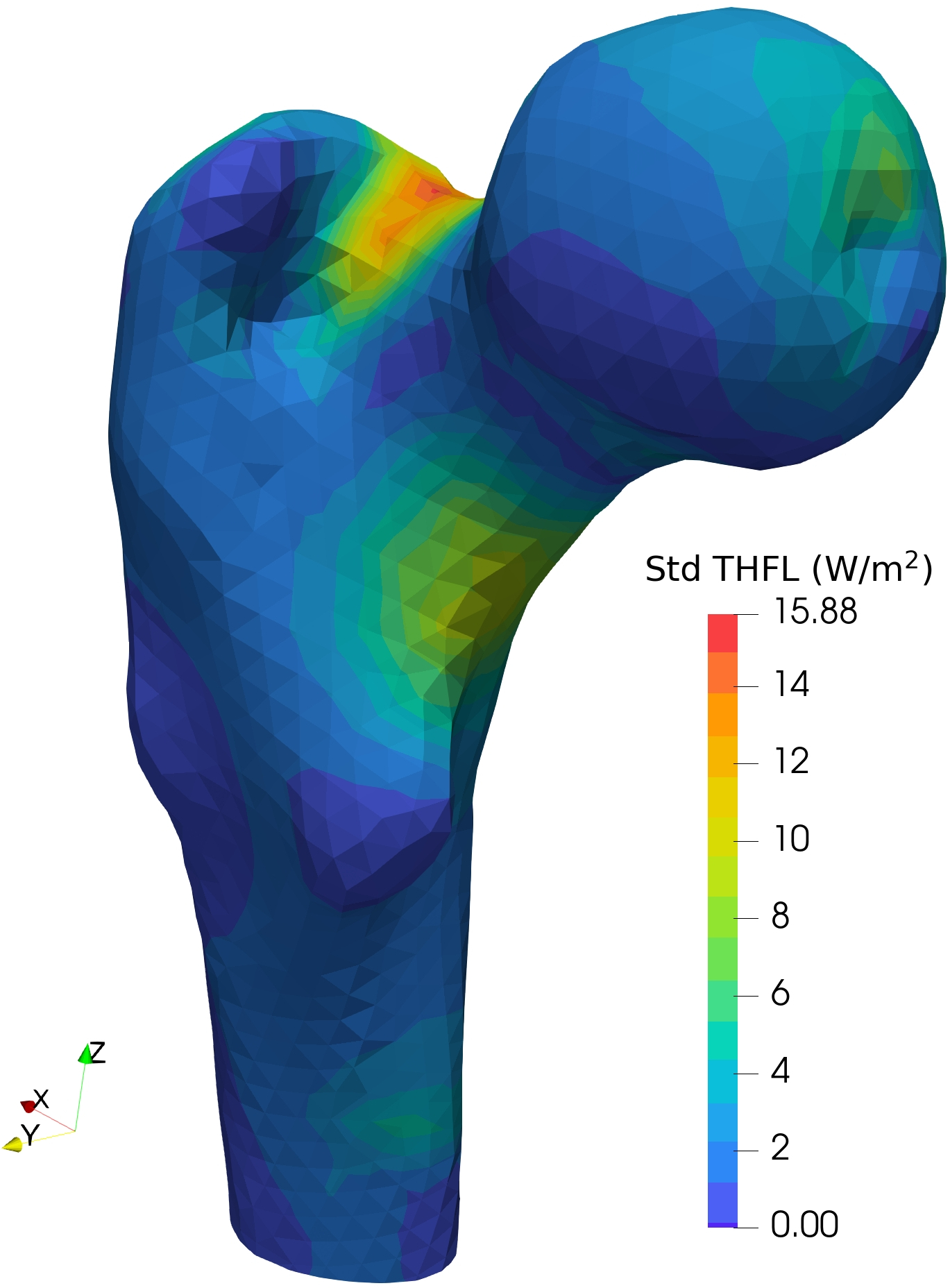}
		\caption{\acrshort{ortho}-\acrshort{ortho}-\acrshort{dir}}
	\end{subfigure}
% 	\begin{subfigure}{0.32\textwidth}
% 		\centering
% 		\includegraphics[width=\widththreeD\linewidth]{3D_rotscl_ortho_ortho_CCX_sideflux_1e5_var_HFL}
% 		\caption{ortho-ortho-rotscl}
% 	\end{subfigure}
% 	\begin{subfigure}{0.32\textwidth}
% 		\centering
% 		\includegraphics[width=\widthtwoD\linewidth]{3D_rotscl_iso_ortho_CCX_sideflux_1e5_var_HFL}
% 		\caption{iso-ortho-rotscl}
% 	\end{subfigure}
	\caption{\acrshort{mc} standard deviation estimate of total heat flux (Euclidean norm)}
	\label{Fig:stdTHFL3D}
\end{figure}

%-----------------------------------------------------------mean var dir-------------------------------------------------------------

%% file: conclusion.tex
% !TEX root = ../paper_tensors.tex
% !TEX encoding = UTF-8 Unicode

\section{Conclusion}\label{Conclusion}
Here the task of modelling and representing random symmetric and positive definite (SPD)
tensors was addressed, and some desiderata for the modelling and representation
were formulated.  Namely, the representation has to be
SPD even under numerical approximation, and has to be such that 
one is able to control the material or spatial invariance resp.\ symmetry, both for each possible 
realisation and for the the expected value.  Regarding the mean or expected value,
it turned out that the usual expectation, which is tied to the structure of a
linear space, may not be appropriate on the open convex cone of SPD tensors.
The more general Fr\'{e}chet mean is based on distance measurements, possibly different from
the usual Euclidean distance, on the set of SPD tensors.  Here some desiderata were
formulated for the distance measure, which take into account the physical purpose and
uncertainty involved in stochastic models of such physically symmetric and positive
definite tensors.

Therefore in Appendix~\ref{App:metrics-SPD} the type of metric that one can consider
for symmetric positive tensors and specifically for the set of SPD 
matrices to compute the correspondin Fr\'{e}chet mean is discussed. 
As the Euclidean mean suffers from the so-called swelling effect, we suggest 
alternatives such as the geometric, log-Euclidean, and scale-rotation separated mean. 
%The inappropriateness of the Euclidean mean for random SPD matrices with random orientations 
%in a 2D case is further explored in more detail in Appendix~\ref{App:derivation}. 
%A special case appears when the 
%random tensor is modelled with random orientations which have isotropic symmetry in the 
%mean and any lower-order material symmetry in the realisation. In this scenario, it is 
%found that the Euclidean mean of random SPD tensor does not suffer from any distortion.

It was also proposed to follow a widely used technique
to model and represent the logarithm of tensor, such that the symmetric positive definite
attribute is automatically satisfied when taking the exponential in the end; the
invariance properties automatically translate to the logarithm.
To expose the ideas in the simplest possible setting, we consider specifically second-order
SPD tensors, where the symmetry classes are quite simple.
%and only fields spatially constant with each realisation, i.e.\ SPD tensor-valued
%random variables.

%We model and generate second-order SPD material tensors with uncertainties in a spatially
%constant random tensor field sense. 
More particularly, the focus is on modelling uncertainties in such a way as to have a
fine control independently over the strength/scaling and directional attributes of the
tensor---given that the material symmetry is fixed. 
As we separate the strength and orientations of the random SPD matrices by spectral 
decomposition, the task boils down to modelling the random scaling parameters of the
logarithm of a random diagonal matrix, as well as the random orthogonal 
matrix of eigenvectors.  Such orthogonal matrices are members of the Lie group $\mrm{SO}(d)$,
and the well-known exponential map from the corresponding Lie algebra $\F{so}(d)$ of
skew-symmetric matrices was used to carry random ensembles in the free vector
space $\F{so}(d)$ onto such ensembles on $\mrm{SO}(d)$.  This results
in the random orientation being modelled by symmetric circular/spherical distributions.
Furthermore, a scenario of fluctuating the material symmetry around the 
Fr\'{e}chet mean is also studied, 
i.e.\ the mean tensor belongs to a higher order of material or spatial symmetry, whereas each 
realisation belongs to a lower order of material symmetry. 

As an application in this study, the random thermal conductivity tensor of a 
steady-state heat conduction on a 2D and a 3D model of a proximal femur bone
was investigated, a highly anisotropic material.
%The scaling values of the thermal conductivity tensor are modelled as positive 
%log-normal random variables, whereas the directional uncertainty of the tensor is 
%accounted by von the von Mises-Fisher spherical distribution in 3D, and the von Mises circular 
%distribution in 2D. 
Three different modelling scenarios are considered: random scalings only, 
random scalings with fluctuating symmetry, and 
random orientations only. 
The material uncertainties are then propagated to the output responses such 
as \acrfull{nt}, \acrfull{thfl} 
and \acrfull{nhfl} using the \acrfull{mc}. Subsequently, statistics such as mean and standard deviation 
estimates of NT and THFL, and circular mean and circular standard deviation estimates 
of NHFL are determined.

From the numerical results of both 2D and 3D models, it is clear that NT is most 
stochastically sensitive to the model with random scalings only, 
and least sensitive to random orientation only. 
On the other hand, the randomness in the model with random scaling only has almost zero stochastic 
influence on THFL---determined by the Euclidean norm---but the impact of the other two 
random models is evident.  One may also observe a similar stochastic sensitivity of the stochastic 
models on NHFL, as well as on THFL.
Therefore, the distinctive impact of the three stochastic models on the output 
response signifies the importance of incorporating different material 
uncertainties---scaling, orientations and material symmetry---independently into 
the constitutive model.

In total, we have proposed a modelling scenario for tensors defining positive-definite
material properties which takes into account possible material resp.\ spatial 
invariances or symmetries.
Our approach allows a fine control of the scaling randomness, and independently of the orientation 
randomness, further separating it into the invariance properties of the mean and
those of the random fluctuations.

%% file: abstract-metrics.tex
% !TEX root = ../paper_tensors.tex
% !TEX encoding = UTF-8 Unicode

\section{Metrics on symmetric positive definite matrices}\label{App:metrics-SPD}
As was already noted in \fsec{Introduction}, positive definite tensors 
$\PosSymSecRankSpace \subset \SymSecRankSpace$ live geometrically on an open convex cone
in the vector space of symmetric tensors $\SymSecRankSpace$.  It is a differentiable
manifold, but it is also part of the ambient space, the Euclidean resp.\  Frobenius or
Hilbert-Schmidt norm.
It has been found that this Euclidean metric on $\SymSecRankSpace$, 
%resulting from the usual Frobenius resp.\ Hilbert-Schmidt norm resp.\ inner product
\begin{linenomath}
\begin{equation}  \label{Eq:Frob-norm}
\nd{\genSecTensor}^2_F := \ip{\genSecTensor}{\genSecTensor}_F := \tr(\genSecTensor^T\genSecTensor),
\end{equation}
\end{linenomath}
namely
\begin{linenomath}
\begin{equation}  \label{Eq:Euclid-dist}
\vartheta_F(\genSecTensor_1,\genSecTensor_2) :=\nd{\genSecTensor_1-\genSecTensor_2}_F,
\end{equation}
\end{linenomath}
when used on the open convex cone $\PosSymSecRankSpace \subset \SymSecRankSpace$
\cite{PennecFillardAyache2004, AndoLiMathias2004, Moakher2005, ArsignyFillardPennecEtAl2006, 
ArsignyFillardPennecEtAl2007, drydenEtal2009, DrydenEtal2010, FujiiSeo2015},
leads to the undesirable \emph{swelling, fattening, and shrinking effects} in interpolation
resp.\ averaging \cite{schwartzman_random_nodate, JungSchwartzmanGroisser2015, 
schwartzman_lognormal_2016, GroisserJungSchwartzman2017, groissJungSchwman2017, FeragenFuster2017}.
In addition, $\PosSymSecRankSpace \subset \SymSecRankSpace$ being an open subset
means that $\PosSymSecRankSpace$ is not metrically complete with the Euclidean metric
$\vartheta_F$, which in turn means that Cauchy sequences in $\PosSymSecRankSpace$ in the 
$\vartheta_F$ metric may have limits which are not in $\PosSymSecRankSpace$, 
i.e.\ they could be singular, which may present a problem in approximation algorithms.

\subsection{Metrics}  \label{Ssec:spd-metrics}
%
%As was already mentioned in \fsec{Introduction}, it has been found that the
%Euclidean metric on $\SymSecRankSpace$ resulting from the usual Frobenius resp.\
%Hilbert-Schmidt norm resp.\ inner product
%\begin{equation}  \label{Eq:Frob-norm}
%\nd{\genSecTensor}^2_F := \ip{\genSecTensor}{\genSecTensor}_F := \tr(\genSecTensor^T\genSecTensor),
%\end{equation}
%namely
%\begin{equation}  \label{Eq:Euclid-dist}
%\vartheta_F(\genSecTensor_1,\genSecTensor_2) :=\nd{\genSecTensor_1-\genSecTensor_2}_F,
%\end{equation}
%when used on the open convex cone $\PosSymSecRankSpace \subset \SymSecRankSpace$
%\cite{PennecFillardAyache2004, AndoLiMathias2004, Moakher2005, ArsignyFillardPennecEtAl2006, 
%ArsignyFillardPennecEtAl2007, drydenEtal2009, DrydenEtal2010, FujiiSeo2015},
%leads to the undesirable \emph{swelling, fattening, and shrinking effects} in interpolation
%resp.\ averaging \cite{schwartzman_random_nodate, JungSchwartzmanGroisser2015, 
%schwartzman_lognormal_2016, GroisserJungSchwartzman2017, groissJungSchwman2017}.
%In addition, $\PosSymSecRankSpace \subset \SymSecRankSpace$ being an open subset
%means that $\PosSymSecRankSpace$ is not metrically complete with the Euclidean metric
%$\vartheta_F$, which in turn means that Cauchy sequences in $\PosSymSecRankSpace$ in the 
%$\vartheta_F$ metric may have limits which are not in $\PosSymSecRankSpace$, 
%i.e.\ they could be singular, which may present a problem in approximation algorithms.

What metric to take instead of the Euclidean one connected with the arithmetic mean
is not so clear.  For an overview of various options one may consult 
\cite{NielsenBhatia2013}, see also the discussion in  \cite{GroisserJungSchwartzman2017, 
PennecSommerFletcher2020}.  Most proposals turn $\PosSymSecRankSpace$, which as an open
subset of the linear vector space $\SymSecRankSpace$ carries a natural structure of a 
differential manifold, into a Riemannian manifold in different ways with different metrics.
This allows to measure the lengths of continuous paths in $\PosSymSecRankSpace$,
thus leading to shortest paths or \emph{geodesics}, on which one can naturally 
define interpolants resp.\ averages.
A good introduction to this very active subject \cite{ThanwerdasPennec2019, ThanwerdasPennec2019-2, 
ThanwerdasPennec2021} relating to metrics on $\PosSymSecRankSpace$ may be found
in the collection \cite{PennecSommerFletcher2020}, 
see also \cite{NielsenBhatia2013, FeragenFuster2017, groissJungSchwman2023}.

As will be seen, the open convex cone $\PosSymSecRankSpace\subset \SymSecRankSpace$ may be considered
as a metric space in ways different from just being a subset of the vector space 
$\SymSecRankSpace$.  In addition to the desiderata listed in \fsec{Sec:Problem} for
a metric, and connected obviously to the first requirement,
one may want that $\vartheta_D(\vek{I},\alpha\genSecTensor) \to \infty$ as
$\alpha\to\infty$ for any $\genSecTensor\in \PosSymSecRankSpace$, 
but the same should also hold as $\alpha\to 0$, in fact singular matrices in $\SymSecRankSpace$
should all be infinitely far away from the identity, i.e.\ all of the boundary  
$\partial\,\PosSymSecRankSpace \subset \SymSecRankSpace$ should be at an infinite
distance from the identity.

\subsection{Affine-invariant distance} \label{Asec:aff-inv}
The best known alternative to the Euclidean distance starts from the observation that
as an open convex cone in a Euclidean space, $\PosSymSecRankSpace$ can be considered
as a Riemannian manifold with tangent space $\SymSecRankSpace$, with the Euclidean
resp.\ Frobenius or Hilbert-Schmidt inner product \feq{Eq:Frob-norm} on the tangent space at each point.
Thus one may define geodesics, and their length gives a distance metric $\vartheta_G$ on 
$\PosSymSecRankSpace$, often termed the  geometric or affine-invariant distance
\begin{linenomath}
\begin{equation}  \label{Eq:aff-inv-dist}
\vartheta_G(\genSecTensor_1,\genSecTensor_2):= \nd{\log(\genSecTensor_1^{-1/2} \genSecTensor_2
     \genSecTensor_1^{-1/2})}_F,
\end{equation}
\end{linenomath}
which in turn leads to the so-called geometric or affine-invariant mean 
\cite{AndoLiMathias2004, Moakher2005, NielsenBhatia2013, PennecSommerFletcher2020}.  
Observe that $\genSecTensor_1^{-1/2} \genSecTensor_2
\genSecTensor_1^{-1/2}\in \PosSymSecRankSpace$, and the existence of a unique logarithm
is guaranteed by the spectral calculus  \cite{segalKunze78, yosida-fa-1980}.
One may remark that $\vartheta_G(\genSecTensor_1,\genSecTensor_2)
= \sum_k (\log\lambda_k)^2$, where $\lambda_k$ are the eigenvalues of the matrix pencil
$(\genSecTensor_1,\genSecTensor_2)$, i.e.\ the eigenvalues of the generalised eigenvalue problem
$\genSecTensor_1 \vek{s}_k = \lambda_k \genSecTensor_2 \vek{s}_k$.  This distance satisfies 
all the desired properties  \cite{AndoLiMathias2004, Moakher2005, NielsenBhatia2013}
listed in \fsec{Sec:Problem}, as well as that the boundary
$\partial\,\PosSymSecRankSpace$ is infinitely far away from the identity---this is
here a consequence of working  with the logarithm.
But other choices for the metric are possible  \cite{FujiiSeo2015, GroisserJungSchwartzman2017}.  
Observe that when $\genSecTensor_1$
and $\genSecTensor_2$ commute, $\vartheta_G(\genSecTensor_1,\genSecTensor_2)= 
\nd{\log \genSecTensor_2 - \log \genSecTensor_1 }_F$, and many of the problems
of defining a suitable metric are connected with the difficulties arising from
non-commuting tensors.

\subsection{Lie-group distance for orthogonal matrices} \label{Asec:Lie-grp-dist}
One may recall that for a connected compact Lie group like the special orthogonal
group $\mrm{SO}(d)$ which is used in our approach, one may define a metric invariant 
under the action of the group on itself (e.g.\ \cite{Moakher2002}),
which is hence a very natural way to obtain a metric.  It takes advantage of the correspondence
between the Lie group $\mrm{SO}(d)$ and its Lie algebra $\F{so}(d)$
of skew-symmetric matrices given by the exponential map resp.\ the logarithm, and is given by
\begin{linenomath}
\begin{equation}  \label{Eq:Lie-grp-dist}
\vartheta_R(\vek{Q}_1,\vek{Q}_2) := \nd{\log(\vek{Q}_1^T \vek{Q}_2)}_F
\end{equation}
\end{linenomath}
for $\vek{Q}_1, \vek{Q}_2 \in \mrm{SO}(d)$.  
Observe that $\vek{Q}_1^T \vek{Q}_2\in\mrm{SO}(d)$
and hence $\log(\vek{Q}_1^T \vek{Q}_2)\in\F{so}(d)$, and if $\vek{Q}_1$
and $\vek{Q}_2$ commute, $\log(\vek{Q}_1^T \vek{Q}_2) = \log\vek{Q}_2 - \log\vek{Q}_1$;
compare this with $\vartheta_G$ \feq{Eq:aff-inv-dist} above in the commutative case.
In any case, straight lines through the origin in the Lie algebra of skew matrices
$\F{so}(d)$ are geodesics with the Frobenius distance, and these are mapped
by the exponential map---the inverse of the logarithm---onto geodesics in 
$\mrm{SO}(d)$.  It may be also noted that the exponential map and its
inverse, the logarithm, establish a correspondence between the Lie algebra
$\F{so}(d)$---which is one-to-one in the ball of $\skewsym\in\F{so}(d)$ defined
as in \feq{Eq:def-skew-mat}
by Euler-vectors in the ball $\{ \axis\in\mathbb{R}^3 \mid \|\axis\| \le \uppi \}$
with antipodal points identified---and 
the Lie group $\mrm{SO}(d)$, see also \cite{cardoso_2010}.

\subsection{Log-Euclidean distance} \label{Asec:log-Euclid}
Unfortunately, when considering these results for Lie groups, one has to note that
$\PosSymSecRankSpace$ is not a group under matrix multiplication,
although---as was already noted---it is stable or invariant under inversion.
But taking the situation for $\mrm{SO}(d)$  as a cue, it has been found 
\cite{ArsignyFillardPennecEtAl2007} that there is a way of making $\PosSymSecRankSpace$ 
into a Lie group with a new and commutative definition of multiplication:
\begin{linenomath}
\begin{equation}  \label{Eq:SymPosDasLieG}
  \forall \genSecTensor_1,\genSecTensor_2\in \PosSymSecRankSpace:\;
  \genSecTensor_1 \boxtimes\genSecTensor_2 := \exp\left(\log \genSecTensor_1 +
  \log \genSecTensor_2 \right).
\end{equation}
\end{linenomath}
This makes $\PosSymSecRankSpace$ into a commutative Lie group 
\cite{ArsignyFillardPennecEtAl2007} with the multiplication `$\boxtimes$'
---indeed, the inverse of $\genSecTensor\in \PosSymSecRankSpace$ w.r.t.\ this multiplication
is the normal matrix inverse $\genSecTensor^{-1}\in \PosSymSecRankSpace$, and the neutral
element is the normal identity matrix $\vek{I}\in \PosSymSecRankSpace$ ---and 
whenever $\genSecTensor_1$ and $\genSecTensor_2$ commute, this new product indeed 
coincides with the normal matrix product.  Further it can be shown 
\cite{ArsignyFillardPennecEtAl2007} that the Lie algebra of this
Lie Group $(\PosSymSecRankSpace,\boxtimes)$ is the vector space of symmetric
matrices $\SymSecRankSpace$, and the exponential map in the sense of Lie algebras
is indeed the matrix exponential.  In fact, the exponential map and its
inverse, the logarithm, establish a one-to-one correspondence between the Lie algebra
$\SymSecRankSpace$ and the Lie group $\PosSymSecRankSpace$, actually
a diffeomorphism and an algebraic isomorphism between the groups $(\PosSymSecRankSpace,
\boxtimes)$ (with the new multiplication) and the additive group $(\SymSecRankSpace,+)$
of the vector space of symmetric matrices.  Thus $\SymSecRankSpace$ may be
seen as a global chart for $\PosSymSecRankSpace$.
This means in particular that straight lines through the origin---geodesics with the 
Euclidean or Frobenius distance in $\SymSecRankSpace$ ---are mapped into geodesics in  
$\PosSymSecRankSpace$ with the logarithmic or log-Euclidean metric
\begin{linenomath}
\begin{equation}  \label{Eq:log-Euclid-dist}
\vartheta_L(\genSecTensor_1,\genSecTensor_2):= \nd{\log \genSecTensor_1 -
  \log \genSecTensor_2 }_F,
\end{equation}
\end{linenomath}
which also satisfies all desiderata listed in 
\fsec{Sec:Problem}, and has the boundary $\partial\,\PosSymSecRankSpace$ 
at an infinite distance from the identity.
Among others, this metric is also investigated in \cite{FujiiSeo2015} (the authors
call $\vartheta_L$ the \emph{chaotic geometric distance}) and
 \cite{schwartzman_lognormal_2016}. 
Notice also the equality with $\vartheta_G$ \feq{Eq:aff-inv-dist} above for the commuting case---if 
$\genSecTensor_1$ and $\genSecTensor_2$ commute, $\vartheta_L(\genSecTensor_1,
\genSecTensor_2)= \nd{\log (\genSecTensor_1 \genSecTensor_2^{-1}}_F =
\vartheta_G(\genSecTensor_1,\genSecTensor_2)$ ---and also the similarity  with the 
Riemannian metric $\vartheta_R$ \feq{Eq:Lie-grp-dist} on $\mrm{SO}(d)$.
This means that the two metrics $\vartheta_G$ \feq{Eq:aff-inv-dist} and 
$\vartheta_L$ \feq{Eq:log-Euclid-dist} agree on commuting subsets of $\PosSymSecRankSpace$.
But such a simple Euclidean structure on the Lie algebra mixes the strength and orientational
information and disregards the effects of non-commutativity.  
In our view the non-commutative aspect is part of what
is the orientational information in the eigenvectors of a SPD matrix, and for that
reason we decide to treat strength and orientation aspects separately.

\subsection{Scaling-rotation distance} \label{Asec:scale-rot}
Thus another possibility arises from the observation that two tensors 
$\genSecTensor_1,\genSecTensor_2\in \PosSymSecRankSpace$ commute only if they
have the same invariant subspaces.  In that case it can be arranged that the 
eigenvector matrices according to the spectral decomposition \feq{Eq:SpecDecSPD}, 
$\genSecTensor_1=\vek{Q}_1\vek{\Lambda}_1\vek{Q}_1^T$ and
$\genSecTensor_2=\vek{Q}_2\vek{\Lambda}_2\vek{Q}_2^T$, are equal, i.e.\ 
$\vek{Q}_1 = \vek{Q}_2$.     In this instance of two commuting tensors, one has
\begin{linenomath}
\begin{multline}  \label{Eq:comm-dist}
\vartheta_L(\genSecTensor_1,\genSecTensor_2)= \vartheta_G(\genSecTensor_1,\genSecTensor_2) 
=\vartheta_L(\vek{\Lambda}_1,\vek{\Lambda}_2) \\ = \nd{\log \vek{\Lambda}_1 -
  \log \vek{\Lambda}_2 }_F = \nd{\log (\vek{\Lambda}_1 \vek{\Lambda}_2^{-1})}_F =
\sqrt{\sum_{k} (\log(\lambda_{k,1}/\lambda_{k,2}))^2} .
\end{multline}
\end{linenomath}

As discussed in \fsec{Ssec:ensuring-spd},
%Let \SymSecRankSpace$ be the vector subspace
%of diagonal matrices, and $\mrm{Diag}^+(d)\subset \PosSymSecRankSpace$ the
%positive definite diagonal matrices, i.e.\ having positive diagonal terms.
%Not only is this again an open convex cone in $\mrm{Diag}(d)$, but obviously also
%an Abelian or commutative Lie group under matrix multiplication.
the spectral decomposition  \feq{Eq:SpecDecSPD} induces the representation
of \feq{Eq:repr-prod-Lie-2} on a product of Lie groups,
namely on: $\mrm{Diag}^+(d) \times \mrm{SO}(d)$, 
which is itself a Lie group; this is the right half of \feq{Eq:repr-prod-Lie-exp}.
%\begin{equation}  \label{Eq:repr-prod-Lie}
%  \mrm{Diag}^+(d) \times \mrm{SO}(d) \ni (\vek{\Lambda},\vek{Q})
%    \mapsto \genSecTensor=\vek{Q}\vek{\Lambda}\vek{Q}^T \in \PosSymSecRankSpace .
%\end{equation}
%One should caution \cite{schwartzman_random_nodate, JungSchwartzmanGroisser2015, 
%schwartzman_lognormal_2016, GroisserJungSchwartzman2017, groissJungSchwman2017}
%that this representation map \feq{Eq:repr-prod-Lie} is onto, but not one-to-one,
%as can be seen when multiple eigenvalues occur.  Observe that the product on the
%left of  \feq{Eq:repr-prod-Lie} is a product of Lie groups, and thus itself a Lie group.
From this and  \feq{Eq:comm-dist} comes the idea 
\cite{JungSchwartzmanGroisser2015, GroisserJungSchwartzman2017}
to measure the distance
on $\PosSymSecRankSpace$ by measuring separately the distance between the 
eigenspaces (the rotation matrices $\vek{Q}_k$), and the distance between 
the eigenvalues (the scaling matrices $\vek{\Lambda}_k$), and define a distance
function by combining $\vartheta_L$ \feq{Eq:log-Euclid-dist} on the commutative factor
$\mrm{Diag}^+(d)$ with $\vartheta_R$ \feq{Eq:Lie-grp-dist} on $\mrm{SO}(d)$
to a product distance
\begin{linenomath}
\begin{equation}  \label{Eq:prod-Lie-dist}
\vartheta_{P}((\vek{\Lambda}_1,\vek{Q}_1),(\vek{\Lambda}_2,\vek{Q}_2))^2= 
  \vartheta_L(\vek{\Lambda}_1,\vek{\Lambda}_2)^2 + c\, \vartheta_R(\vek{Q}_1,\vek{Q}_2)^2
\end{equation}
\end{linenomath}
for some $c>0$ on the product Lie group $\mrm{Diag}^+(d) \times \mrm{SO}(d)$,
i.e.\ one measures the distance squared between
the scales and the orientations separately and adds them both up.

As was pointed in \fsec{Ssec:ensuring-spd}, the representation map in 
\feq{Eq:repr-prod-Lie-2} is not one-to-one, so this does not carry
directly over to a metric on $\PosSymSecRankSpace$ 
\cite{GroisserJungSchwartzman2017, FeragenFuster2017}.
Thus one sets, \rvsi{by minimising over all possible decompositions,}
\begin{linenomath}
\begin{equation}  \label{Eq:scal-dist1}
\rvsi{\widetilde{\vartheta}_{S}}((\genSecTensor_1,\genSecTensor_2)) := \min \{
  \vartheta_{P}((\vek{\Lambda}_1,\vek{Q}_1),(\vek{\Lambda}_2,\vek{Q}_2)) \mid
  \genSecTensor_1=\vek{Q}_1\vek{\Lambda}_1\vek{Q}_1^T \, , \,
    \genSecTensor_2=\vek{Q}_2\vek{\Lambda}_2\vek{Q}_2^T \},
\end{equation}
\end{linenomath}
where $(\vek{\Lambda}_1,\vek{Q}_1),(\vek{\Lambda}_2,\vek{Q}_2) \in \mrm{Diag}^+(d) \times \mrm{SO}(d)$.  
\rvsi{Unfortunately, $\tilde{\vartheta}_{S}$, although it satisfies all our desiderata listed in 
\fsec{Sec:Problem}, may fail the triangle inequality---cf.\ the discussion in 
\cite{FeragenFuster2017} ---and is thus not necessarily a distance resp.\ a metric.  
Luckily, this can be repaired by setting \cite{FeragenFuster2017}}
\rvsi{
\begin{linenomath}
\begin{multline}  \label{Eq:scal-dist}
  \vartheta_{S}((\genSecTensor_a,\genSecTensor_b)) := \inf \{ \sum_{j=1}^{n-1}
  \widetilde{\vartheta}_{S}(\genSecTensor_j,\genSecTensor_{j+1}) \mid n \ge 2,
  \genSecTensor_a = \genSecTensor_1, \,     \genSecTensor_b=\genSecTensor_n, \\ 
  \genSecTensor_j \in \PosSymSecRankSpace,  j=1,\dots,n \}.
\end{multline}
\end{linenomath}
}
\rvsi{This expression, $\vartheta_{S}$ in \feq{Eq:scal-dist}, now is a metric 
and satisfies all our desiderata.  It is termed the scaling resp.\ scaling-rotation distance
and has the boundary $\partial\,\PosSymSecRankSpace$ infinitely far away from the identity,
as this already holds for the product distance $\vartheta_{P}$ in \feq{Eq:prod-Lie-dist}.}

As we propose to rather work with $\vek{H}=\log(\genSecTensor) \in\SymSecRankSpace$,
where $\log(\vek{Q}\vek{\Lambda}\vek{Q}^T) = \vek{Q}\log(\vek{\Lambda})\vek{Q}^T
=\vek{Q}\vek{Y}\vek{Q}^T$, it is important to note that
the metrics mentioned before can easily be computed on $\vek{H}=\log(\genSecTensor)$.  One has
\begin{linenomath}
\begin{align}  \label{Eq:metric_G}
  \vartheta_G(\genSecTensor_1,\genSecTensor_2) &= \nd{\log(\exp(-(1/2)\vek{Y}_1) \vek{Q}_1^T
  \vek{Q}_2 \exp(\vek{Y}_2) \vek{Q}_2^T \vek{Q}_1 \exp(-(1/2)\vek{Y}_1))}_F, \\  \label{Eq:metric_L}
  \vartheta_L(\genSecTensor_1,\genSecTensor_2) &= \nd{\vek{H}_1-\vek{H}_2}_F, \\  \label{Eq:metric_S}
  \vartheta_P(\genSecTensor_1,\genSecTensor_2) &= 
%    \sqrt{\nd{\vek{Y}_1-\vek{Y}_2}^2_F  +c\,\vartheta_R(\vek{Q}_1,\vek{Q}_2)^2} = 
    \sqrt{\nd{\vek{Y}_1-\vek{Y}_2}^2_F+c\,\nd{\log(\vek{Q}_1^T \vek{Q}_2)}_F^2}.
\end{align}
\end{linenomath}